
\ifx\shlhetal\undefinedcontrolsequence\let\shlhetal\relax\fi
\def\fmtname{AmS-TeX}

\def\fmtversion{2.1}
\catcode`\@=11
\ifx\amstexloaded@\relax\catcode`\@=\active
  \endinput\else\let\amstexloaded@\relax\fi
\newlinechar=`\^^J
\def\W@{\immediate\write\sixt@@n}
\def\CR@{\W@{^^J\fmtname - Version \fmtversion^^J%
COPYRIGHT 1985, 1990, 1991 - AMERICAN MATHEMATICAL SOCIETY^^J%
Use of this macro package is not restricted provided^^J%
each use is acknowledged upon publication.^^J}}
\CR@ \everyjob{\CR@}
\message{Loading definitions for}
\message{misc utility macros,}
\toksdef\toks@@=2
\long\def\rightappend@#1\to#2{\toks@{\\{#1}}\toks@@
 =\expandafter{#2}\xdef#2{\the\toks@@\the\toks@}\toks@{}\toks@@{}}
\def\alloclist@{}
\newif\ifalloc@
\def\showallocations{{\def\\{\immediate\write\m@ne}\alloclist@}\alloc@true}
\def\alloc@#1#2#3#4#5{\global\advance\count1#1by\@ne
 \ch@ck#1#4#2\allocationnumber=\count1#1
 \global#3#5=\allocationnumber
 \edef\next@{\string#5=\string#2\the\allocationnumber}%
 \expandafter\rightappend@\next@\to\alloclist@}
\newcount\count@@
\newcount\count@@@
\def\FN@{\futurelet\next}
\def\DN@{\def\next@}
\def\DNii@{\def\nextii@}
\def\RIfM@{\relax\ifmmode}
\def\RIfMIfI@{\relax\ifmmode\ifinner}
\def\setboxz@h{\setbox\z@\hbox}
\def\wdz@{\wd\z@}
\def\boxz@{\box\z@}
\def\setbox@ne{\setbox\@ne}
\def\wd@ne{\wd\@ne}
\def\iterate{\body\expandafter\iterate\else\fi}
\def\err@#1{\errmessage{AmS-TeX error: #1}}
\newhelp\defaulthelp@{Sorry, I already gave what help I could...^^J
Maybe you should try asking a human?^^J
An error might have occurred before I noticed any problems.^^J
``If all else fails, read the instructions.''}
\def\Err@{\errhelp\defaulthelp@\err@}
\def\eat@#1{}
\def\in@#1#2{\def\in@@##1#1##2##3\in@@{\ifx\in@##2\in@false\else\in@true\fi}%
 \in@@#2#1\in@\in@@}
\newif\ifin@
\def\space@.{\futurelet\space@\relax}
\space@. %
\newhelp\athelp@
{Only certain combinations beginning with @ make sense to me.^^J
Perhaps you wanted \string\@\space for a printed @?^^J
I've ignored the character or group after @.}
{\catcode`\~=\active 
 \lccode`\~=`\@ \lowercase{\gdef~{\FN@\at@}}}
\def\at@{\let\next@\at@@
 \ifcat\noexpand\next a\else\ifcat\noexpand\next0\else
 \ifcat\noexpand\next\relax\else
   \let\next\at@@@\fi\fi\fi
 \next@}
\def\at@@#1{\expandafter
 \ifx\csname\space @\string#1\endcsname\relax
  \expandafter\at@@@ \else
  \csname\space @\string#1\expandafter\endcsname\fi}
\def\at@@@#1{\errhelp\athelp@ \err@{\Invalid@@ @}}
\def\atdef@#1{\expandafter\def\csname\space @\string#1\endcsname}
\newhelp\defahelp@{If you typed \string\define\space cs instead of
\string\define\string\cs\space^^J
I've substituted an inaccessible control sequence so that your^^J
definition will be completed without mixing me up too badly.^^J
If you typed \string\define{\string\cs} the inaccessible control sequence^^J
was defined to be \string\cs, and the rest of your^^J
definition appears as input.}
\newhelp\defbhelp@{I've ignored your definition, because it might^^J
conflict with other uses that are important to me.}
\def\define{\FN@\define@}
\def\define@{\ifcat\noexpand\next\relax
 \expandafter\define@@\else\errhelp\defahelp@                               
 \err@{\string\define\space must be followed by a control
 sequence}\expandafter\def\expandafter\nextii@\fi}                          
\def\undefined@@@@@@@@@@{}
\def\preloaded@@@@@@@@@@{}
\def\next@@@@@@@@@@{}
\def\define@@#1{\ifx#1\relax\errhelp\defbhelp@                              
 \err@{\string#1\space is already defined}\DN@{\DNii@}\else
 \expandafter\ifx\csname\expandafter\eat@\string                            
 #1@@@@@@@@@@\endcsname\undefined@@@@@@@@@@\errhelp\defbhelp@
 \err@{\string#1\space can't be defined}\DN@{\DNii@}\else
 \expandafter\ifx\csname\expandafter\eat@\string#1\endcsname\relax          
 \global\let#1\undefined\DN@{\def#1}\else\errhelp\defbhelp@
 \err@{\string#1\space is already defined}\DN@{\DNii@}\fi
 \fi\fi\next@}

\def\predefine#1#2{\let#1#2}
\def\undefine#1{\let#1\undefined}
\message{page layout,}
\newdimen\captionwidth@
\captionwidth@\hsize
\advance\captionwidth@-1.5in
\def\pagewidth#1{\hsize#1\relax
 \captionwidth@\hsize\advance\captionwidth@-1.5in}
\def\pageheight#1{\vsize#1\relax}
\def\hcorrection#1{\advance\hoffset#1\relax}
\def\vcorrection#1{\advance\voffset#1\relax}
\message{accents/punctuation,}

\let\graveaccent\`
\let\acuteaccent\'
\let\tildeaccent\~
\let\hataccent\^
\let\underscore\_
\let\B\=
\let\D\.
\let\ic@\/
\def\/{\unskip\ic@}
\def\textfonti{\the\textfont\@ne}
\def\t#1#2{{\edef\next@{\the\font}\textfonti\accent"7F \next@#1#2}}
\def~{\unskip\nobreak\ \ignorespaces}
\def\.{.\spacefactor\@m}
\atdef@;{\leavevmode\null;}
\atdef@:{\leavevmode\null:}
\atdef@?{\leavevmode\null?}
\edef\@{\string @}
\def\relaxnext@{\let\next\relax}
\atdef@-{\relaxnext@\leavevmode
 \DN@{\ifx\next-\DN@-{\FN@\nextii@}\else
  \DN@{\leavevmode\hbox{-}}\fi\next@}%
 \DNii@{\ifx\next-\DN@-{\leavevmode\hbox{---}}\else
  \DN@{\leavevmode\hbox{--}}\fi\next@}%
 \FN@\next@}
\def\srdr@{\kern.16667em}
\def\drsr@{\kern.02778em}
\def\sldl@{\drsr@}
\def\dlsl@{\srdr@}
\atdef@"{\unskip\relaxnext@
 \DN@{\ifx\next\space@\DN@. {\FN@\nextii@}\else
  \DN@.{\FN@\nextii@}\fi\next@.}%
 \DNii@{\ifx\next`\DN@`{\FN@\nextiii@}\else
  \ifx\next\lq\DN@\lq{\FN@\nextiii@}\else
  \DN@####1{\FN@\nextiv@}\fi\fi\next@}%
 \def\nextiii@{\ifx\next`\DN@`{\sldl@``}\else\ifx\next\lq
  \DN@\lq{\sldl@``}\else\DN@{\dlsl@`}\fi\fi\next@}%
 \def\nextiv@{\ifx\next'\DN@'{\srdr@''}\else
  \ifx\next\rq\DN@\rq{\srdr@''}\else\DN@{\drsr@'}\fi\fi\next@}%
 \FN@\next@}

\def\textfontii{\the\textfont\tw@}
\def\lbrace@{\delimiter"4266308 }
\def\rbrace@{\delimiter"5267309 }
\def\{{\RIfM@\lbrace@\else{\textfontii f}\spacefactor\@m\fi}
\def\}{\RIfM@\rbrace@\else
 \let\@sf\empty\ifhmode\edef\@sf{\spacefactor\the\spacefactor}\fi
 {\textfontii g}\@sf\relax\fi}
\let\lbrace\{
\let\rbrace\}
\def\AmSTeX{{\textfontii A\kern-.1667em%
  \lower.5ex\hbox{M}\kern-.125emS}-\TeX}
\message{line and page breaks,}
\def\vmodeerr@#1{\Err@{\string#1\space not allowed between paragraphs}}
\def\mathmodeerr@#1{\Err@{\string#1\space not allowed in math mode}}
\def\linebreak{\RIfM@\mathmodeerr@\linebreak\else
 \ifhmode\unskip\unkern\break\else\vmodeerr@\linebreak\fi\fi}

\newskip\saveskip@
\def\allowlinebreak{\RIfM@\mathmodeerr@\allowlinebreak\else
 \ifhmode\saveskip@\lastskip\unskip
 \allowbreak\ifdim\saveskip@>\z@\hskip\saveskip@\fi
 \else\vmodeerr@\allowlinebreak\fi\fi}
\def\nolinebreak{\RIfM@\mathmodeerr@\nolinebreak\else
 \ifhmode\saveskip@\lastskip\unskip
 \nobreak\ifdim\saveskip@>\z@\hskip\saveskip@\fi
 \else\vmodeerr@\nolinebreak\fi\fi}
\def\newline{\relaxnext@
 \DN@{\RIfM@\expandafter\mathmodeerr@\expandafter\newline\else
  \ifhmode\ifx\next\par\else
  \expandafter\unskip\expandafter\null\expandafter\hfill\expandafter\break\fi
  \else
  \expandafter\vmodeerr@\expandafter\newline\fi\fi}%
 \FN@\next@}
\def\dmatherr@#1{\Err@{\string#1\space not allowed in display math mode}}
\def\nondmatherr@#1{\Err@{\string#1\space not allowed in non-display math
 mode}}
\def\onlydmatherr@#1{\Err@{\string#1\space allowed only in display math mode}}
\def\nonmatherr@#1{\Err@{\string#1\space allowed only in math mode}}
\def\mathbreak{\RIfMIfI@\break\else
 \dmatherr@\mathbreak\fi\else\nonmatherr@\mathbreak\fi}
\def\nomathbreak{\RIfMIfI@\nobreak\else
 \dmatherr@\nomathbreak\fi\else\nonmatherr@\nomathbreak\fi}
\def\allowmathbreak{\RIfMIfI@\allowbreak\else
 \dmatherr@\allowmathbreak\fi\else\nonmatherr@\allowmathbreak\fi}
\def\pagebreak{\RIfM@
 \ifinner\nondmatherr@\pagebreak\else\postdisplaypenalty-\@M\fi
 \else\ifvmode\removelastskip\break\else\vadjust{\break}\fi\fi}
\def\nopagebreak{\RIfM@
 \ifinner\nondmatherr@\nopagebreak\else\postdisplaypenalty\@M\fi
 \else\ifvmode\nobreak\else\vadjust{\nobreak}\fi\fi}
\def\nonvmodeerr@#1{\Err@{\string#1\space not allowed within a paragraph
 or in math}}
\def\vnonvmode@#1#2{\relaxnext@\DNii@{\ifx\next\par\DN@{#1}\else
 \DN@{#2}\fi\next@}%
 \ifvmode\DN@{#1}\else
 \DN@{\FN@\nextii@}\fi\next@}
\def\newpage{\vnonvmode@{\vfill\break}{\nonvmodeerr@\newpage}}
\def\smallpagebreak{\vnonvmode@\smallbreak{\nonvmodeerr@\smallpagebreak}}
\def\medpagebreak{\vnonvmode@\medbreak{\nonvmodeerr@\medpagebreak}}
\def\bigpagebreak{\vnonvmode@\bigbreak{\nonvmodeerr@\bigpagebreak}}
\def\NoBlackBoxes{\global\overfullrule\z@}
\def\BlackBoxes{\global\overfullrule5\p@}
\def\Invalid@#1{\def#1{\Err@{\Invalid@@\string#1}}}
\def\Invalid@@{Invalid use of }
\message{figures,}
\Invalid@\caption
\Invalid@\captionwidth
\newdimen\smallcaptionwidth@
\def\topspace{\mid@false\ins@}
\def\midspace{\mid@true\ins@}
\newif\ifmid@
\def\captionfont@{}
\def\ins@#1{\relaxnext@\allowbreak
 \smallcaptionwidth@\captionwidth@\gdef\thespace@{#1}%
 \DN@{\ifx\next\space@\DN@. {\FN@\nextii@}\else
  \DN@.{\FN@\nextii@}\fi\next@.}%
 \DNii@{\ifx\next\caption\DN@\caption{\FN@\nextiii@}%
  \else\let\next@\nextiv@\fi\next@}%
 \def\nextiv@{\vnonvmode@
  {\ifmid@\expandafter\midinsert\else\expandafter\topinsert\fi
   \vbox to\thespace@{}\endinsert}
  {\ifmid@\nonvmodeerr@\midspace\else\nonvmodeerr@\topspace\fi}}%
 \def\nextiii@{\ifx\next\captionwidth\expandafter\nextv@
  \else\expandafter\nextvi@\fi}%
 \def\nextv@\captionwidth##1##2{\smallcaptionwidth@##1\relax\nextvi@{##2}}%
 \def\nextvi@##1{\def\thecaption@{\captionfont@##1}%
  \DN@{\ifx\next\space@\DN@. {\FN@\nextvii@}\else
   \DN@.{\FN@\nextvii@}\fi\next@.}%
  \FN@\next@}%
 \def\nextvii@{\vnonvmode@
  {\ifmid@\expandafter\midinsert\else
  \expandafter\topinsert\fi\vbox to\thespace@{}\nobreak\smallskip
  \setboxz@h{\noindent\ignorespaces\thecaption@\unskip}%
  \ifdim\wdz@>\smallcaptionwidth@\centerline{\vbox{\hsize\smallcaptionwidth@
   \noindent\ignorespaces\thecaption@\unskip}}%
  \else\centerline{\boxz@}\fi\endinsert}
  {\ifmid@\nonvmodeerr@\midspace
  \else\nonvmodeerr@\topspace\fi}}%
 \FN@\next@}
\message{comments,}
\def\newcodes@{\catcode`\\12\catcode`\{12\catcode`\}12\catcode`\#12%
 \catcode`\%12\relax}
\def\oldcodes@{\catcode`\\0\catcode`\{1\catcode`\}2\catcode`\#6%
 \catcode`\%14\relax}
\def\comment{\newcodes@\endlinechar=10 \comment@}
{\lccode`\0=`\\
\lowercase{\gdef\comment@#1^^J{\comment@@#10endcomment\comment@@@}%
\gdef\comment@@#10endcomment{\FN@\comment@@@}%
\gdef\comment@@@#1\comment@@@{\ifx\next\comment@@@\let\next\comment@
 \else\def\next{\oldcodes@\endlinechar=`\^^M\relax}%
 \fi\next}}}
\def\pr@m@s{\ifx'\next\DN@##1{\prim@s}\else\let\next@\egroup\fi\next@}
\def\prime{{\null\prime@\null}}
\mathchardef\prime@="0230
\let\dsize\displaystyle

\let\ssize\scriptstyle

\message{math spacing,}
\def\,{\RIfM@\mskip\thinmuskip\relax\else\kern.16667em\fi}
\def\!{\RIfM@\mskip-\thinmuskip\relax\else\kern-.16667em\fi}
\let\thinspace\,
\let\negthinspace\!
\def\medspace{\RIfM@\mskip\medmuskip\relax\else\kern.222222em\fi}
\def\negmedspace{\RIfM@\mskip-\medmuskip\relax\else\kern-.222222em\fi}
\def\thickspace{\RIfM@\mskip\thickmuskip\relax\else\kern.27777em\fi}
\let\;\thickspace
\def\negthickspace{\RIfM@\mskip-\thickmuskip\relax\else
 \kern-.27777em\fi}
\atdef@,{\RIfM@\mskip.1\thinmuskip\else\leavevmode\null,\fi}
\atdef@!{\RIfM@\mskip-.1\thinmuskip\else\leavevmode\null!\fi}
\atdef@.{\RIfM@&&\else\leavevmode.\spacefactor3000 \fi}
\def\and{\DOTSB\;\mathchar"3026 \;}

\message{fractions,}
\def\frac#1#2{{#1\over#2}}

\newdimen\ex@
\ex@.2326ex
\Invalid@\thickness
\def\thickfrac{\relaxnext@
 \DN@{\ifx\next\thickness\let\next@\nextii@\else
 \DN@{\nextii@\thickness1}\fi\next@}%
 \DNii@\thickness##1##2##3{{##2\above##1\ex@##3}}%
 \FN@\next@}

\def\thickfracwithdelims#1#2{\relaxnext@\def\ldelim@{#1}\def\rdelim@{#2}%
 \DN@{\ifx\next\thickness\let\next@\nextii@\else
 \DN@{\nextii@\thickness1}\fi\next@}%
 \DNii@\thickness##1##2##3{{##2\abovewithdelims
 \ldelim@\rdelim@##1\ex@##3}}%
 \FN@\next@}
\def\binom#1#2{{#1\choose#2}}

\def\:{\nobreak\hskip.1111em\mathpunct{}\nonscript\mkern-\thinmuskip{:}\hskip
 .3333emplus.0555em\relax}
\def\snug{\unskip\kern-\mathsurround}
\message{smash commands,}
\def\topsmash{\top@true\bot@false\smash@}
\def\botsmash{\top@false\bot@true\smash@}
\newif\iftop@
\newif\ifbot@
\def\smash{\top@true\bot@true\smash@}
\def\smash@{\RIfM@\expandafter\mathpalette\expandafter\mathsm@sh\else
 \expandafter\makesm@sh\fi}
\def\finsm@sh{\iftop@\ht\z@\z@\fi\ifbot@\dp\z@\z@\fi\leavevmode\boxz@}
\message{large operator symbols,}
\def\LimitsOnSums{\global\let\slimits@\displaylimits}
\def\NoLimitsOnSums{\global\let\slimits@\nolimits}
\LimitsOnSums
\mathchardef\coprod@="1360       \def\coprod{\DOTSB\coprod@\slimits@}
\mathchardef\bigvee@="1357       \def\bigvee{\DOTSB\bigvee@\slimits@}
\mathchardef\bigwedge@="1356     \def\bigwedge{\DOTSB\bigwedge@\slimits@}
\mathchardef\biguplus@="1355     \def\biguplus{\DOTSB\biguplus@\slimits@}
\mathchardef\bigcap@="1354       \def\bigcap{\DOTSB\bigcap@\slimits@}
\mathchardef\bigcup@="1353       \def\bigcup{\DOTSB\bigcup@\slimits@}
\mathchardef\prod@="1351         \def\prod{\DOTSB\prod@\slimits@}
\mathchardef\sum@="1350          \def\sum{\DOTSB\sum@\slimits@}
\mathchardef\bigotimes@="134E    \def\bigotimes{\DOTSB\bigotimes@\slimits@}
\mathchardef\bigoplus@="134C     \def\bigoplus{\DOTSB\bigoplus@\slimits@}
\mathchardef\bigodot@="134A      \def\bigodot{\DOTSB\bigodot@\slimits@}
\mathchardef\bigsqcup@="1346     \def\bigsqcup{\DOTSB\bigsqcup@\slimits@}
\message{integrals,}
\def\LimitsOnInts{\global\let\ilimits@\displaylimits}
\def\NoLimitsOnInts{\global\let\ilimits@\nolimits}
\NoLimitsOnInts
\def\int{\DOTSI\intop\ilimits@}
\def\oint{\DOTSI\ointop\ilimits@}
\def\intic@{\mathchoice{\hskip.5em}{\hskip.4em}{\hskip.4em}{\hskip.4em}}
\def\negintic@{\mathchoice
 {\hskip-.5em}{\hskip-.4em}{\hskip-.4em}{\hskip-.4em}}
\def\intkern@{\mathchoice{\!\!\!}{\!\!}{\!\!}{\!\!}}
\def\intdots@{\mathchoice{\plaincdots@}
 {{\cdotp}\mkern1.5mu{\cdotp}\mkern1.5mu{\cdotp}}
 {{\cdotp}\mkern1mu{\cdotp}\mkern1mu{\cdotp}}
 {{\cdotp}\mkern1mu{\cdotp}\mkern1mu{\cdotp}}}
\newcount\intno@
\def\iint{\DOTSI\intno@\tw@\FN@\ints@}
\def\iiint{\DOTSI\intno@\thr@@\FN@\ints@}
\def\iiiint{\DOTSI\intno@4 \FN@\ints@}
\def\idotsint{\DOTSI\intno@\z@\FN@\ints@}
\def\ints@{\findlimits@\ints@@}
\newif\iflimtoken@
\newif\iflimits@
\def\findlimits@{\limtoken@true\ifx\next\limits\limits@true
 \else\ifx\next\nolimits\limits@false\else
 \limtoken@false\ifx\ilimits@\nolimits\limits@false\else
 \ifinner\limits@false\else\limits@true\fi\fi\fi\fi}
\def\multint@{\int\ifnum\intno@=\z@\intdots@                                
 \else\intkern@\fi                                                          
 \ifnum\intno@>\tw@\int\intkern@\fi                                         
 \ifnum\intno@>\thr@@\int\intkern@\fi                                       
 \int}                                                                      
\def\multintlimits@{\intop\ifnum\intno@=\z@\intdots@\else\intkern@\fi
 \ifnum\intno@>\tw@\intop\intkern@\fi
 \ifnum\intno@>\thr@@\intop\intkern@\fi\intop}
\def\ints@@{\iflimtoken@                                                    
 \def\ints@@@{\iflimits@\negintic@\mathop{\intic@\multintlimits@}\limits    
  \else\multint@\nolimits\fi                                                
  \eat@}                                                                    
 \else                                                                      
 \def\ints@@@{\iflimits@\negintic@
  \mathop{\intic@\multintlimits@}\limits\else
  \multint@\nolimits\fi}\fi\ints@@@}
\def\LimitsOnNames{\global\let\nlimits@\displaylimits}
\def\NoLimitsOnNames{\global\let\nlimits@\nolimits@}
\LimitsOnNames
\def\nolimits@{\relaxnext@
 \DN@{\ifx\next\limits\DN@\limits{\nolimits}\else
  \let\next@\nolimits\fi\next@}%
 \FN@\next@}
\message{operator names,}
\def\newmcodes@{\mathcode`\'"27\mathcode`\*"2A\mathcode`\."613A%
 \mathcode`\-"2D\mathcode`\/"2F\mathcode`\:"603A }
\def\operatorname#1{\mathop{\newmcodes@\kern\z@\fam\z@#1}\nolimits@}
\def\operatornamewithlimits#1{\mathop{\newmcodes@\kern\z@\fam\z@#1}\nlimits@}
\def\qopname@#1{\mathop{\fam\z@#1}\nolimits@}
\def\qopnamewl@#1{\mathop{\fam\z@#1}\nlimits@}
\def\arccos{\qopname@{arccos}}
\def\arcsin{\qopname@{arcsin}}
\def\arctan{\qopname@{arctan}}
\def\arg{\qopname@{arg}}
\def\cos{\qopname@{cos}}
\def\cosh{\qopname@{cosh}}
\def\cot{\qopname@{cot}}
\def\coth{\qopname@{coth}}
\def\csc{\qopname@{csc}}
\def\deg{\qopname@{deg}}
\def\det{\qopnamewl@{det}}
\def\dim{\qopname@{dim}}
\def\exp{\qopname@{exp}}
\def\gcd{\qopnamewl@{gcd}}
\def\hom{\qopname@{hom}}
\def\inf{\qopnamewl@{inf}}
\def\injlim{\qopnamewl@{inj\,lim}}
\def\ker{\qopname@{ker}}
\def\lg{\qopname@{lg}}
\def\lim{\qopnamewl@{lim}}
\def\liminf{\qopnamewl@{lim\,inf}}
\def\limsup{\qopnamewl@{lim\,sup}}
\def\ln{\qopname@{ln}}
\def\log{\qopname@{log}}
\def\max{\qopnamewl@{max}}
\def\min{\qopnamewl@{min}}
\def\Pr{\qopnamewl@{Pr}}
\def\projlim{\qopnamewl@{proj\,lim}}
\def\sec{\qopname@{sec}}
\def\sin{\qopname@{sin}}
\def\sinh{\qopname@{sinh}}
\def\sup{\qopnamewl@{sup}}
\def\tan{\qopname@{tan}}
\def\tanh{\qopname@{tanh}}
\def\varinjlim{\mathop{\vtop{\ialign{##\crcr
 \hfil\rm lim\hfil\crcr\noalign{\nointerlineskip}\rightarrowfill\crcr
 \noalign{\nointerlineskip\kern-\ex@}\crcr}}}}
\def\varprojlim{\mathop{\vtop{\ialign{##\crcr
 \hfil\rm lim\hfil\crcr\noalign{\nointerlineskip}\leftarrowfill\crcr
 \noalign{\nointerlineskip\kern-\ex@}\crcr}}}}
\def\varliminf{\mathop{\underline{\vrule height\z@ depth.2exwidth\z@
 \hbox{\rm lim}}}}

\newdimen\buffer@
\buffer@\fontdimen13 \tenex
\newdimen\buffer
\buffer\buffer@

\def\ResetBuffer{\fontdimen13 \tenex\buffer@\global\buffer\buffer@}
\def\shave#1{\mathop{\hbox{$\m@th\fontdimen13 \tenex\z@                     
 \displaystyle{#1}$}}\fontdimen13 \tenex\buffer}

\message{multilevel sub/superscripts,}
\Invalid@\\
\def\Let@{\relax\iffalse{\fi\let\\=\cr\iffalse}\fi}
\Invalid@\vspace
\def\vspace@{\def\vspace##1{\crcr\noalign{\vskip##1\relax}}}
\def\multilimits@{\bgroup\vspace@\Let@
 \baselineskip\fontdimen10 \scriptfont\tw@
 \advance\baselineskip\fontdimen12 \scriptfont\tw@
 \lineskip\thr@@\fontdimen8 \scriptfont\thr@@
 \lineskiplimit\lineskip
 \vbox\bgroup\ialign\bgroup\hfil$\m@th\scriptstyle{##}$\hfil\crcr}
\def\Sb{_\multilimits@}
\def\endSb{\crcr\egroup\egroup\egroup}
\def\Sp{^\multilimits@}

\def\spreadlines#1{\RIfMIfI@\onlydmatherr@\spreadlines\else
 \openup#1\relax\fi\else\onlydmatherr@\spreadlines\fi}
\def\Mathstrut@{\copy\Mathstrutbox@}
\newbox\Mathstrutbox@
\setbox\Mathstrutbox@\null
\setboxz@h{$\m@th($}
\ht\Mathstrutbox@\ht\z@
\dp\Mathstrutbox@\dp\z@
\message{matrices,}
\newdimen\spreadmlines@
\def\spreadmatrixlines#1{\RIfMIfI@
 \onlydmatherr@\spreadmatrixlines\else
 \spreadmlines@#1\relax\fi\else\onlydmatherr@\spreadmatrixlines\fi}
\def\matrix{\null\,\vcenter\bgroup\Let@\vspace@
 \normalbaselines\openup\spreadmlines@\ialign
 \bgroup\hfil$\m@th##$\hfil&&\quad\hfil$\m@th##$\hfil\crcr
 \Mathstrut@\crcr\noalign{\kern-\baselineskip}}
\def\endmatrix{\crcr\Mathstrut@\crcr\noalign{\kern-\baselineskip}\egroup
 \egroup\,}
\def\format{\crcr\egroup\iffalse{\fi\ifnum`}=0 \fi\format@}
\newtoks\hashtoks@
\hashtoks@{#}
\def\format@#1\\{\def\preamble@{#1}%
 \def\l{$\m@th\the\hashtoks@$\hfil}%
 \def\c{\hfil$\m@th\the\hashtoks@$\hfil}%
 \def\r{\hfil$\m@th\the\hashtoks@$}%
 \edef\preamble@@{\preamble@}\ifnum`{=0 \fi\iffalse}\fi
 \ialign\bgroup\span\preamble@@\crcr}
\def\smallmatrix{\null\,\vcenter\bgroup\vspace@\Let@
 \baselineskip9\ex@\lineskip\ex@
 \ialign\bgroup\hfil$\m@th\scriptstyle{##}$\hfil&&\thickspace\hfil
 $\m@th\scriptstyle{##}$\hfil\crcr}
\def\endsmallmatrix{\crcr\egroup\egroup\,}

\newmuskip\dotsspace@
\dotsspace@1.5mu
\def\strip@#1 {#1}
\def\spacehdots#1\for#2{\multispan{#2}\xleaders
 \hbox{$\m@th\mkern\strip@#1 \dotsspace@.\mkern\strip@#1 \dotsspace@$}\hfill}
\def\hdotsfor#1{\spacehdots\@ne\for{#1}}
\def\multispan@#1{\omit\mscount#1\unskip\loop\ifnum\mscount>\@ne\sp@n\repeat}
\def\spaceinnerhdots#1\for#2\after#3{\multispan@{\strip@#2 }#3\xleaders
 \hbox{$\m@th\mkern\strip@#1 \dotsspace@.\mkern\strip@#1 \dotsspace@$}\hfill}
\def\innerhdotsfor#1\after#2{\spaceinnerhdots\@ne\for#1\after{#2}}
\def\cases{\bgroup\spreadmlines@\jot\left\{\,\matrix\format\l&\quad\l\\}
\def\endcases{\endmatrix\right.\egroup}
\message{multiline displays,}
\newif\ifinany@
\newif\ifinalign@
\newif\ifingather@
\def\strut@{\copy\strutbox@}
\newbox\strutbox@
\setbox\strutbox@\hbox{\vrule height8\p@ depth3\p@ width\z@}
\def\topaligned{\null\,\vtop\aligned@}
\def\botaligned{\null\,\vbox\aligned@}
\def\aligned{\null\,\vcenter\aligned@}
\def\aligned@{\bgroup\vspace@\Let@
 \ifinany@\else\openup\jot\fi\ialign
 \bgroup\hfil\strut@$\m@th\displaystyle{##}$&
 $\m@th\displaystyle{{}##}$\hfil\crcr}
\def\endaligned{\crcr\egroup\egroup}

\def\alignedat#1{\null\,\vcenter\bgroup\doat@{#1}\vspace@\Let@
 \ifinany@\else\openup\jot\fi\ialign\bgroup\span\preamble@@\crcr}
\newcount\atcount@
\def\doat@#1{\toks@{\hfil\strut@$\m@th
 \displaystyle{\the\hashtoks@}$&$\m@th\displaystyle
 {{}\the\hashtoks@}$\hfil}
 \atcount@#1\relax\advance\atcount@\m@ne                                    
 \loop\ifnum\atcount@>\z@\toks@=\expandafter{\the\toks@&\hfil$\m@th
 \displaystyle{\the\hashtoks@}$&$\m@th
 \displaystyle{{}\the\hashtoks@}$\hfil}\advance
  \atcount@\m@ne\repeat                                                     
 \xdef\preamble@{\the\toks@}\xdef\preamble@@{\preamble@}}

\def\gathered{\null\,\vcenter\bgroup\vspace@\Let@
 \ifinany@\else\openup\jot\fi\ialign
 \bgroup\hfil\strut@$\m@th\displaystyle{##}$\hfil\crcr}
\def\endgathered{\crcr\egroup\egroup}
\newif\iftagsleft@
\def\TagsOnLeft{\global\tagsleft@true}
\def\TagsOnRight{\global\tagsleft@false}
\TagsOnLeft
\newif\ifmathtags@
\def\TagsAsMath{\global\mathtags@true}
\def\TagsAsText{\global\mathtags@false}
\TagsAsText
\def\tagform@#1{\hbox{\rm(\ignorespaces#1\unskip)}}
\def\thetag{\leavevmode\tagform@}
\def\tag#1$${\iftagsleft@\leqno\else\eqno\fi                                
 \maketag@#1\maketag@                                                       
 $$}                                                                        
\def\maketag@{\FN@\maketag@@}
\def\maketag@@{\ifx\next"\expandafter\maketag@@@\else\expandafter\maketag@@@@
 \fi}
\def\maketag@@@"#1"#2\maketag@{\hbox{\rm#1}}                                
\def\maketag@@@@#1\maketag@{\ifmathtags@\tagform@{$\m@th#1$}\else
 \tagform@{#1}\fi}
\interdisplaylinepenalty\@M
\def\allowdisplaybreaks{\RIfMIfI@
 \onlydmatherr@\allowdisplaybreaks\else
 \interdisplaylinepenalty\z@\fi\else\onlydmatherr@\allowdisplaybreaks\fi}
\Invalid@\allowdisplaybreak
\Invalid@\displaybreak
\Invalid@\intertext
\def\allowdisplaybreak@{\def\allowdisplaybreak{\crcr\noalign{\allowbreak}}}
\def\displaybreak@{\def\displaybreak{\crcr\noalign{\break}}}
\def\intertext@{\def\intertext##1{\crcr\noalign{%
 \penalty\postdisplaypenalty \vskip\belowdisplayskip
 \vbox{\normalbaselines\noindent##1}%
 \penalty\predisplaypenalty \vskip\abovedisplayskip}}}
\newskip\centering@
\centering@\z@ plus\@m\p@
\def\align{\relax\ifingather@\DN@{\csname align (in
  \string\gather)\endcsname}\else
 \ifmmode\ifinner\DN@{\onlydmatherr@\align}\else
  \let\next@\align@\fi
 \else\DN@{\onlydmatherr@\align}\fi\fi\next@}
\newhelp\andhelp@
{An extra & here is so disastrous that you should probably exit^^J
and fix things up.}
\newif\iftag@
\newcount\and@
\def\align@{\inalign@true\inany@true
 \vspace@\allowdisplaybreak@\displaybreak@\intertext@
 \def\tag{\global\tag@true\ifnum\and@=\z@\DN@{&&}\else
          \DN@{&}\fi\next@}%
 \iftagsleft@\DN@{\csname align \endcsname}\else
  \DN@{\csname align \space\endcsname}\fi\next@}
\def\Tag@{\iftag@\else\errhelp\andhelp@\err@{Extra & on this line}\fi}
\newdimen\lwidth@
\newdimen\rwidth@
\newdimen\maxlwidth@
\newdimen\maxrwidth@
\newdimen\totwidth@
\def\measure@#1\endalign{\lwidth@\z@\rwidth@\z@\maxlwidth@\z@\maxrwidth@\z@
 \global\and@\z@                                                            
 \setbox@ne\vbox                                                            
  {\everycr{\noalign{\global\tag@false\global\and@\z@}}\Let@                
  \halign{\setboxz@h{$\m@th\displaystyle{\@lign##}$}
   \global\lwidth@\wdz@                                                     
   \ifdim\lwidth@>\maxlwidth@\global\maxlwidth@\lwidth@\fi                  
   \global\advance\and@\@ne                                                 
   &\setboxz@h{$\m@th\displaystyle{{}\@lign##}$}\global\rwidth@\wdz@        
   \ifdim\rwidth@>\maxrwidth@\global\maxrwidth@\rwidth@\fi                  
   \global\advance\and@\@ne                                                
   &\Tag@
   \eat@{##}\crcr#1\crcr}}
 \totwidth@\maxlwidth@\advance\totwidth@\maxrwidth@}                       
\def\displ@y@{\global\dt@ptrue\openup\jot
 \everycr{\noalign{\global\tag@false\global\and@\z@\ifdt@p\global\dt@pfalse
 \vskip-\lineskiplimit\vskip\normallineskiplimit\else
 \penalty\interdisplaylinepenalty\fi}}}
\def\black@#1{\noalign{\ifdim#1>\displaywidth
 \dimen@\prevdepth\nointerlineskip                                          
 \vskip-\ht\strutbox@\vskip-\dp\strutbox@                                   
 \vbox{\noindent\hbox to#1{\strut@\hfill}}
 \prevdepth\dimen@                                                          
 \fi}}
\expandafter\def\csname align \space\endcsname#1\endalign
 {\measure@#1\endalign\global\and@\z@                                       
 \ifingather@\everycr{\noalign{\global\and@\z@}}\else\displ@y@\fi           
 \Let@\tabskip\centering@                                                   
 \halign to\displaywidth
  {\hfil\strut@\setboxz@h{$\m@th\displaystyle{\@lign##}$}
  \global\lwidth@\wdz@\boxz@\global\advance\and@\@ne                        
  \tabskip\z@skip                                                           
  &\setboxz@h{$\m@th\displaystyle{{}\@lign##}$}
  \global\rwidth@\wdz@\boxz@\hfill\global\advance\and@\@ne                  
  \tabskip\centering@                                                       
  &\setboxz@h{\@lign\strut@\maketag@##\maketag@}
  \dimen@\displaywidth\advance\dimen@-\totwidth@
  \divide\dimen@\tw@\advance\dimen@\maxrwidth@\advance\dimen@-\rwidth@     
  \ifdim\dimen@<\tw@\wdz@\llap{\vtop{\normalbaselines\null\boxz@}}
  \else\llap{\boxz@}\fi                                                    
  \tabskip\z@skip                                                          
  \crcr#1\crcr                                                             
  \black@\totwidth@}}                                                      
\newdimen\lineht@
\expandafter\def\csname align \endcsname#1\endalign{\measure@#1\endalign
 \global\and@\z@
 \ifdim\totwidth@>\displaywidth\let\displaywidth@\totwidth@\else
  \let\displaywidth@\displaywidth\fi                                        
 \ifingather@\everycr{\noalign{\global\and@\z@}}\else\displ@y@\fi
 \Let@\tabskip\centering@\halign to\displaywidth
  {\hfil\strut@\setboxz@h{$\m@th\displaystyle{\@lign##}$}%
  \global\lwidth@\wdz@\global\lineht@\ht\z@                                 
  \boxz@\global\advance\and@\@ne
  \tabskip\z@skip&\setboxz@h{$\m@th\displaystyle{{}\@lign##}$}%
  \global\rwidth@\wdz@\ifdim\ht\z@>\lineht@\global\lineht@\ht\z@\fi         
  \boxz@\hfil\global\advance\and@\@ne
  \tabskip\centering@&\kern-\displaywidth@                                  
  \setboxz@h{\@lign\strut@\maketag@##\maketag@}%
  \dimen@\displaywidth\advance\dimen@-\totwidth@
  \divide\dimen@\tw@\advance\dimen@\maxlwidth@\advance\dimen@-\lwidth@
  \ifdim\dimen@<\tw@\wdz@
   \rlap{\vbox{\normalbaselines\boxz@\vbox to\lineht@{}}}\else
   \rlap{\boxz@}\fi
  \tabskip\displaywidth@\crcr#1\crcr\black@\totwidth@}}
\expandafter\def\csname align (in \string\gather)\endcsname
  #1\endalign{\vcenter{\align@#1\endalign}}
\Invalid@\endalign
\newif\ifxat@
\def\alignat{\RIfMIfI@\DN@{\onlydmatherr@\alignat}\else
 \DN@{\csname alignat \endcsname}\fi\else
 \DN@{\onlydmatherr@\alignat}\fi\next@}
\newif\ifmeasuring@
\newbox\savealignat@
\expandafter\def\csname alignat \endcsname#1#2\endalignat                   
 {\inany@true\xat@false
 \def\tag{\global\tag@true\count@#1\relax\multiply\count@\tw@
  \xdef\tag@{}\loop\ifnum\count@>\and@\xdef\tag@{&\tag@}\advance\count@\m@ne
  \repeat\tag@}%
 \vspace@\allowdisplaybreak@\displaybreak@\intertext@
 \displ@y@\measuring@true                                                   
 \setbox\savealignat@\hbox{$\m@th\displaystyle\Let@
  \attag@{#1}
  \vbox{\halign{\span\preamble@@\crcr#2\crcr}}$}%
 \measuring@false                                                           
 \Let@\attag@{#1}
 \tabskip\centering@\halign to\displaywidth
  {\span\preamble@@\crcr#2\crcr                                             
  \black@{\wd\savealignat@}}}                                               
\Invalid@\endalignat
\def\xalignat{\RIfMIfI@
 \DN@{\onlydmatherr@\xalignat}\else
 \DN@{\csname xalignat \endcsname}\fi\else
 \DN@{\onlydmatherr@\xalignat}\fi\next@}
\expandafter\def\csname xalignat \endcsname#1#2\endxalignat
 {\inany@true\xat@true
 \def\tag{\global\tag@true\def\tag@{}\count@#1\relax\multiply\count@\tw@
  \loop\ifnum\count@>\and@\xdef\tag@{&\tag@}\advance\count@\m@ne\repeat\tag@}%
 \vspace@\allowdisplaybreak@\displaybreak@\intertext@
 \displ@y@\measuring@true\setbox\savealignat@\hbox{$\m@th\displaystyle\Let@
 \attag@{#1}\vbox{\halign{\span\preamble@@\crcr#2\crcr}}$}%
 \measuring@false\Let@
 \attag@{#1}\tabskip\centering@\halign to\displaywidth
 {\span\preamble@@\crcr#2\crcr\black@{\wd\savealignat@}}}
\def\attag@#1{\let\Maketag@\maketag@\let\TAG@\Tag@                          
 \let\Tag@=0\let\maketag@=0
 \ifmeasuring@\def\llap@##1{\setboxz@h{##1}\hbox to\tw@\wdz@{}}%
  \def\rlap@##1{\setboxz@h{##1}\hbox to\tw@\wdz@{}}\else
  \let\llap@\llap\let\rlap@\rlap\fi                                         
 \toks@{\hfil\strut@$\m@th\displaystyle{\@lign\the\hashtoks@}$\tabskip\z@skip
  \global\advance\and@\@ne&$\m@th\displaystyle{{}\@lign\the\hashtoks@}$\hfil
  \ifxat@\tabskip\centering@\fi\global\advance\and@\@ne}
 \iftagsleft@
  \toks@@{\tabskip\centering@&\Tag@\kern-\displaywidth
   \rlap@{\@lign\maketag@\the\hashtoks@\maketag@}%
   \global\advance\and@\@ne\tabskip\displaywidth}\else
  \toks@@{\tabskip\centering@&\Tag@\llap@{\@lign\maketag@
   \the\hashtoks@\maketag@}\global\advance\and@\@ne\tabskip\z@skip}\fi      
 \atcount@#1\relax\advance\atcount@\m@ne
 \loop\ifnum\atcount@>\z@
 \toks@=\expandafter{\the\toks@&\hfil$\m@th\displaystyle{\@lign
  \the\hashtoks@}$\global\advance\and@\@ne
  \tabskip\z@skip&$\m@th\displaystyle{{}\@lign\the\hashtoks@}$\hfil\ifxat@
  \tabskip\centering@\fi\global\advance\and@\@ne}\advance\atcount@\m@ne
 \repeat                                                                    
 \xdef\preamble@{\the\toks@\the\toks@@}
 \xdef\preamble@@{\preamble@}
 \let\maketag@\Maketag@\let\Tag@\TAG@}                                      
\Invalid@\endxalignat
\def\xxalignat{\RIfMIfI@
 \DN@{\onlydmatherr@\xxalignat}\else\DN@{\csname xxalignat
  \endcsname}\fi\else
 \DN@{\onlydmatherr@\xxalignat}\fi\next@}
\expandafter\def\csname xxalignat \endcsname#1#2\endxxalignat{\inany@true
 \vspace@\allowdisplaybreak@\displaybreak@\intertext@
 \displ@y\setbox\savealignat@\hbox{$\m@th\displaystyle\Let@
 \xxattag@{#1}\vbox{\halign{\span\preamble@@\crcr#2\crcr}}$}%
 \Let@\xxattag@{#1}\tabskip\z@skip\halign to\displaywidth
 {\span\preamble@@\crcr#2\crcr\black@{\wd\savealignat@}}}
\def\xxattag@#1{\toks@{\tabskip\z@skip\hfil\strut@
 $\m@th\displaystyle{\the\hashtoks@}$&%
 $\m@th\displaystyle{{}\the\hashtoks@}$\hfil\tabskip\centering@&}%
 \atcount@#1\relax\advance\atcount@\m@ne\loop\ifnum\atcount@>\z@
 \toks@=\expandafter{\the\toks@&\hfil$\m@th\displaystyle{\the\hashtoks@}$%
  \tabskip\z@skip&$\m@th\displaystyle{{}\the\hashtoks@}$\hfil
  \tabskip\centering@}\advance\atcount@\m@ne\repeat
 \xdef\preamble@{\the\toks@\tabskip\z@skip}\xdef\preamble@@{\preamble@}}
\Invalid@\endxxalignat
\newdimen\gwidth@
\newdimen\gmaxwidth@
\def\gmeasure@#1\endgather{\gwidth@\z@\gmaxwidth@\z@\setbox@ne\vbox{\Let@
 \halign{\setboxz@h{$\m@th\displaystyle{##}$}\global\gwidth@\wdz@
 \ifdim\gwidth@>\gmaxwidth@\global\gmaxwidth@\gwidth@\fi
 &\eat@{##}\crcr#1\crcr}}}
\def\gather{\RIfMIfI@\DN@{\onlydmatherr@\gather}\else
 \ingather@true\inany@true\def\tag{&}%
 \vspace@\allowdisplaybreak@\displaybreak@\intertext@
 \displ@y\Let@
 \iftagsleft@\DN@{\csname gather \endcsname}\else
  \DN@{\csname gather \space\endcsname}\fi\fi
 \else\DN@{\onlydmatherr@\gather}\fi\next@}
\expandafter\def\csname gather \space\endcsname#1\endgather
 {\gmeasure@#1\endgather\tabskip\centering@
 \halign to\displaywidth{\hfil\strut@\setboxz@h{$\m@th\displaystyle{##}$}%
 \global\gwidth@\wdz@\boxz@\hfil&
 \setboxz@h{\strut@{\maketag@##\maketag@}}%
 \dimen@\displaywidth\advance\dimen@-\gwidth@
 \ifdim\dimen@>\tw@\wdz@\llap{\boxz@}\else
 \llap{\vtop{\normalbaselines\null\boxz@}}\fi
 \tabskip\z@skip\crcr#1\crcr\black@\gmaxwidth@}}
\newdimen\glineht@
\expandafter\def\csname gather \endcsname#1\endgather{\gmeasure@#1\endgather
 \ifdim\gmaxwidth@>\displaywidth\let\gdisplaywidth@\gmaxwidth@\else
 \let\gdisplaywidth@\displaywidth\fi\tabskip\centering@\halign to\displaywidth
 {\hfil\strut@\setboxz@h{$\m@th\displaystyle{##}$}%
 \global\gwidth@\wdz@\global\glineht@\ht\z@\boxz@\hfil&\kern-\gdisplaywidth@
 \setboxz@h{\strut@{\maketag@##\maketag@}}%
 \dimen@\displaywidth\advance\dimen@-\gwidth@
 \ifdim\dimen@>\tw@\wdz@\rlap{\boxz@}\else
 \rlap{\vbox{\normalbaselines\boxz@\vbox to\glineht@{}}}\fi
 \tabskip\gdisplaywidth@\crcr#1\crcr\black@\gmaxwidth@}}
\newif\ifctagsplit@
\def\CenteredTagsOnSplits{\global\ctagsplit@true}
\def\TopOrBottomTagsOnSplits{\global\ctagsplit@false}
\TopOrBottomTagsOnSplits
\def\split{\relax\ifinany@\let\next@\insplit@\else
 \ifmmode\ifinner\def\next@{\onlydmatherr@\split}\else
 \let\next@\outsplit@\fi\else
 \def\next@{\onlydmatherr@\split}\fi\fi\next@}
\def\insplit@{\global\setbox\z@\vbox\bgroup\vspace@\Let@\ialign\bgroup
 \hfil\strut@$\m@th\displaystyle{##}$&$\m@th\displaystyle{{}##}$\hfill\crcr}
\def\endsplit{\crcr\egroup\egroup\iftagsleft@\expandafter\lendsplit@\else
 \expandafter\rendsplit@\fi}
\def\rendsplit@{\global\setbox9 \vbox
 {\unvcopy\z@\global\setbox8 \lastbox\unskip}
 \setbox@ne\hbox{\unhcopy8 \unskip\global\setbox\tw@\lastbox
 \unskip\global\setbox\thr@@\lastbox}
 \global\setbox7 \hbox{\unhbox\tw@\unskip}
 \ifinalign@\ifctagsplit@                                                   
  \gdef\split@{\hbox to\wd\thr@@{}&
   \vcenter{\vbox{\moveleft\wd\thr@@\boxz@}}}
 \else\gdef\split@{&\vbox{\moveleft\wd\thr@@\box9}\crcr
  \box\thr@@&\box7}\fi                                                      
 \else                                                                      
  \ifctagsplit@\gdef\split@{\vcenter{\boxz@}}\else
  \gdef\split@{\box9\crcr\hbox{\box\thr@@\box7}}\fi
 \fi
 \split@}                                                                   
\def\lendsplit@{\global\setbox9\vtop{\unvcopy\z@}
 \setbox@ne\vbox{\unvcopy\z@\global\setbox8\lastbox}
 \setbox@ne\hbox{\unhcopy8\unskip\setbox\tw@\lastbox
  \unskip\global\setbox\thr@@\lastbox}
 \ifinalign@\ifctagsplit@                                                   
  \gdef\split@{\hbox to\wd\thr@@{}&
  \vcenter{\vbox{\moveleft\wd\thr@@\box9}}}
  \else                                                                     
  \gdef\split@{\hbox to\wd\thr@@{}&\vbox{\moveleft\wd\thr@@\box9}}\fi
 \else
  \ifctagsplit@\gdef\split@{\vcenter{\box9}}\else
  \gdef\split@{\box9}\fi
 \fi\split@}
\def\outsplit@#1$${\align\insplit@#1\endalign$$}
\newdimen\multlinegap@
\multlinegap@1em
\newdimen\multlinetaggap@
\multlinetaggap@1em
\def\MultlineGap#1{\global\multlinegap@#1\relax}
\def\multlinegap#1{\RIfMIfI@\onlydmatherr@\multlinegap\else
 \multlinegap@#1\relax\fi\else\onlydmatherr@\multlinegap\fi}
\def\nomultlinegap{\multlinegap{\z@}}
\def\multline{\RIfMIfI@
 \DN@{\onlydmatherr@\multline}\else
 \DN@{\multline@}\fi\else
 \DN@{\onlydmatherr@\multline}\fi\next@}
\newif\iftagin@
\def\tagin@#1{\tagin@false\in@\tag{#1}\ifin@\tagin@true\fi}
\def\multline@#1$${\inany@true\vspace@\allowdisplaybreak@\displaybreak@
 \tagin@{#1}\iftagsleft@\DN@{\multline@l#1$$}\else
 \DN@{\multline@r#1$$}\fi\next@}
\newdimen\mwidth@
\def\rmmeasure@#1\endmultline{%
 \def\shoveleft##1{##1}\def\shoveright##1{##1}
 \setbox@ne\vbox{\Let@\halign{\setboxz@h
  {$\m@th\@lign\displaystyle{}##$}\global\mwidth@\wdz@
  \crcr#1\crcr}}}
\newdimen\mlineht@
\newif\ifzerocr@
\newif\ifonecr@
\def\lmmeasure@#1\endmultline{\global\zerocr@true\global\onecr@false
 \everycr{\noalign{\ifonecr@\global\onecr@false\fi
  \ifzerocr@\global\zerocr@false\global\onecr@true\fi}}
  \def\shoveleft##1{##1}\def\shoveright##1{##1}%
 \setbox@ne\vbox{\Let@\halign{\setboxz@h
  {$\m@th\@lign\displaystyle{}##$}\ifonecr@\global\mwidth@\wdz@
  \global\mlineht@\ht\z@\fi\crcr#1\crcr}}}
\newbox\mtagbox@
\newdimen\ltwidth@
\newdimen\rtwidth@
\def\multline@l#1$${\iftagin@\DN@{\lmultline@@#1$$}\else
 \DN@{\setbox\mtagbox@\null\ltwidth@\z@\rtwidth@\z@
  \lmultline@@@#1$$}\fi\next@}
\def\lmultline@@#1\endmultline\tag#2$${%
 \setbox\mtagbox@\hbox{\maketag@#2\maketag@}
 \lmmeasure@#1\endmultline\dimen@\mwidth@\advance\dimen@\wd\mtagbox@
 \advance\dimen@\multlinetaggap@                                            
 \ifdim\dimen@>\displaywidth\ltwidth@\z@\else\ltwidth@\wd\mtagbox@\fi       
 \lmultline@@@#1\endmultline$$}
\def\lmultline@@@{\displ@y
 \def\shoveright##1{##1\hfilneg\hskip\multlinegap@}%
 \def\shoveleft##1{\setboxz@h{$\m@th\displaystyle{}##1$}%
  \setbox@ne\hbox{$\m@th\displaystyle##1$}%
  \hfilneg
  \iftagin@
   \ifdim\ltwidth@>\z@\hskip\ltwidth@\hskip\multlinetaggap@\fi
  \else\hskip\multlinegap@\fi\hskip.5\wd@ne\hskip-.5\wdz@##1}
  \halign\bgroup\Let@\hbox to\displaywidth
   {\strut@$\m@th\displaystyle\hfil{}##\hfil$}\crcr
   \hfilneg                                                                 
   \iftagin@                                                                
    \ifdim\ltwidth@>\z@                                                     
     \box\mtagbox@\hskip\multlinetaggap@                                    
    \else
     \rlap{\vbox{\normalbaselines\hbox{\strut@\box\mtagbox@}%
     \vbox to\mlineht@{}}}\fi                                               
   \else\hskip\multlinegap@\fi}                                             
\def\multline@r#1$${\iftagin@\DN@{\rmultline@@#1$$}\else
 \DN@{\setbox\mtagbox@\null\ltwidth@\z@\rtwidth@\z@
  \rmultline@@@#1$$}\fi\next@}
\def\rmultline@@#1\endmultline\tag#2$${\ltwidth@\z@
 \setbox\mtagbox@\hbox{\maketag@#2\maketag@}%
 \rmmeasure@#1\endmultline\dimen@\mwidth@\advance\dimen@\wd\mtagbox@
 \advance\dimen@\multlinetaggap@
 \ifdim\dimen@>\displaywidth\rtwidth@\z@\else\rtwidth@\wd\mtagbox@\fi
 \rmultline@@@#1\endmultline$$}
\def\rmultline@@@{\displ@y
 \def\shoveright##1{##1\hfilneg\iftagin@\ifdim\rtwidth@>\z@
  \hskip\rtwidth@\hskip\multlinetaggap@\fi\else\hskip\multlinegap@\fi}%
 \def\shoveleft##1{\setboxz@h{$\m@th\displaystyle{}##1$}%
  \setbox@ne\hbox{$\m@th\displaystyle##1$}%
  \hfilneg\hskip\multlinegap@\hskip.5\wd@ne\hskip-.5\wdz@##1}%
 \halign\bgroup\Let@\hbox to\displaywidth
  {\strut@$\m@th\displaystyle\hfil{}##\hfil$}\crcr
 \hfilneg\hskip\multlinegap@}
\def\endmultline{\iftagsleft@\expandafter\lendmultline@\else
 \expandafter\rendmultline@\fi}
\def\lendmultline@{\hfilneg\hskip\multlinegap@\crcr\egroup}
\def\rendmultline@{\iftagin@                                                
 \ifdim\rtwidth@>\z@                                                        
  \hskip\multlinetaggap@\box\mtagbox@                                       
 \else\llap{\vtop{\normalbaselines\null\hbox{\strut@\box\mtagbox@}}}\fi     
 \else\hskip\multlinegap@\fi                                                
 \hfilneg\crcr\egroup}
\def\bmod{\mskip-\medmuskip\mkern5mu\mathbin{\fam\z@ mod}\penalty900
 \mkern5mu\mskip-\medmuskip}
\def\pmod#1{\allowbreak\ifinner\mkern8mu\else\mkern18mu\fi
 ({\fam\z@ mod}\,\,#1)}
\def\pod#1{\allowbreak\ifinner\mkern8mu\else\mkern18mu\fi(#1)}
\def\mod#1{\allowbreak\ifinner\mkern12mu\else\mkern18mu\fi{\fam\z@ mod}\,\,#1}
\message{continued fractions,}
\newcount\cfraccount@
\def\cfrac{\bgroup\bgroup\advance\cfraccount@\@ne\strut
 \iffalse{\fi\def\\{\over\displaystyle}\iffalse}\fi}
\def\lcfrac{\bgroup\bgroup\advance\cfraccount@\@ne\strut
 \iffalse{\fi\def\\{\hfill\over\displaystyle}\iffalse}\fi}
\def\rcfrac{\bgroup\bgroup\advance\cfraccount@\@ne\strut\hfill
 \iffalse{\fi\def\\{\over\displaystyle}\iffalse}\fi}
\def\gloop@#1\repeat{\gdef\body{#1}\iterate}
\def\endcfrac{\gloop@\ifnum\cfraccount@>\z@\global\advance\cfraccount@\m@ne
 \egroup\hskip-\nulldelimiterspace\egroup\repeat}
\message{compound symbols,}
\def\binrel@#1{\setboxz@h{\thinmuskip0mu
  \medmuskip\m@ne mu\thickmuskip\@ne mu$#1\m@th$}%
 \setbox@ne\hbox{\thinmuskip0mu\medmuskip\m@ne mu\thickmuskip
  \@ne mu${}#1{}\m@th$}%
 \setbox\tw@\hbox{\hskip\wd@ne\hskip-\wdz@}}
\def\overset#1\to#2{\binrel@{#2}\ifdim\wd\tw@<\z@
 \mathbin{\mathop{\kern\z@#2}\limits^{#1}}\else\ifdim\wd\tw@>\z@
 \mathrel{\mathop{\kern\z@#2}\limits^{#1}}\else
 {\mathop{\kern\z@#2}\limits^{#1}}{}\fi\fi}
\def\underset#1\to#2{\binrel@{#2}\ifdim\wd\tw@<\z@
 \mathbin{\mathop{\kern\z@#2}\limits_{#1}}\else\ifdim\wd\tw@>\z@
 \mathrel{\mathop{\kern\z@#2}\limits_{#1}}\else
 {\mathop{\kern\z@#2}\limits_{#1}}{}\fi\fi}
\def\oversetbrace#1\to#2{\overbrace{#2}^{#1}}
\def\undersetbrace#1\to#2{\underbrace{#2}_{#1}}
\def\sideset#1\and#2\to#3{%
 \setbox@ne\hbox{$\dsize{\vphantom{#3}}#1{#3}\m@th$}%
 \setbox\tw@\hbox{$\dsize{#3}#2\m@th$}%
 \hskip\wd@ne\hskip-\wd\tw@\mathop{\hskip\wd\tw@\hskip-\wd@ne
  {\vphantom{#3}}#1{#3}#2}}
\def\rightarrowfill@#1{\setboxz@h{$#1-\m@th$}\ht\z@\z@
  $#1\m@th\copy\z@\mkern-6mu\cleaders
  \hbox{$#1\mkern-2mu\box\z@\mkern-2mu$}\hfill
  \mkern-6mu\mathord\rightarrow$}
\def\leftarrowfill@#1{\setboxz@h{$#1-\m@th$}\ht\z@\z@
  $#1\m@th\mathord\leftarrow\mkern-6mu\cleaders
  \hbox{$#1\mkern-2mu\copy\z@\mkern-2mu$}\hfill
  \mkern-6mu\box\z@$}
\def\leftrightarrowfill@#1{\setboxz@h{$#1-\m@th$}\ht\z@\z@
  $#1\m@th\mathord\leftarrow\mkern-6mu\cleaders
  \hbox{$#1\mkern-2mu\box\z@\mkern-2mu$}\hfill
  \mkern-6mu\mathord\rightarrow$}
\def\overrightarrow{\mathpalette\overrightarrow@}
\def\overrightarrow@#1#2{\vbox{\ialign{##\crcr\rightarrowfill@#1\crcr
 \noalign{\kern-\ex@\nointerlineskip}$\m@th\hfil#1#2\hfil$\crcr}}}

\def\overleftarrow{\mathpalette\overleftarrow@}
\def\overleftarrow@#1#2{\vbox{\ialign{##\crcr\leftarrowfill@#1\crcr
 \noalign{\kern-\ex@\nointerlineskip}$\m@th\hfil#1#2\hfil$\crcr}}}
\def\overleftrightarrow{\mathpalette\overleftrightarrow@}
\def\overleftrightarrow@#1#2{\vbox{\ialign{##\crcr\leftrightarrowfill@#1\crcr
 \noalign{\kern-\ex@\nointerlineskip}$\m@th\hfil#1#2\hfil$\crcr}}}
\def\underrightarrow{\mathpalette\underrightarrow@}
\def\underrightarrow@#1#2{\vtop{\ialign{##\crcr$\m@th\hfil#1#2\hfil$\crcr
 \noalign{\nointerlineskip}\rightarrowfill@#1\crcr}}}

\def\underleftarrow{\mathpalette\underleftarrow@}
\def\underleftarrow@#1#2{\vtop{\ialign{##\crcr$\m@th\hfil#1#2\hfil$\crcr
 \noalign{\nointerlineskip}\leftarrowfill@#1\crcr}}}
\def\underleftrightarrow{\mathpalette\underleftrightarrow@}
\def\underleftrightarrow@#1#2{\vtop{\ialign{##\crcr$\m@th\hfil#1#2\hfil$\crcr
 \noalign{\nointerlineskip}\leftrightarrowfill@#1\crcr}}}
\message{various kinds of dots,}
\let\DOTSI\relax
\let\DOTSB\relax

\newif\ifmath@
{\uccode`7=`\\ \uccode`8=`m \uccode`9=`a \uccode`0=`t \uccode`!=`h
 \uppercase{\gdef\math@#1#2#3#4#5#6\math@{\global\math@false\ifx 7#1\ifx 8#2%
 \ifx 9#3\ifx 0#4\ifx !#5\xdef\meaning@{#6}\global\math@true\fi\fi\fi\fi\fi}}}
\newif\ifmathch@
{\uccode`7=`c \uccode`8=`h \uccode`9=`\"
 \uppercase{\gdef\mathch@#1#2#3#4#5#6\mathch@{\global\mathch@false
  \ifx 7#1\ifx 8#2\ifx 9#5\global\mathch@true\xdef\meaning@{9#6}\fi\fi\fi}}}
\newcount\classnum@
\def\getmathch@#1.#2\getmathch@{\classnum@#1 \divide\classnum@4096
 \ifcase\number\classnum@\or\or\gdef\thedots@{\dotsb@}\or
 \gdef\thedots@{\dotsb@}\fi}
\newif\ifmathbin@
{\uccode`4=`b \uccode`5=`i \uccode`6=`n
 \uppercase{\gdef\mathbin@#1#2#3{\relaxnext@
  \DNii@##1\mathbin@{\ifx\space@\next\global\mathbin@true\fi}%
 \global\mathbin@false\DN@##1\mathbin@{}%
 \ifx 4#1\ifx 5#2\ifx 6#3\DN@{\FN@\nextii@}\fi\fi\fi\next@}}}
\newif\ifmathrel@
{\uccode`4=`r \uccode`5=`e \uccode`6=`l
 \uppercase{\gdef\mathrel@#1#2#3{\relaxnext@
  \DNii@##1\mathrel@{\ifx\space@\next\global\mathrel@true\fi}%
 \global\mathrel@false\DN@##1\mathrel@{}%
 \ifx 4#1\ifx 5#2\ifx 6#3\DN@{\FN@\nextii@}\fi\fi\fi\next@}}}
\newif\ifmacro@
{\uccode`5=`m \uccode`6=`a \uccode`7=`c
 \uppercase{\gdef\macro@#1#2#3#4\macro@{\global\macro@false
  \ifx 5#1\ifx 6#2\ifx 7#3\global\macro@true
  \xdef\meaning@{\macro@@#4\macro@@}\fi\fi\fi}}}
\def\macro@@#1->#2\macro@@{#2}
\newif\ifDOTS@
\newcount\DOTSCASE@
{\uccode`6=`\\ \uccode`7=`D \uccode`8=`O \uccode`9=`T \uccode`0=`S
 \uppercase{\gdef\DOTS@#1#2#3#4#5{\global\DOTS@false\DN@##1\DOTS@{}%
  \ifx 6#1\ifx 7#2\ifx 8#3\ifx 9#4\ifx 0#5\let\next@\DOTS@@\fi\fi\fi\fi\fi
  \next@}}}
{\uccode`3=`B \uccode`4=`I \uccode`5=`X
 \uppercase{\gdef\DOTS@@#1{\relaxnext@
  \DNii@##1\DOTS@{\ifx\space@\next\global\DOTS@true\fi}%
  \DN@{\FN@\nextii@}%
  \ifx 3#1\global\DOTSCASE@\z@\else
  \ifx 4#1\global\DOTSCASE@\@ne\else
  \ifx 5#1\global\DOTSCASE@\tw@\else\DN@##1\DOTS@{}%
  \fi\fi\fi\next@}}}
\newif\ifnot@
{\uccode`5=`\\ \uccode`6=`n \uccode`7=`o \uccode`8=`t
 \uppercase{\gdef\not@#1#2#3#4{\relaxnext@
  \DNii@##1\not@{\ifx\space@\next\global\not@true\fi}%
 \global\not@false\DN@##1\not@{}%
 \ifx 5#1\ifx 6#2\ifx 7#3\ifx 8#4\DN@{\FN@\nextii@}\fi\fi\fi
 \fi\next@}}}
\newif\ifkeybin@
\def\keybin@{\keybin@true
 \ifx\next+\else\ifx\next=\else\ifx\next<\else\ifx\next>\else\ifx\next-\else
 \ifx\next*\else\ifx\next:\else\keybin@false\fi\fi\fi\fi\fi\fi\fi}
\def\dots{\RIfM@\expandafter\mdots@\else\expandafter\tdots@\fi}
\def\tdots@{\unskip\relaxnext@
 \DN@{$\m@th\mathinner{\ldotp\ldotp\ldotp}\,
   \ifx\next,\,$\else\ifx\next.\,$\else\ifx\next;\,$\else\ifx\next:\,$\else
   \ifx\next?\,$\else\ifx\next!\,$\else$ \fi\fi\fi\fi\fi\fi}%
 \ \FN@\next@}
\def\mdots@{\FN@\mdots@@}
\def\mdots@@{\gdef\thedots@{\dotso@}
 \ifx\next\boldkey\gdef\thedots@\boldkey{\boldkeydots@}\else                
 \ifx\next\boldsymbol\gdef\thedots@\boldsymbol{\boldsymboldots@}\else       
 \ifx,\next\gdef\thedots@{\dotsc}
 \else\ifx\not\next\gdef\thedots@{\dotsb@}
 \else\keybin@
 \ifkeybin@\gdef\thedots@{\dotsb@}
 \else\xdef\meaning@{\meaning\next..........}\xdef\meaning@@{\meaning@}
  \expandafter\math@\meaning@\math@
  \ifmath@
   \expandafter\mathch@\meaning@\mathch@
   \ifmathch@\expandafter\getmathch@\meaning@\getmathch@\fi                 
  \else\expandafter\macro@\meaning@@\macro@                                 
  \ifmacro@                                                                
   \expandafter\not@\meaning@\not@\ifnot@\gdef\thedots@{\dotsb@}
  \else\expandafter\DOTS@\meaning@\DOTS@
  \ifDOTS@
   \ifcase\number\DOTSCASE@\gdef\thedots@{\dotsb@}%
    \or\gdef\thedots@{\dotsi}\else\fi                                      
  \else\expandafter\math@\meaning@\math@                                   
  \ifmath@\expandafter\mathbin@\meaning@\mathbin@
  \ifmathbin@\gdef\thedots@{\dotsb@}
  \else\expandafter\mathrel@\meaning@\mathrel@
  \ifmathrel@\gdef\thedots@{\dotsb@}
  \fi\fi\fi\fi\fi\fi\fi\fi\fi\fi\fi\fi
 \thedots@}
\def\plainldots@{\mathinner{\ldotp\ldotp\ldotp}}
\def\plaincdots@{\mathinner{\cdotp\cdotp\cdotp}}
\def\dotsi{\!\plaincdots@}
\let\dotsb@\plaincdots@
\newif\ifextra@
\newif\ifrightdelim@
\def\rightdelim@{\global\rightdelim@true                                    
 \ifx\next)\else                                                            
 \ifx\next]\else
 \ifx\next\rbrack\else
 \ifx\next\}\else
 \ifx\next\rbrace\else
 \ifx\next\rangle\else
 \ifx\next\rceil\else
 \ifx\next\rfloor\else
 \ifx\next\rgroup\else
 \ifx\next\rmoustache\else
 \ifx\next\right\else
 \ifx\next\bigr\else
 \ifx\next\biggr\else
 \ifx\next\Bigr\else                                                        
 \ifx\next\Biggr\else\global\rightdelim@false
 \fi\fi\fi\fi\fi\fi\fi\fi\fi\fi\fi\fi\fi\fi\fi}
\def\extra@{%
 \global\extra@false\rightdelim@\ifrightdelim@\global\extra@true            
 \else\ifx\next$\global\extra@true                                          
 \else\xdef\meaning@{\meaning\next..........}
 \expandafter\macro@\meaning@\macro@\ifmacro@                               
 \expandafter\DOTS@\meaning@\DOTS@
 \ifDOTS@
 \ifnum\DOTSCASE@=\tw@\global\extra@true                                    
 \fi\fi\fi\fi\fi}
\newif\ifbold@
\def\dotso@{\relaxnext@
 \ifbold@
  \let\next\delayed@
  \DNii@{\extra@\plainldots@\ifextra@\,\fi}%
 \else
  \DNii@{\DN@{\extra@\plainldots@\ifextra@\,\fi}\FN@\next@}%
 \fi
 \nextii@}
\def\extrap@#1{%
 \ifx\next,\DN@{#1\,}\else
 \ifx\next;\DN@{#1\,}\else
 \ifx\next.\DN@{#1\,}\else\extra@
 \ifextra@\DN@{#1\,}\else
 \let\next@#1\fi\fi\fi\fi\next@}
\def\ldots{\DN@{\extrap@\plainldots@}%
 \FN@\next@}
\def\cdots{\DN@{\extrap@\plaincdots@}%
 \FN@\next@}

\def\dotsc{\relaxnext@
 \DN@{\ifx\next;\plainldots@\,\else
  \ifx\next.\plainldots@\,\else\extra@\plainldots@
  \ifextra@\,\fi\fi\fi}%
 \FN@\next@}
\def\cdot{\mathchar"2201 }

\message{special superscripts,}
\def\dddot#1{{\mathop{#1}\limits^{\vbox to-1.4\ex@{\kern-\tw@\ex@
 \hbox{\rm...}\vss}}}}
\def\ddddot#1{{\mathop{#1}\limits^{\vbox to-1.4\ex@{\kern-\tw@\ex@
 \hbox{\rm....}\vss}}}}
\def\sphat{^{\mathchoice{}{}%
 {\,\,\botsmash{\hbox{\lower4\ex@\hbox{$\m@th\widehat{\null}$}}}}%
 {\,\botsmash{\hbox{\lower3\ex@\hbox{$\m@th\hat{\null}$}}}}}}

\def\spacute{^{\!\botsmash{\hbox{\lower\@ne ex\hbox{\'{}}}}}}
\def\spgrave{^{\mathchoice{}{}{}{\!}%
 \botsmash{\hbox{\lower\@ne ex\hbox{\`{}}}}}}
\def\spdot{^{\hbox{\raise\ex@\hbox{\rm.}}}}
\def\spddot{^{\hbox{\raise\ex@\hbox{\rm..}}}}
\def\spdddot{^{\hbox{\raise\ex@\hbox{\rm...}}}}
\def\spddddot{^{\hbox{\raise\ex@\hbox{\rm....}}}}
\def\spbreve{^{\!\botsmash{\hbox{\lower4\ex@\hbox{\u{}}}}}}

\message{\string\text,}
\def\textonlyfont@#1#2{\def#1{\RIfM@
 \Err@{Use \string#1\space only in text}\else#2\fi}}
\textonlyfont@\rm\tenrm
\textonlyfont@\it\tenit
\textonlyfont@\sl\tensl
\textonlyfont@\bf\tenbf
\def\oldnos#1{\RIfM@{\mathcode`\,="013B \fam\@ne#1}\else
 \leavevmode\hbox{$\m@th\mathcode`\,="013B \fam\@ne#1$}\fi}
\def\text{\RIfM@\expandafter\text@\else\expandafter\text@@\fi}
\def\text@@#1{\leavevmode\hbox{#1}}
\def\mathhexbox@#1#2#3{\text{$\m@th\mathchar"#1#2#3$}}
\def\dag{{\mathhexbox@279}}
\def\ddag{{\mathhexbox@27A}}
\def\S{{\mathhexbox@278}}
\def\P{{\mathhexbox@27B}}
\newif\iffirstchoice@
\firstchoice@true
\def\text@#1{\mathchoice
 {\hbox{\everymath{\displaystyle}\def\textfonti{\the\textfont\@ne}%
  \def\textfontii{\the\textfont\tw@}\textdef@@ T#1}}
 {\hbox{\firstchoice@false
  \everymath{\textstyle}\def\textfonti{\the\textfont\@ne}%
  \def\textfontii{\the\textfont\tw@}\textdef@@ T#1}}
 {\hbox{\firstchoice@false
  \everymath{\scriptstyle}\def\textfonti{\the\scriptfont\@ne}%
  \def\textfontii{\the\scriptfont\tw@}\textdef@@ S\rm#1}}
 {\hbox{\firstchoice@false
  \everymath{\scriptscriptstyle}\def\textfonti
  {\the\scriptscriptfont\@ne}%
  \def\textfontii{\the\scriptscriptfont\tw@}\textdef@@ s\rm#1}}}
\def\textdef@@#1{\textdef@#1\rm\textdef@#1\bf\textdef@#1\sl\textdef@#1\it}
\def\rmfam{0}
\def\textdef@#1#2{%
 \DN@{\csname\expandafter\eat@\string#2fam\endcsname}%
 \if S#1\edef#2{\the\scriptfont\next@\relax}%
 \else\if s#1\edef#2{\the\scriptscriptfont\next@\relax}%
 \else\edef#2{\the\textfont\next@\relax}\fi\fi}
\scriptfont\itfam\tenit \scriptscriptfont\itfam\tenit
\scriptfont\slfam\tensl \scriptscriptfont\slfam\tensl
\newif\iftopfolded@
\newif\ifbotfolded@
\def\topfoldedtext{\topfolded@true\botfolded@false\foldedtext@}
\def\botfoldedtext{\botfolded@true\topfolded@false\foldedtext@}
\def\foldedtext{\topfolded@false\botfolded@false\foldedtext@}
\Invalid@\foldedwidth
\def\foldedtext@{\relaxnext@
 \DN@{\ifx\next\foldedwidth\let\next@\nextii@\else
  \DN@{\nextii@\foldedwidth{.3\hsize}}\fi\next@}%
 \DNii@\foldedwidth##1##2{\setbox\z@\vbox
  {\normalbaselines\hsize##1\relax
  \tolerance1600 \noindent\ignorespaces##2}\ifbotfolded@\boxz@\else
  \iftopfolded@\vtop{\unvbox\z@}\else\vcenter{\boxz@}\fi\fi}%
 \FN@\next@}
\message{math font commands,}
\def\bold{\RIfM@\expandafter\bold@\else
 \expandafter\nonmatherr@\expandafter\bold\fi}
\def\bold@#1{{\bold@@{#1}}}
\def\bold@@#1{\fam\bffam\relax#1}
\def\slanted{\RIfM@\expandafter\slanted@\else
 \expandafter\nonmatherr@\expandafter\slanted\fi}
\def\slanted@#1{{\slanted@@{#1}}}
\def\slanted@@#1{\fam\slfam\relax#1}
\def\roman{\RIfM@\expandafter\roman@\else
 \expandafter\nonmatherr@\expandafter\roman\fi}
\def\roman@#1{{\roman@@{#1}}}
\def\roman@@#1{\fam\rmfam\relax#1}
\def\italic{\RIfM@\expandafter\italic@\else
 \expandafter\nonmatherr@\expandafter\italic\fi}
\def\italic@#1{{\italic@@{#1}}}
\def\italic@@#1{\fam\itfam\relax#1}
\def\Cal{\RIfM@\expandafter\Cal@\else
 \expandafter\nonmatherr@\expandafter\Cal\fi}
\def\Cal@#1{{\Cal@@{#1}}}
\def\Cal@@#1{\noaccents@\fam\tw@#1}
\mathchardef\Gamma="0000
\mathchardef\Delta="0001
\mathchardef\Theta="0002
\mathchardef\Lambda="0003
\mathchardef\Xi="0004
\mathchardef\Pi="0005
\mathchardef\Sigma="0006
\mathchardef\Upsilon="0007
\mathchardef\Phi="0008
\mathchardef\Psi="0009
\mathchardef\Omega="000A
\mathchardef\varGamma="0100
\mathchardef\varDelta="0101
\mathchardef\varTheta="0102
\mathchardef\varLambda="0103
\mathchardef\varXi="0104
\mathchardef\varPi="0105
\mathchardef\varSigma="0106
\mathchardef\varUpsilon="0107
\mathchardef\varPhi="0108
\mathchardef\varPsi="0109
\mathchardef\varOmega="010A
\let\alloc@@\alloc@
\def\hexnumber@#1{\ifcase#1 0\or 1\or 2\or 3\or 4\or 5\or 6\or 7\or 8\or
 9\or A\or B\or C\or D\or E\or F\fi}
\def\loadmsam{%
 \font@\tenmsa=msam10
 \font@\sevenmsa=msam7
 \font@\fivemsa=msam5
 \alloc@@8\fam\chardef\sixt@@n\msafam
 \textfont\msafam=\tenmsa
 \scriptfont\msafam=\sevenmsa
 \scriptscriptfont\msafam=\fivemsa
 \edef\next{\hexnumber@\msafam}%
 \mathchardef\dabar@"0\next39
 \edef\dashrightarrow{\mathrel{\dabar@\dabar@\mathchar"0\next4B}}%
 \edef\dashleftarrow{\mathrel{\mathchar"0\next4C\dabar@\dabar@}}%
 \let\dasharrow\dashrightarrow
 \edef\ulcorner{\delimiter"4\next70\next70 }%
 \edef\urcorner{\delimiter"5\next71\next71 }%
 \edef\llcorner{\delimiter"4\next78\next78 }%
 \edef\lrcorner{\delimiter"5\next79\next79 }%
 \edef\yen{{\noexpand\mathhexbox@\next55}}%
 \edef\checkmark{{\noexpand\mathhexbox@\next58}}%
 \edef\circledR{{\noexpand\mathhexbox@\next72}}%
 \edef\maltese{{\noexpand\mathhexbox@\next7A}}%
 \global\let\loadmsam\empty}%
\def\loadmsbm{%
 \font@\tenmsb=msbm10 \font@\sevenmsb=msbm7 \font@\fivemsb=msbm5
 \alloc@@8\fam\chardef\sixt@@n\msbfam
 \textfont\msbfam=\tenmsb
 \scriptfont\msbfam=\sevenmsb \scriptscriptfont\msbfam=\fivemsb
 \global\let\loadmsbm\empty
 }
\def\widehat#1{\ifx\undefined\msbfam \DN@{362}%
  \else \setboxz@h{$\m@th#1$}%
    \edef\next@{\ifdim\wdz@>\tw@ em%
        \hexnumber@\msbfam 5B%
      \else 362\fi}\fi
  \mathaccent"0\next@{#1}}
\def\widetilde#1{\ifx\undefined\msbfam \DN@{365}%
  \else \setboxz@h{$\m@th#1$}%
    \edef\next@{\ifdim\wdz@>\tw@ em%
        \hexnumber@\msbfam 5D%
      \else 365\fi}\fi
  \mathaccent"0\next@{#1}}
\message{\string\newsymbol,}
\def\newsymbol#1#2#3#4#5{\define#1{}%
  \count@#2\relax \advance\count@\m@ne 
 \ifcase\count@
   \ifx\undefined\msafam\loadmsam\fi \let\next@\msafam
 \or \ifx\undefined\msbfam\loadmsbm\fi \let\next@\msbfam
 \else  \Err@{\Invalid@@\string\newsymbol}\let\next@\tw@\fi
 \mathchardef#1="#3\hexnumber@\next@#4#5\space}

\def\Bbb{\RIfM@\expandafter\Bbb@\else
 \expandafter\nonmatherr@\expandafter\Bbb\fi}
\def\Bbb@#1{{\Bbb@@{#1}}}
\def\Bbb@@#1{\noaccents@\fam\msbfam\relax#1}
\message{bold Greek and bold symbols,}
\def\loadbold{%
 \font@\tencmmib=cmmib10 \font@\sevencmmib=cmmib7 \font@\fivecmmib=cmmib5
 \skewchar\tencmmib'177 \skewchar\sevencmmib'177 \skewchar\fivecmmib'177
 \alloc@@8\fam\chardef\sixt@@n\cmmibfam
 \textfont\cmmibfam\tencmmib
 \scriptfont\cmmibfam\sevencmmib \scriptscriptfont\cmmibfam\fivecmmib
 \font@\tencmbsy=cmbsy10 \font@\sevencmbsy=cmbsy7 \font@\fivecmbsy=cmbsy5
 \skewchar\tencmbsy'60 \skewchar\sevencmbsy'60 \skewchar\fivecmbsy'60
 \alloc@@8\fam\chardef\sixt@@n\cmbsyfam
 \textfont\cmbsyfam\tencmbsy
 \scriptfont\cmbsyfam\sevencmbsy \scriptscriptfont\cmbsyfam\fivecmbsy
 \let\loadbold\empty
}
\def\boldnotloaded#1{\Err@{\ifcase#1\or First\else Second\fi
       bold symbol font not loaded}}
\def\mathchari@#1#2#3{\ifx\undefined\cmmibfam
    \boldnotloaded@\@ne
  \else\mathchar"#1\hexnumber@\cmmibfam#2#3\space \fi}
\def\mathcharii@#1#2#3{\ifx\undefined\cmbsyfam
    \boldnotloaded\tw@
  \else \mathchar"#1\hexnumber@\cmbsyfam#2#3\space\fi}
\edef\bffam@{\hexnumber@\bffam}
\def\boldkey#1{\ifcat\noexpand#1A%
  \ifx\undefined\cmmibfam \boldnotloaded\@ne
  \else {\fam\cmmibfam#1}\fi
 \else
 \ifx#1!\mathchar"5\bffam@21 \else
 \ifx#1(\mathchar"4\bffam@28 \else\ifx#1)\mathchar"5\bffam@29 \else
 \ifx#1+\mathchar"2\bffam@2B \else\ifx#1:\mathchar"3\bffam@3A \else
 \ifx#1;\mathchar"6\bffam@3B \else\ifx#1=\mathchar"3\bffam@3D \else
 \ifx#1?\mathchar"5\bffam@3F \else\ifx#1[\mathchar"4\bffam@5B \else
 \ifx#1]\mathchar"5\bffam@5D \else
 \ifx#1,\mathchari@63B \else
 \ifx#1-\mathcharii@200 \else
 \ifx#1.\mathchari@03A \else
 \ifx#1/\mathchari@03D \else
 \ifx#1<\mathchari@33C \else
 \ifx#1>\mathchari@33E \else
 \ifx#1*\mathcharii@203 \else
 \ifx#1|\mathcharii@06A \else
 \ifx#10\bold0\else\ifx#11\bold1\else\ifx#12\bold2\else\ifx#13\bold3\else
 \ifx#14\bold4\else\ifx#15\bold5\else\ifx#16\bold6\else\ifx#17\bold7\else
 \ifx#18\bold8\else\ifx#19\bold9\else
  \Err@{\string\boldkey\space can't be used with #1}%
 \fi\fi\fi\fi\fi\fi\fi\fi\fi\fi\fi\fi\fi\fi\fi
 \fi\fi\fi\fi\fi\fi\fi\fi\fi\fi\fi\fi\fi\fi}
\def\boldsymbol#1{%
 \DN@{\Err@{You can't use \string\boldsymbol\space with \string#1}#1}%
 \ifcat\noexpand#1A%
   \let\next@\relax
   \ifx\undefined\cmmibfam \boldnotloaded\@ne
   \else {\fam\cmmibfam#1}\fi
 \else
  \xdef\meaning@{\meaning#1.........}%
  \expandafter\math@\meaning@\math@
  \ifmath@
   \expandafter\mathch@\meaning@\mathch@
   \ifmathch@
    \expandafter\boldsymbol@@\meaning@\boldsymbol@@
   \fi
  \else
   \expandafter\macro@\meaning@\macro@
   \expandafter\delim@\meaning@\delim@
   \ifdelim@
    \expandafter\delim@@\meaning@\delim@@
   \else
    \boldsymbol@{#1}%
   \fi
  \fi
 \fi
 \next@}
\def\mathhexboxii@#1#2{\ifx\undefined\cmbsyfam
    \boldnotloaded\tw@
  \else \mathhexbox@{\hexnumber@\cmbsyfam}{#1}{#2}\fi}
\def\boldsymbol@#1{\let\next@\relax\let\next#1%
 \ifx\next\cdot\mathcharii@201 \else
 \ifx\next\prime{{\null\mathcharii@030 \null}}\else
 \ifx\next\lbrack\mathchar"4\bffam@5B \else
 \ifx\next\rbrack\mathchar"5\bffam@5D \else
 \ifx\next\{\mathcharii@466 \else
 \ifx\next\lbrace\mathcharii@466 \else
 \ifx\next\}\mathcharii@567 \else
 \ifx\next\rbrace\mathcharii@567 \else
 \ifx\next\surd{{\mathcharii@170}}\else
 \ifx\next\S{{\mathhexboxii@78}}\else
 \ifx\next\P{{\mathhexboxii@7B}}\else
 \ifx\next\dag{{\mathhexboxii@79}}\else
 \ifx\next\ddag{{\mathhexboxii@7A}}\else
 \DN@{\Err@{You can't use \string\boldsymbol\space with \string#1}#1}%
 \fi\fi\fi\fi\fi\fi\fi\fi\fi\fi\fi\fi\fi}
\def\boldsymbol@@#1.#2\boldsymbol@@{\classnum@#1 \count@@@\classnum@        
 \divide\classnum@4096 \count@\classnum@                                    
 \multiply\count@4096 \advance\count@@@-\count@ \count@@\count@@@           
 \divide\count@@@\@cclvi \count@\count@@                                    
 \multiply\count@@@\@cclvi \advance\count@@-\count@@@                       
 \divide\count@@@\@cclvi                                                    
 \multiply\classnum@4096 \advance\classnum@\count@@                         
 \ifnum\count@@@=\z@                                                        
  \count@"\bffam@ \multiply\count@\@cclvi
  \advance\classnum@\count@
  \DN@{\mathchar\number\classnum@}%
 \else
  \ifnum\count@@@=\@ne                                                      
   \ifx\undefined\cmmibfam \DN@{\boldnotloaded\@ne}%
   \else \count@\cmmibfam \multiply\count@\@cclvi
     \advance\classnum@\count@
     \DN@{\mathchar\number\classnum@}\fi
  \else
   \ifnum\count@@@=\tw@                                                    
     \ifx\undefined\cmbsyfam
       \DN@{\boldnotloaded\tw@}%
     \else
       \count@\cmbsyfam \multiply\count@\@cclvi
       \advance\classnum@\count@
       \DN@{\mathchar\number\classnum@}%
     \fi
  \fi
 \fi
\fi}
\newif\ifdelim@
\newcount\delimcount@
{\uccode`6=`\\ \uccode`7=`d \uccode`8=`e \uccode`9=`l
 \uppercase{\gdef\delim@#1#2#3#4#5\delim@
  {\delim@false\ifx 6#1\ifx 7#2\ifx 8#3\ifx 9#4\delim@true
   \xdef\meaning@{#5}\fi\fi\fi\fi}}}
\def\delim@@#1"#2#3#4#5#6\delim@@{\if#32%
\let\next@\relax
 \ifx\undefined\cmbsyfam \boldnotloaded\@ne
 \else \mathcharii@#2#4#5\space \fi\fi}
\def\vert{\delimiter"026A30C }
\def\Vert{\delimiter"026B30D }
\let\|\Vert
\def\backslash{\delimiter"026E30F }
\def\boldkeydots@#1{\bold@true\let\next=#1\let\delayed@=#1\mdots@@
 \boldkey#1\bold@false}  
\def\boldsymboldots@#1{\bold@true\let\next#1\let\delayed@#1\mdots@@
 \boldsymbol#1\bold@false}
\message{Euler fonts,}

\def\frak{\mathfont@\frak}

\def\loadmathfont#1{%
   \expandafter\font@\csname ten#1\endcsname=#110
   \expandafter\font@\csname seven#1\endcsname=#17
   \expandafter\font@\csname five#1\endcsname=#15
   \edef\next{\noexpand\alloc@@8\fam\chardef\sixt@@n
     \expandafter\noexpand\csname#1fam\endcsname}%
   \next
   \textfont\csname#1fam\endcsname \csname ten#1\endcsname
   \scriptfont\csname#1fam\endcsname \csname seven#1\endcsname
   \scriptscriptfont\csname#1fam\endcsname \csname five#1\endcsname
   \expandafter\def\csname #1\expandafter\endcsname\expandafter{%
      \expandafter\mathfont@\csname#1\endcsname}%
 \expandafter\gdef\csname load#1\endcsname{}%
}
\def\mathfont@#1{\RIfM@\expandafter\mathfont@@\expandafter#1\else
  \expandafter\nonmatherr@\expandafter#1\fi}
\def\mathfont@@#1#2{{\mathfont@@@#1{#2}}}
\def\mathfont@@@#1#2{\noaccents@
   \fam\csname\expandafter\eat@\string#1fam\endcsname
   \relax#2}
\message{math accents,}
\def\accentclass@{7}
\def\noaccents@{\def\accentclass@{0}}
\def\makeacc@#1#2{\def#1{\mathaccent"\accentclass@#2 }}
\makeacc@\hat{05E}
\makeacc@\check{014}
\makeacc@\tilde{07E}
\makeacc@\acute{013}
\makeacc@\grave{012}
\makeacc@\dot{05F}
\makeacc@\ddot{07F}
\makeacc@\breve{015}
\makeacc@\bar{016}

\newcount\skewcharcount@
\newcount\familycount@
\def\theskewchar@{\familycount@\@ne
 \global\skewcharcount@\the\skewchar\textfont\@ne                           
 \ifnum\fam>\m@ne\ifnum\fam<16
  \global\familycount@\the\fam\relax
  \global\skewcharcount@\the\skewchar\textfont\the\fam\relax\fi\fi          
 \ifnum\skewcharcount@>\m@ne
  \ifnum\skewcharcount@<128
  \multiply\familycount@256
  \global\advance\skewcharcount@\familycount@
  \global\advance\skewcharcount@28672
  \mathchar\skewcharcount@\else
  \global\skewcharcount@\m@ne\fi\else
 \global\skewcharcount@\m@ne\fi}                                            
\newcount\pointcount@
\def\getpoints@#1.#2\getpoints@{\pointcount@#1 }
\newdimen\accentdimen@
\newcount\accentmu@
\def\dimentomu@{\multiply\accentdimen@ 100
 \expandafter\getpoints@\the\accentdimen@\getpoints@
 \multiply\pointcount@18
 \divide\pointcount@\@m
 \global\accentmu@\pointcount@}
\def\Makeacc@#1#2{\def#1{\RIfM@\DN@{\mathaccent@
 {"\accentclass@#2 }}\else\DN@{\nonmatherr@{#1}}\fi\next@}}
\def\unbracefonts@{\let\Cal@\Cal@@\let\roman@\roman@@\let\bold@\bold@@
 \let\slanted@\slanted@@}
\def\mathaccent@#1#2{\ifnum\fam=\m@ne\xdef\thefam@{1}\else
 \xdef\thefam@{\the\fam}\fi                                                 
 \accentdimen@\z@                                                           
 \setboxz@h{\unbracefonts@$\m@th\fam\thefam@\relax#2$}
 \ifdim\accentdimen@=\z@\DN@{\mathaccent#1{#2}}
  \setbox@ne\hbox{\unbracefonts@$\m@th\fam\thefam@\relax#2\theskewchar@$}
  \setbox\tw@\hbox{$\m@th\ifnum\skewcharcount@=\m@ne\else
   \mathchar\skewcharcount@\fi$}
  \global\accentdimen@\wd@ne\global\advance\accentdimen@-\wdz@
  \global\advance\accentdimen@-\wd\tw@                                     
  \global\multiply\accentdimen@\tw@
  \dimentomu@\global\advance\accentmu@\@ne                                 
 \else\DN@{{\mathaccent#1{#2\mkern\accentmu@ mu}%
    \mkern-\accentmu@ mu}{}}\fi                                             
 \next@}\Makeacc@\Hat{05E}
\Makeacc@\Check{014}
\Makeacc@\Tilde{07E}
\Makeacc@\Acute{013}
\Makeacc@\Grave{012}
\Makeacc@\Dot{05F}
\Makeacc@\Ddot{07F}
\Makeacc@\Breve{015}
\Makeacc@\Bar{016}
\def\Vec{\RIfM@\DN@{\mathaccent@{"017E }}\else
 \DN@{\nonmatherr@\Vec}\fi\next@}
\def\accentedsymbol#1#2{\csname newbox\expandafter\endcsname
  \csname\expandafter\eat@\string#1@box\endcsname
 \expandafter\setbox\csname\expandafter\eat@
  \string#1@box\endcsname\hbox{$\m@th#2$}\define
  #1{\copy\csname\expandafter\eat@\string#1@box\endcsname{}}}
\message{roots,}
\def\sqrt#1{\radical"270370 {#1}}
\let\underline@\underline
\let\overline@\overline
\def\underline#1{\underline@{#1}}
\def\overline#1{\overline@{#1}}
\Invalid@\leftroot
\Invalid@\uproot
\newcount\uproot@
\newcount\leftroot@
\def\root{\relaxnext@
  \DN@{\ifx\next\uproot\let\next@\nextii@\else
   \ifx\next\leftroot\let\next@\nextiii@\else
   \let\next@\plainroot@\fi\fi\next@}%
  \DNii@\uproot##1{\uproot@##1\relax\FN@\nextiv@}%
  \def\nextiv@{\ifx\next\space@\DN@. {\FN@\nextv@}\else
   \DN@.{\FN@\nextv@}\fi\next@.}%
  \def\nextv@{\ifx\next\leftroot\let\next@\nextvi@\else
   \let\next@\plainroot@\fi\next@}%
  \def\nextvi@\leftroot##1{\leftroot@##1\relax\plainroot@}%
   \def\nextiii@\leftroot##1{\leftroot@##1\relax\FN@\nextvii@}%
  \def\nextvii@{\ifx\next\space@
   \DN@. {\FN@\nextviii@}\else
   \DN@.{\FN@\nextviii@}\fi\next@.}%
  \def\nextviii@{\ifx\next\uproot\let\next@\nextix@\else
   \let\next@\plainroot@\fi\next@}%
  \def\nextix@\uproot##1{\uproot@##1\relax\plainroot@}%
  \bgroup\uproot@\z@\leftroot@\z@\FN@\next@}
\def\plainroot@#1\of#2{\setbox\rootbox\hbox{$\m@th\scriptscriptstyle{#1}$}%
 \mathchoice{\r@@t\displaystyle{#2}}{\r@@t\textstyle{#2}}
 {\r@@t\scriptstyle{#2}}{\r@@t\scriptscriptstyle{#2}}\egroup}
\def\r@@t#1#2{\setboxz@h{$\m@th#1\sqrt{#2}$}%
 \dimen@\ht\z@\advance\dimen@-\dp\z@
 \setbox@ne\hbox{$\m@th#1\mskip\uproot@ mu$}\advance\dimen@ 1.667\wd@ne
 \mkern-\leftroot@ mu\mkern5mu\raise.6\dimen@\copy\rootbox
 \mkern-10mu\mkern\leftroot@ mu\boxz@}
\def\boxed#1{\setboxz@h{$\m@th\displaystyle{#1}$}\dimen@.4\ex@
 \advance\dimen@3\ex@\advance\dimen@\dp\z@
 \hbox{\lower\dimen@\hbox{%
 \vbox{\hrule height.4\ex@
 \hbox{\vrule width.4\ex@\hskip3\ex@\vbox{\vskip3\ex@\boxz@\vskip3\ex@}%
 \hskip3\ex@\vrule width.4\ex@}\hrule height.4\ex@}%
 }}}
\message{commutative diagrams,}
\let\ampersand@\relax
\newdimen\minaw@
\minaw@11.11128\ex@
\newdimen\minCDaw@
\minCDaw@2.5pc
\def\minCDarrowwidth#1{\RIfMIfI@\onlydmatherr@\minCDarrowwidth
 \else\minCDaw@#1\relax\fi\else\onlydmatherr@\minCDarrowwidth\fi}
\newif\ifCD@
\def\CD{\bgroup\vspace@\relax\let\ampersand@&\iffalse}\fi
 \CD@true\vcenter\bgroup\Let@\tabskip\z@skip\baselineskip20\ex@
 \lineskip3\ex@\lineskiplimit3\ex@\halign\bgroup
 &\hfill$\m@th##$\hfill\crcr}
\def\endCD{\crcr\egroup\egroup\egroup}
\newdimen\bigaw@
\atdef@>#1>#2>{\ampersand@                                                  
 \setboxz@h{$\m@th\ssize\;{#1}\;\;$}
 \setbox@ne\hbox{$\m@th\ssize\;{#2}\;\;$}
 \setbox\tw@\hbox{$\m@th#2$}
 \ifCD@\global\bigaw@\minCDaw@\else\global\bigaw@\minaw@\fi                 
 \ifdim\wdz@>\bigaw@\global\bigaw@\wdz@\fi
 \ifdim\wd@ne>\bigaw@\global\bigaw@\wd@ne\fi                                
 \ifCD@\enskip\fi                                                           
 \ifdim\wd\tw@>\z@
  \mathrel{\mathop{\hbox to\bigaw@{\rightarrowfill@\displaystyle}}%
    \limits^{#1}_{#2}}
 \else\mathrel{\mathop{\hbox to\bigaw@{\rightarrowfill@\displaystyle}}%
    \limits^{#1}}\fi                                                        
 \ifCD@\enskip\fi                                                          
 \ampersand@}                                                              
\atdef@<#1<#2<{\ampersand@\setboxz@h{$\m@th\ssize\;\;{#1}\;$}%
 \setbox@ne\hbox{$\m@th\ssize\;\;{#2}\;$}\setbox\tw@\hbox{$\m@th#2$}%
 \ifCD@\global\bigaw@\minCDaw@\else\global\bigaw@\minaw@\fi
 \ifdim\wdz@>\bigaw@\global\bigaw@\wdz@\fi
 \ifdim\wd@ne>\bigaw@\global\bigaw@\wd@ne\fi
 \ifCD@\enskip\fi
 \ifdim\wd\tw@>\z@
  \mathrel{\mathop{\hbox to\bigaw@{\leftarrowfill@\displaystyle}}%
       \limits^{#1}_{#2}}\else
  \mathrel{\mathop{\hbox to\bigaw@{\leftarrowfill@\displaystyle}}%
       \limits^{#1}}\fi
 \ifCD@\enskip\fi\ampersand@}
\begingroup
 \catcode`\~=\active \lccode`\~=`\@
 \lowercase{%
  \global\atdef@)#1)#2){~>#1>#2>}
  \global\atdef@(#1(#2({~<#1<#2<}}
\endgroup
\atdef@ A#1A#2A{\llap{$\m@th\vcenter{\hbox
 {$\ssize#1$}}$}\Big\uparrow\rlap{$\m@th\vcenter{\hbox{$\ssize#2$}}$}&&}
\atdef@ V#1V#2V{\llap{$\m@th\vcenter{\hbox
 {$\ssize#1$}}$}\Big\downarrow\rlap{$\m@th\vcenter{\hbox{$\ssize#2$}}$}&&}
\atdef@={&\enskip\mathrel
 {\vbox{\hrule width\minCDaw@\vskip3\ex@\hrule width
 \minCDaw@}}\enskip&}
\atdef@|{\Big\Vert&&}
\atdef@\vert{\Big\Vert&&}
\def\pretend#1\haswidth#2{\setboxz@h{$\m@th\scriptstyle{#2}$}\hbox
 to\wdz@{\hfill$\m@th\scriptstyle{#1}$\hfill}}
\message{poor man's bold,}
\def\pmb{\RIfM@\expandafter\mathpalette\expandafter\pmb@\else
 \expandafter\pmb@@\fi}
\def\pmb@@#1{\leavevmode\setboxz@h{#1}%
   \dimen@-\wdz@
   \kern-.5\ex@\copy\z@
   \kern\dimen@\kern.25\ex@\raise.4\ex@\copy\z@
   \kern\dimen@\kern.25\ex@\box\z@
}
\def\binrel@@#1{\ifdim\wd2<\z@\mathbin{#1}\else\ifdim\wd\tw@>\z@
 \mathrel{#1}\else{#1}\fi\fi}
\newdimen\pmbraise@
\def\pmb@#1#2{\setbox\thr@@\hbox{$\m@th#1{#2}$}%
 \setbox4\hbox{$\m@th#1\mkern.5mu$}\pmbraise@\wd4\relax
 \binrel@{#2}%
 \dimen@-\wd\thr@@
   \binrel@@{%
   \mkern-.8mu\copy\thr@@
   \kern\dimen@\mkern.4mu\raise\pmbraise@\copy\thr@@
   \kern\dimen@\mkern.4mu\box\thr@@
}}
\def\documentstyle#1{\W@{}\input #1.sty\relax}
\message{syntax check,}
\font\dummyft@=dummy
\fontdimen1 \dummyft@=\z@
\fontdimen2 \dummyft@=\z@
\fontdimen3 \dummyft@=\z@
\fontdimen4 \dummyft@=\z@
\fontdimen5 \dummyft@=\z@
\fontdimen6 \dummyft@=\z@
\fontdimen7 \dummyft@=\z@
\fontdimen8 \dummyft@=\z@
\fontdimen9 \dummyft@=\z@
\fontdimen10 \dummyft@=\z@
\fontdimen11 \dummyft@=\z@
\fontdimen12 \dummyft@=\z@
\fontdimen13 \dummyft@=\z@
\fontdimen14 \dummyft@=\z@
\fontdimen15 \dummyft@=\z@
\fontdimen16 \dummyft@=\z@
\fontdimen17 \dummyft@=\z@
\fontdimen18 \dummyft@=\z@
\fontdimen19 \dummyft@=\z@
\fontdimen20 \dummyft@=\z@
\fontdimen21 \dummyft@=\z@
\fontdimen22 \dummyft@=\z@
\def\fontlist@{\\{\tenrm}\\{\sevenrm}\\{\fiverm}\\{\teni}\\{\seveni}%
 \\{\fivei}\\{\tensy}\\{\sevensy}\\{\fivesy}\\{\tenex}\\{\tenbf}\\{\sevenbf}%
 \\{\fivebf}\\{\tensl}\\{\tenit}}
\def\font@#1=#2 {\rightappend@#1\to\fontlist@\font#1=#2 }
\def\dodummy@{{\def\\##1{\global\let##1\dummyft@}\fontlist@}}
\def\nopages@{\output{\setbox\z@\box\@cclv \deadcycles\z@}%
 \alloc@5\toks\toksdef\@cclvi\output}
\let\galleys\nopages@
\newif\ifsyntax@
\newcount\countxviii@
\def\syntax{\syntax@true\dodummy@\countxviii@\count18
 \loop\ifnum\countxviii@>\m@ne\textfont\countxviii@=\dummyft@
 \scriptfont\countxviii@=\dummyft@\scriptscriptfont\countxviii@=\dummyft@
 \advance\countxviii@\m@ne\repeat                                           
 \dummyft@\tracinglostchars\z@\nopages@\frenchspacing\hbadness\@M}
\def\first@#1#2\end{#1}
\def\printoptions{\W@{Do you want S(yntax check),
  G(alleys) or P(ages)?}%
 \message{Type S, G or P, followed by <return>: }%
 \begingroup 
 \endlinechar\m@ne 
 \read\m@ne to\ans@
 \edef\ans@{\uppercase{\def\noexpand\ans@{%
   \expandafter\first@\ans@ P\end}}}%
 \expandafter\endgroup\ans@
 \if\ans@ P
 \else \if\ans@ S\syntax
 \else \if\ans@ G\galleys
 \else\message{? Unknown option: \ans@; using the `pages' option.}%
 \fi\fi\fi}
\def\alloc@#1#2#3#4#5{\global\advance\count1#1by\@ne
 \ch@ck#1#4#2\allocationnumber=\count1#1
 \global#3#5=\allocationnumber
 \ifalloc@\wlog{\string#5=\string#2\the\allocationnumber}\fi}
\def\document{\def\alloclist@{}\def\fontlist@{}}
\let\enddocument\bye

\let\proclaim\undefined
\let\footnote\undefined
\let\=\undefined
\let\>\undefined

\catcode`\@=\active
\message{... finished}

\expandafter\ifx\csname mathdefs.tex\endcsname\relax
  \expandafter\gdef\csname mathdefs.tex\endcsname{}
\else \message{Hey!  Apparently you were trying to
  \string\input{mathdefs.tex} twice.   This does not make sense.} 
\errmessage{Please edit your file (probably \jobname.tex) and remove
any duplicate ``\string\input'' lines}\endinput\fi




\catcode`\X=12\catcode`\@=11

\def\n@wcount{\alloc@0\count\countdef\insc@unt}
\def\n@wwrite{\alloc@7\write\chardef\sixt@@n}
\def\n@wread{\alloc@6\read\chardef\sixt@@n}
\def\r@s@t{\relax}\def\v@idline{\par}\def\@mputate#1/{#1}
\def\l@c@l#1X{\firstpart.#1}\def\gl@b@l#1X{#1}\def\t@d@l#1X{{}}

\def\crossrefs#1{\ifx\all#1\let\tr@ce=\all\else\def\tr@ce{#1,}\fi
   \n@wwrite\cit@tionsout\openout\cit@tionsout=\jobname.cit 
   \write\cit@tionsout{\tr@ce}\expandafter\setfl@gs\tr@ce,}
\def\setfl@gs#1,{\def\@{#1}\ifx\@\empty\let\next=\relax
   \else\let\next=\setfl@gs\expandafter\xdef
   \csname#1tr@cetrue\endcsname{}\fi\next}
\def\m@ketag#1#2{\expandafter\n@wcount\csname#2tagno\endcsname
     \csname#2tagno\endcsname=0\let\tail=\all\xdef\all{\tail#2,}
   \ifx#1\l@c@l\let\tail=\r@s@t\xdef\r@s@t{\csname#2tagno\endcsname=0\tail}\fi
   \expandafter\gdef\csname#2cite\endcsname##1{\expandafter
     \ifx\csname#2tag##1\endcsname\relax?\else\csname#2tag##1\endcsname\fi
     \expandafter\ifx\csname#2tr@cetrue\endcsname\relax\else
     \write\cit@tionsout{#2tag ##1 cited on page \folio.}\fi}
   \expandafter\gdef\csname#2page\endcsname##1{\expandafter
     \ifx\csname#2page##1\endcsname\relax?\else\csname#2page##1\endcsname\fi
     \expandafter\ifx\csname#2tr@cetrue\endcsname\relax\else
     \write\cit@tionsout{#2tag ##1 cited on page \folio.}\fi}
   \expandafter\gdef\csname#2tag\endcsname##1{\expandafter
      \ifx\csname#2check##1\endcsname\relax
      \expandafter\xdef\csname#2check##1\endcsname{}%
      \else\immediate\write16{Warning: #2tag ##1 used more than once.}\fi
      \multit@g{#1}{#2}##1/X%
      \write\t@gsout{#2tag ##1 assigned number \csname#2tag##1\endcsname\space
      on page \number\count0.}%
   \csname#2tag##1\endcsname}}

\def\multit@g#1#2#3/#4X{\def\t@mp{#4}\ifx\t@mp\empty%
      \global\advance\csname#2tagno\endcsname by 1 
      \expandafter\xdef\csname#2tag#3\endcsname
      {#1\number\csname#2tagno\endcsnameX}%
   \else\expandafter\ifx\csname#2last#3\endcsname\relax
      \expandafter\n@wcount\csname#2last#3\endcsname
      \global\advance\csname#2tagno\endcsname by 1 
      \expandafter\xdef\csname#2tag#3\endcsname
      {#1\number\csname#2tagno\endcsnameX}
      \write\t@gsout{#2tag #3 assigned number \csname#2tag#3\endcsname\space
      on page \number\count0.}\fi
   \global\advance\csname#2last#3\endcsname by 1
   \def\t@mp{\expandafter\xdef\csname#2tag#3/}%
   \expandafter\t@mp\@mputate#4\endcsname
   {\csname#2tag#3\endcsname\lastpart{\csname#2last#3\endcsname}}\fi}
\def\t@gs#1{\def\all{}\m@ketag#1e\m@ketag#1s\m@ketag\t@d@l p
\let\realscite\scite
\let\realstag\stag
   \m@ketag\gl@b@l r \n@wread\t@gsin
   \openin\t@gsin=\jobname.tgs \re@der \closein\t@gsin
   \n@wwrite\t@gsout\openout\t@gsout=\jobname.tgs }
\outer\def\localtags{\t@gs\l@c@l}
\outer\def\globaltags{\t@gs\gl@b@l}
\outer\def\newlocaltag#1{\m@ketag\l@c@l{#1}}
\outer\def\newglobaltag#1{\m@ketag\gl@b@l{#1}}

\newif\ifpr@ 
\def\m@kecs #1tag #2 assigned number #3 on page #4.%
   {\expandafter\gdef\csname#1tag#2\endcsname{#3}
   \expandafter\gdef\csname#1page#2\endcsname{#4}
   \ifpr@\expandafter\xdef\csname#1check#2\endcsname{}\fi}
\def\re@der{\ifeof\t@gsin\let\next=\relax\else
   \read\t@gsin to\t@gline\ifx\t@gline\v@idline\else
   \expandafter\m@kecs \t@gline\fi\let \next=\re@der\fi\next}
\def\pretags#1{\pr@true\pret@gs#1,,}
\def\pret@gs#1,{\def\@{#1}\ifx\@\empty\let\n@xtfile=\relax
   \else\let\n@xtfile=\pret@gs \openin\t@gsin=#1.tgs \message{#1} \re@der 
   \closein\t@gsin\fi \n@xtfile}

\newcount\sectno\sectno=0\newcount\subsectno\subsectno=0
\newif\ifultr@local \def\ultralocal{\ultr@localtrue}
\def\firstpart{\number\sectno}
\def\lastpart#1{\ifcase#1 \or a\or b\or c\or d\or e\or f\or g\or h\or 
   i\or k\or l\or m\or n\or o\or p\or q\or r\or s\or t\or u\or v\or w\or 
   x\or y\or z \fi}

\def\resetall{\global\advance\sectno by 1\subsectno=0
   \gdef\firstpart{\number\sectno}\r@s@t}
\def\resetsub{\global\advance\subsectno by 1
   \gdef\firstpart{\number\sectno.\number\subsectno}\r@s@t}
\def\newsection#1\par{\resetall\vskip0pt plus.3\vsize\penalty-250
   \vskip0pt plus-.3\vsize\bigskip\bigskip
   \message{#1}\leftline{\bf#1}\nobreak\bigskip}
\def\subsection#1\par{\ifultr@local\resetsub\fi
   \vskip0pt plus.2\vsize\penalty-250\vskip0pt plus-.2\vsize
   \bigskip\smallskip\message{#1}\leftline{\bf#1}\nobreak\medskip}


\newdimen\marginshift

\newdimen\margindelta
\newdimen\marginmax
\newdimen\marginmin

\def\margininit{       
\marginmax=3 true cm                  
				      
\margindelta=0.1 true cm              
\marginmin=0.1true cm                 
\marginshift=\marginmin
}    

\def\t@gsjj#1,{\def\@{#1}\ifx\@\empty\let\next=\relax\else\let\next=\t@gsjj
   \def\@@{p}\ifx\@\@@\else
   \expandafter\gdef\csname#1cite\endcsname##1{\citejj{##1}}
   \expandafter\gdef\csname#1page\endcsname##1{?}
   \expandafter\gdef\csname#1tag\endcsname##1{\tagjj{##1}}\fi\fi\next}
\newif\ifshowstuffinmargin
\showstuffinmarginfalse
\def\jjtags{\ifx\shlhetal\relax 
  \else
\ifx\shlhetal\undefinedcontrolseq
\else
\showstuffinmargintrue
\ifx\all\relax\else\expandafter\t@gsjj\all,\fi\fi \fi
}

\def\tagjj#1{\realstag{#1}\mginpar{\zeigen{#1}}}
\def\citejj#1{\zeigen{#1}\mginpar{\rechnen{#1}}}

\def\rechnen#1{\expandafter\ifx\csname stag#1\endcsname\relax ??\else
                           \csname stag#1\endcsname\fi}

\newdimen\theight

\def\marginfont{\sevenrm}

\def\trymarginbox#1{\setbox0=\hbox{\marginfont\hskip\marginshift #1}%
		\global\marginshift\wd0 
		\global\advance\marginshift\margindelta}

\def \mginpar#1{%
\ifvmode\setbox0\hbox to \hsize{\hfill\rlap{\marginfont\quad#1}}%
\ht0 0cm
\dp0 0cm
\box0\vskip-\baselineskip
\else 
             \vadjust{\trymarginbox{#1}%
		\ifdim\marginshift>\marginmax \global\marginshift\marginmin
			\trymarginbox{#1}%
                \fi
             \theight=\ht0
             \advance\theight by \dp0    \advance\theight by \lineskip
             \kern -\theight \vbox to \theight{\rightline{\rlap{\box0}}%
\vss}}\fi}


\def\t@gsoff#1,{\def\@{#1}\ifx\@\empty\let\next=\relax\else\let\next=\t@gsoff
   \def\@@{p}\ifx\@\@@\else
   \expandafter\gdef\csname#1cite\endcsname##1{\zeigen{##1}}
   \expandafter\gdef\csname#1page\endcsname##1{?}
   \expandafter\gdef\csname#1tag\endcsname##1{\zeigen{##1}}\fi\fi\next}
\def\verbatimtags{\showstuffinmarginfalse
\ifx\all\relax\else\expandafter\t@gsoff\all,\fi}
\def\zeigen#1{\hbox{$\langle$}#1\hbox{$\rangle$}}
\def\margincite#1{\ifshowstuffinmargin\mginpar{\rechnen{#1}}\fi}

\def\(#1){\edef\dot@g{\ifmmode\ifinner(\hbox{\noexpand\etag{#1}})
   \else\noexpand\eqno(\hbox{\noexpand\etag{#1}})\fi
   \else(\noexpand\ecite{#1})\fi}\dot@g}

\newif\ifbr@ck
\def\eat#1{}
\def\[#1]{\br@cktrue[\br@cket#1'X]}
\def\br@cket#1'#2X{\def\temp{#2}\ifx\temp\empty\let\next\eat
   \else\let\next\br@cket\fi
   \ifbr@ck\br@ckfalse\br@ck@t#1,X\else\br@cktrue#1\fi\next#2X}
\def\br@ck@t#1,#2X{\def\temp{#2}\ifx\temp\empty\let\neext\eat
   \else\let\neext\br@ck@t\def\temp{,}\fi
   \def\teemp{#1}\ifx\teemp\empty\else\rcite{#1}\fi\temp\neext#2X}
\def\resetbr@cket{\gdef\[##1]{[\rtag{##1}]}}
\def\references{\resetbr@cket\newsection References\par}

\newtoks\symb@ls\newtoks\s@mb@ls\newtoks\p@gelist\n@wcount\ftn@mber
    \ftn@mber=1\newif\ifftn@mbers\ftn@mbersfalse\newif\ifbyp@ge\byp@gefalse
\def\defm@rk{\ifftn@mbers\n@mberm@rk\else\symb@lm@rk\fi}
\def\n@mberm@rk{\xdef\m@rk{{\the\ftn@mber}}%
    \global\advance\ftn@mber by 1 }
\def\rot@te#1{\let\temp=#1\global#1=\expandafter\r@t@te\the\temp,X}
\def\r@t@te#1,#2X{{#2#1}\xdef\m@rk{{#1}}}
\def\b@@st#1{{$^{#1}$}}\def\str@p#1{#1}
\def\symb@lm@rk{\ifbyp@ge\rot@te\p@gelist\ifnum\expandafter\str@p\m@rk=1 
    \s@mb@ls=\symb@ls\fi\write\f@nsout{\number\count0}\fi \rot@te\s@mb@ls}
\def\byp@ge{\byp@getrue\n@wwrite\f@nsin\openin\f@nsin=\jobname.fns 
    \n@wcount\currentp@ge\currentp@ge=0\p@gelist={0}
    \re@dfns\closein\f@nsin\rot@te\p@gelist
    \n@wread\f@nsout\openout\f@nsout=\jobname.fns }
\def\m@kelist#1X#2{{#1,#2}}
\def\re@dfns{\ifeof\f@nsin\let\next=\relax\else\read\f@nsin to \f@nline
    \ifx\f@nline\v@idline\else\let\t@mplist=\p@gelist
    \ifnum\currentp@ge=\f@nline
    \global\p@gelist=\expandafter\m@kelist\the\t@mplistX0
    \else\currentp@ge=\f@nline
    \global\p@gelist=\expandafter\m@kelist\the\t@mplistX1\fi\fi
    \let\next=\re@dfns\fi\next}
\def\symbols#1{\symb@ls={#1}\s@mb@ls=\symb@ls} 
\def\bigsymbol{\textstyle}
\symbols{\bigsymbol\ast,\dagger,\ddagger,\sharp,\flat,\natural,\star}
\def\ftnumbers{\ftn@mberstrue} \def\ftsymbols{\ftn@mbersfalse}
\def\paginal{\byp@ge} \def\resetftnumbers{\ftn@mber=1}
\def\ftnote#1{\defm@rk\expandafter\expandafter\expandafter\footnote
    \expandafter\b@@st\m@rk{#1}}

\long\def\jump#1\endjump{}
\def\ssum{\mathop{\lower .1em\hbox{$\textstyle\Sigma$}}\nolimits}

\def\qed{\nobreak\kern 1em \vrule height .5em width .5em depth 0em}
\def\newneq{\hbox{\rlap{\hbox to 1\wd9{\hss$=$\hss}}\raise .1em 
   \hbox to 1\wd9{\hss$\scriptscriptstyle/$\hss}}}
\def\subsetne{\setbox9 = \hbox{$\subset$}\mathrel{\hbox{\rlap
   {\lower .4em \newneq}\raise .13em \hbox{$\subset$}}}}
\def\supsetne{\setbox9 = \hbox{$\subset$}\mathrel{\hbox{\rlap
   {\lower .4em \newneq}\raise .13em \hbox{$\supset$}}}}

\def\vbar{\mathchoice{\vrule height6.3ptdepth-.5ptwidth.8pt\kern-.8pt}
   {\vrule height6.3ptdepth-.5ptwidth.8pt\kern-.8pt}
   {\vrule height4.1ptdepth-.35ptwidth.6pt\kern-.6pt}
   {\vrule height3.1ptdepth-.25ptwidth.5pt\kern-.5pt}}
\def\f@dge{\mathchoice{}{}{\mkern.5mu}{\mkern.8mu}}
\def\b@c#1#2{{\rm \mkern#2mu\vbar\mkern-#2mu#1}}
\def\b@b#1{{\rm I\mkern-3.5mu #1}}
\def\b@a#1#2{{\rm #1\mkern-#2mu\f@dge #1}}
\def\bb#1{{\count4=`#1 \advance\count4by-64 \ifcase\count4\or\b@a A{11.5}\or
   \b@b B\or\b@c C{5}\or\b@b D\or\b@b E\or\b@b F \or\b@c G{5}\or\b@b H\or
   \b@b I\or\b@c J{3}\or\b@b K\or\b@b L \or\b@b M\or\b@b N\or\b@c O{5} \or
   \b@b P\or\b@c Q{5}\or\b@b R\or\b@a S{8}\or\b@a T{10.5}\or\b@c U{5}\or
   \b@a V{12}\or\b@a W{16.5}\or\b@a X{11}\or\b@a Y{11.7}\or\b@a Z{7.5}\fi}}

\catcode`\X=11 \catcode`\@=12


\expandafter\ifx\csname citeadd.tex\endcsname\relax
\expandafter\gdef\csname citeadd.tex\endcsname{}
\else \message{Hey!  Apparently you were trying to
\string\input{citeadd.tex} twice.   This does not make sense.} 
\errmessage{Please edit your file (probably \jobname.tex) and remove
any duplicate ``\string\input'' lines}\endinput\fi

\def\sciteu{\sciteerror{undefined}}
\def\sciteuphantom{\complainaboutcitation{undefined}}

\def\sciteerror#1#2{{\mathortextbf{\scite{#2}}}\complainaboutcitation{#1}{#2}}
\def\mathortextbf#1{\hbox{\bf #1}}
\def\complainaboutcitation#1#2{%
\vadjust{\line{\llap{---$\!\!>$ }\qquad scite$\{$#2$\}$ #1\hfil}}}

\sectno=-1   
\localtags
\jjtags
\NoBlackBoxes
\define\mr{\medskip\roster}
\define\sn{\smallskip\noindent}
\define\mn{\medskip\noindent}
\define\bn{\bigskip\noindent}
\define\ub{\underbar}
\define\wilog{\text{without loss of generality}}
\define\ermn{\endroster\medskip\noindent}
\define\dbca{\dsize\bigcap}
\define\dbcu{\dsize\bigcup}
\define \nl{\newline}
\magnification=\magstep 1
\documentstyle{amsppt}

{    
\catcode`@11

\ifx\alicetwothousandloaded@\relax
  \endinput\else\global\let\alicetwothousandloaded@\relax\fi

\gdef\subjclass{\let\savedef@\subjclass
 \def\subjclass##1\endsubjclass{\let\subjclass\savedef@
   \toks@{\def\usualspace{{\rm\enspace}}\eightpoint}%
   \toks@@{##1\unskip.}%
   \edef\thesubjclass@{\the\toks@
     \frills@{{\noexpand\rm2000 {\noexpand\it Mathematics Subject
       Classification}.\noexpand\enspace}}%
     \the\toks@@}}%
  \nofrillscheck\subjclass}
} 


\expandafter\ifx\csname alice2jlem.tex\endcsname\relax
  \expandafter\xdef\csname alice2jlem.tex\endcsname{\the\catcode`@}
\else \message{Hey!  Apparently you were trying to
\string\input{alice2jlem.tex}  twice.   This does not make sense.}
\errmessage{Please edit your file (probably \jobname.tex) and remove
any duplicate ``\string\input'' lines}\endinput\fi

\expandafter\ifx\csname bib4plain.tex\endcsname\relax
  \expandafter\gdef\csname bib4plain.tex\endcsname{}
\else \message{Hey!  Apparently you were trying to \string\input
  bib4plain.tex twice.   This does not make sense.}
\errmessage{Please edit your file (probably \jobname.tex) and remove
any duplicate ``\string\input'' lines}\endinput\fi

\def\renewcommand{\newcommand}	       
\edef\cite{\the\catcode`@}%
\catcode`@ = 11
\let\@oldatcatcode = \cite
\chardef\@letter = 11
\chardef\@other = 12
%
%
%
%
\def\@innerdef#1#2{\edef#1{\expandafter\noexpand\csname #2\endcsname}}%
%
%
\@innerdef\@innernewcount{newcount}%
\@innerdef\@innernewdimen{newdimen}%
\@innerdef\@innernewif{newif}%
\@innerdef\@innernewwrite{newwrite}%
%
%
%
\def\@gobble#1{}%
%
%
%
\ifx\inputlineno\@undefined
   \let\@linenumber = \empty 
\else
   \def\@linenumber{\the\inputlineno:\space}%
\fi
%
%
%
\def\@futurenonspacelet#1{\def\cs{#1}%
   \afterassignment\@stepone\let\@nexttoken=
}%
\begingroup 
\def\\{\global\let\@stoken= }%
\\ 
\endgroup
\def\@stepone{\expandafter\futurelet\cs\@steptwo}%
\def\@steptwo{\expandafter\ifx\cs\@stoken\let\@@next=\@stepthree
   \else\let\@@next=\@nexttoken\fi \@@next}%
\def\@stepthree{\afterassignment\@stepone\let\@@next= }%
%
%
%
\def\@getoptionalarg#1{%
   \let\@optionaltemp = #1%
   \let\@optionalnext = \relax
   \@futurenonspacelet\@optionalnext\@bracketcheck
}%
%
%
\def\@bracketcheck{%
   \ifx [\@optionalnext
      \expandafter\@@getoptionalarg
   \else
      \let\@optionalarg = \empty
      \expandafter\@optionaltemp
   \fi
}%
\def\@@getoptionalarg[#1]{%
   \def\@optionalarg{#1}%
   \@optionaltemp
}%
%
%
%
\def\@nnil{\@nil}%
\def\@fornoop#1\@@#2#3{}%
\def\@for#1:=#2\do#3{%
   \edef\@fortmp{#2}%
   \ifx\@fortmp\empty \else
      \expandafter\@forloop#2,\@nil,\@nil\@@#1{#3}%
   \fi
}%
\def\@forloop#1,#2,#3\@@#4#5{\def#4{#1}\ifx #4\@nnil \else
       #5\def#4{#2}\ifx #4\@nnil \else#5\@iforloop #3\@@#4{#5}\fi\fi
}%
\def\@iforloop#1,#2\@@#3#4{\def#3{#1}\ifx #3\@nnil
       \let\@nextwhile=\@fornoop \else
      #4\relax\let\@nextwhile=\@iforloop\fi\@nextwhile#2\@@#3{#4}%
}%
%
%
%
\@innernewif\if@fileexists
\def\@testfileexistence{\@getoptionalarg\@finishtestfileexistence}%
\def\@finishtestfileexistence#1{%
   \begingroup
      \def\extension{#1}%
      \immediate\openin0 =
         \ifx\@optionalarg\empty\jobname\else\@optionalarg\fi
         \ifx\extension\empty \else .#1\fi
         \space
      \ifeof 0
         \global\@fileexistsfalse
      \else
         \global\@fileexiststrue
      \fi
      \immediate\closein0
   \endgroup
}%
%
%
%
%
\def\bibliographystyle#1{%
   \@readauxfile
   \@writeaux{\string\bibstyle{#1}}%
}%
\let\bibstyle = \@gobble
%
%
\let\bblfilebasename = \jobname
\def\bibliography#1{%
   \@readauxfile
   \@writeaux{\string\bibdata{#1}}%
   \@testfileexistence[\bblfilebasename]{bbl}%
   \if@fileexists
      \nobreak
      \@readbblfile
   \fi
}%
\let\bibdata = \@gobble
%
%
\def\nocite#1{%
   \@readauxfile
   \@writeaux{\string\citation{#1}}%
}%
\@innernewif\if@notfirstcitation
%
%
\def\cite{\@getoptionalarg\@cite}%
%
%
\def\@cite#1{%
   \let\@citenotetext = \@optionalarg
   \printcitestart
   \nocite{#1}%
   \@notfirstcitationfalse
   \@for \@citation :=#1\do
   {%
      \expandafter\@onecitation\@citation\@@
   }%
   \ifx\empty\@citenotetext\else
      \printcitenote{\@citenotetext}%
   \fi
   \printcitefinish
}%
\def\@onecitation#1\@@{%
   \if@notfirstcitation
      \printbetweencitations
   \fi
   \expandafter \ifx \csname\@citelabel{#1}\endcsname \relax
      \if@citewarning
         \message{\@linenumber Undefined citation `#1'.}%
      \fi
      \expandafter\gdef\csname\@citelabel{#1}\endcsname{%
\strut
\vadjust{\vskip-\dp\strutbox
\vbox to 0pt{\vss\parindent0cm \leftskip=\hsize 
\advance\leftskip3mm
\advance\hsize 4cm\strut\openup-4pt 
\rightskip 0cm plus 1cm minus 0.5cm ?  #1 ?\strut}}
         {\tt
            \escapechar = -1
            \nobreak\hskip0pt
            \expandafter\string\csname#1\endcsname
            \nobreak\hskip0pt
         }%
      }%
   \fi
   \csname\@citelabel{#1}\endcsname
   \@notfirstcitationtrue
}%
%
%
\def\@citelabel#1{b@#1}%
%
%
\def\@citedef#1#2{\expandafter\gdef\csname\@citelabel{#1}\endcsname{#2}}%
%
%
%
\def\@readbblfile{%
   \ifx\@itemnum\@undefined
      \@innernewcount\@itemnum
   \fi
   \begingroup
      \def\begin##1##2{%
         \setbox0 = \hbox{\biblabelcontents{##2}}%
         \biblabelwidth = \wd0
      }%
      \def\end##1{}
      %
      %
      \@itemnum = 0
      \def\bibitem{\@getoptionalarg\@bibitem}%
      \def\@bibitem{%
         \ifx\@optionalarg\empty
            \expandafter\@numberedbibitem
         \else
            \expandafter\@alphabibitem
         \fi
      }%
      \def\@alphabibitem##1{%
         \expandafter \xdef\csname\@citelabel{##1}\endcsname {\@optionalarg}%
         \ifx\biblabelprecontents\@undefined
            \let\biblabelprecontents = \relax
         \fi
         \ifx\biblabelpostcontents\@undefined
            \let\biblabelpostcontents = \hss
         \fi
         \@finishbibitem{##1}%
      }%
      \def\@numberedbibitem##1{%
         \advance\@itemnum by 1
         \expandafter \xdef\csname\@citelabel{##1}\endcsname{\number\@itemnum}%
         \ifx\biblabelprecontents\@undefined
            \let\biblabelprecontents = \hss
         \fi
         \ifx\biblabelpostcontents\@undefined
            \let\biblabelpostcontents = \relax
         \fi
         \@finishbibitem{##1}%
      }%
      \def\@finishbibitem##1{%
         \biblabelprint{\csname\@citelabel{##1}\endcsname}%
         \@writeaux{\string\@citedef{##1}{\csname\@citelabel{##1}\endcsname}}%
         \ignorespaces
      }%
      %
      %
      \let\em = \bblem
      \let\newblock = \bblnewblock
      \let\sc = \bblsc
      \frenchspacing
      \clubpenalty = 4000 \widowpenalty = 4000
      \tolerance = 10000 \hfuzz = .5pt
      \everypar = {\hangindent = \biblabelwidth
                      \advance\hangindent by \biblabelextraspace}%
      \bblrm
      \parskip = 1.5ex plus .5ex minus .5ex
      \biblabelextraspace = .5em
      \bblhook
      \input \bblfilebasename.bbl
   \endgroup
}%
%
%
\@innernewdimen\biblabelwidth
\@innernewdimen\biblabelextraspace
%
%
%
\def\biblabelprint#1{%
   \noindent
   \hbox to \biblabelwidth{%
      \biblabelprecontents
      \biblabelcontents{#1}%
      \biblabelpostcontents
   }%
   \kern\biblabelextraspace
}%
%
%
%
\def\biblabelcontents#1{{\bblrm [#1]}}%
%
%
\def\bblrm{\rm}%
%
%
\def\bblem{\it}%
%
%
\def\bblsc{\ifx\@scfont\@undefined
              \font\@scfont = cmcsc10
           \fi
           \@scfont
}%
%
%
\def\bblnewblock{\hskip .11em plus .33em minus .07em }%
%
%
\let\bblhook = \empty
%
%
%
\def\printcitestart{[}
\def\printcitefinish{]}
\def\printbetweencitations{, }
\def\printcitenote#1{, #1}
%
%
%
\let\citation = \@gobble
%
%
%
\@innernewcount\@numparams
%
%
\def\newcommand#1{%
   \def\@commandname{#1}%
   \@getoptionalarg\@continuenewcommand
}%
%
%
\def\@continuenewcommand{%
   \@numparams = \ifx\@optionalarg\empty 0\else\@optionalarg \fi \relax
   \@newcommand
}%
%
%
\def\@newcommand#1{%
   \def\@startdef{\expandafter\edef\@commandname}%
   \ifnum\@numparams=0
      \let\@paramdef = \empty
   \else
      \ifnum\@numparams>9
         \errmessage{\the\@numparams\space is too many parameters}%
      \else
         \ifnum\@numparams<0
            \errmessage{\the\@numparams\space is too few parameters}%
         \else
            \edef\@paramdef{%
               \ifcase\@numparams
                  \empty  No arguments.
               \or ####1%
               \or ####1####2%
               \or ####1####2####3%
               \or ####1####2####3####4%
               \or ####1####2####3####4####5%
               \or ####1####2####3####4####5####6%
               \or ####1####2####3####4####5####6####7%
               \or ####1####2####3####4####5####6####7####8%
               \or ####1####2####3####4####5####6####7####8####9%
               \fi
            }%
         \fi
      \fi
   \fi
   \expandafter\@startdef\@paramdef{#1}%
}%
%
%
%
%
\def\@readauxfile{%
   \if@auxfiledone \else 
      \global\@auxfiledonetrue
      \@testfileexistence{aux}%
      \if@fileexists
         \begingroup
            \endlinechar = -1
            \catcode`@ = 11
            \input \jobname.aux
         \endgroup
      \else
         \message{\@undefinedmessage}%
         \global\@citewarningfalse
      \fi
      \immediate\openout\@auxfile = \jobname.aux
   \fi
}%
%
%
\newif\if@auxfiledone
\ifx\noauxfile\@undefined \else \@auxfiledonetrue\fi
%
%
%
%
\@innernewwrite\@auxfile
\def\@writeaux#1{\ifx\noauxfile\@undefined \write\@auxfile{#1}\fi}%
%
%
%
\ifx\@undefinedmessage\@undefined
   \def\@undefinedmessage{No .aux file; I won't give you warnings about
                          undefined citations.}%
\fi
%
%
\@innernewif\if@citewarning
\ifx\noauxfile\@undefined \@citewarningtrue\fi
%
%
%
\catcode`@ = \@oldatcatcode


\def\widestnumber#1#2{}

\def\rm{\fam0 \tenrm}

\def\fakesubhead#1\endsubhead{\bigskip\noindent{\bf#1}\par}



%
%
%

%

\font\textrsfs=rsfs10
\font\scriptrsfs=rsfs7
\font\scriptscriptrsfs=rsfs5

\newfam\rsfsfam
\textfont\rsfsfam=\textrsfs
\scriptfont\rsfsfam=\scriptrsfs
\scriptscriptfont\rsfsfam=\scriptscriptrsfs

\edef\oldcatcodeofat{\the\catcode`\@}
\catcode`\@11

\def\Cal@@#1{\noaccents@ \fam \rsfsfam #1}

\catcode`\@\oldcatcodeofat


\expandafter\ifx \csname margininit\endcsname \relax\else\margininit\fi

\pageheight{8.5truein}
\topmatter
\title{On the arrow property \\
 Sh782} \endtitle
\author {Saharon Shelah \thanks {\null\newline I would like to thank 
Alice Leonhardt for the beautiful typing. \null\newline
 Done 2/4/01 \null\newline
 Latest Revision - 01/Dec/15} \endthanks} \endauthor 


\endtopmatter
\document  
 
\newpage

\head {\S0 Introduction } \endhead  \resetall \sectno=0
\bigskip



Let $X$ be a finite set of alternatives.  A choice function $c$ is a
mapping which assigns to nonempty subsets $S$ of $X$ an element $c(S)$
of $S$.  A {\it rational} choice function ins one for which there is a
linear ordering on the alternatives such that $c(S)$ is the maximal
element of $S$ according to that ordering.  (We will consider choice
functions which are defined on subsets of $X$ of fixed cardinality
$k$.)

Arrow's impossibility theorem \cite{Arr50} asserts that under certain
natural conditions, if there are at least three alternatives then
every non-dictatorial social choice gives rise to a non-rational
choice function, i.e., there exist profiles such that the social
choice is not rational.  A profile is a finite list of linear orders
on the alternatives which represent the individual choices.  For
general references on Arrow's theorem and social choice functions see
\cite{Fis73}, \cite{Pel84} and \cite{Sen86}.

Non-rational classes of choice functions which may represent
individual behavior where considered in \cite{KRS01} and
\cite{Kal01}.  For example: $c(S)$ is the second largest element in
$S$ according to some ordering, or $c(S)$ is the median element of $S$
(assume $|S|$ is odd) according to some ordering.  Note that the
classes of choice functions in these classes are {\it symmetric} namely
are invariant under permutations of the alternatives.  Gil Kalai asked
if Arrow's theorem can be extended to the case when the individual
choices are not rational but rather belong to an arbitrary non-trivial
symmetric class of choice functions.  (A class is non-trivial if it
does not contain all choice functions.)  The main theorem of this
paper gives an affirmative answer in a very general setting.  See also
\cite{RF86} for general forms of Arrow's and related theorem.

The main result is related to clones which are studied in universal
algebras (but we do not use this theory).  On clones see \cite{Szb99}
and \cite{Szn96}.
\bn
\ub{Notation}:
\sn
1) $n,m,k,\ell,r,s,t,i,j$ natural numbers always $k$, \nl
many times $r$ are constant (there may be some misuses of $k$). \nl
2) $X$ a finite set. \nl
3) ${\frak C}$ a family of choice function on $\binom X k = \{Y:Y
\subseteq X,|Y| = k\}$. \nl
4) ${\Cal F}$ is a clone on $X$ (see Definition \scite{2.1}(2)). \nl
5) $a,b,e \in X$. \nl
6) $c,d, \in {\frak C}$. \nl
7) $f,g \in {\Cal F}$.
\newpage

\head {Annotated content} \endhead  \resetall 
\bn
\S1 Framework
\mr
\item "{${{}}$}"  [What are $X,{\frak C},{\Cal F} = Av({\frak C})$, the
Arrow property restriction to $\binom X k,{\frak C}$ is $(X,k) = FCF$
(debt: no connection for different $k-S$) and the Main theorem.  For
${\frak C},{\Cal F},r=r({\Cal F})$.]
\ermn
\ub{Part A}:  The simple case.
\mn
\S2
\mr
\item "{${{}}$}"  [Define a clone, $r({\Cal F})$.   If $f \in {\Cal
F}_{(r)}$ is not a monarchy, $r \ge 4$ on the family of not one-to-one
sequences $\bar a \in {}^r X$ \ub{then} $f$ is a projection, \scite{2.5}. \nl
Define $f_{r;\ell,k}$, basic implications on $f_{r;\ell,k} \in {\Cal
F}$, \scite{2.6}, \scite{2.7}. \nl
If $r=3,f \in {\Cal F}_{[s]}$ is not a monarchy on one-to-one triples,
then $f$ without loss of generality it is $f_{r;1,2}$ or $g_{r;1,2}$ on then
\scite{7.8}. \nl
If $r=3,f$ is not a semi monarchy on permutations of $\bar a$. \nl
If $r=3$, there are some ``useful" $f$, \scite{7.9B}.  Implications on
$f_{r;\ell,k} \in {\Cal F}$.] 
\ermn
\S3 Getting ${\frak C}$ is full
\mr
\item "{${{}}$}"  [Sufficient condition for $r \ge 4$ with $f_{r;1,2}$
or so (\scite{12.1}), end?  $r=3$? \nl
Sufficient condition for $r=3$ with $g_{r;1,2}$ or so
(\scite{12.2}). \nl
A pure sufficient condition for ${\frak C}$ full \scite{12.2a}. \nl
Subset $\binom X 3$, closed under a distance \scite{12.3}. \nl
Getting the final conclusion (relying on 3, \scite{12.4}).]
\ermn
\S4 The $r=2$ case.
\mr
\item "{${{}}$}"  [By stages we get a $f \in {\Cal F}_{[r]}$ which is a
monarchy with exactly one exceptional pair \scite{13.2} -
\scite{13.4}.  Then by composition we get $g \in {\Cal F}_2$ similar
to $f_{r;1,2}$.]
\ermn
\ub{Part B}:  Non simple case.
\mr
\item "{${{}}$}"  [Maybe we can start with this and prove from part A
only what is used; important question from presentation does not
interest me.]
\ermn
\S5 Fullness - the nonsimple case
\mr
\item "{${{}}$}"  [We derive ``${\frak C}$ is full" from various
assumptions, and then prove the main theorem.]
\ermn
\S6 The case $r=2$.
\mn
\S7 The case $r \ge 4$.
\newpage

\head {\S1 Framework} \endhead  \resetall \sectno=1
\bigskip

\demo{\stag{1.0} Context}  We fix a finite $X$ and $r = \{0,\dotsc,r-1\}$.
\enddemo
\bigskip

\definition{\stag{1.1} Definition}  1) An $(X,r)$-election rule is a
function $c$ such that: for every ``vote" $\bar t = \langle t_a:a \in
X \rangle \in {}^X r$ we have $c(\bar t) \in r = \{0,\dotsc,r-1\}$. \nl
2) $c$ is a monarchy if $(\exists a \in X)(\forall \bar t \in {}^X
\bar r)[c(\bar t) = t_a]$. \nl
3) $c$ is reasonable if $(\forall \bar t)(c(t) \in \{t_a:a \in X\})$.
\enddefinition
\bigskip

\definition{\stag{1.2} Definition}  1) We say ${\frak C}$ is a family of
choice functions for $X$ ($X-FCF$ in short) if:

$$
\align
{\frak C} \subseteq \{c:&c \text{ is a function} \\
  &\text{Dom}(c) = {\Cal P}^-(X) (= \text{ family of nonempty subsets
of } X) \\
  &\text{and } (\forall Y \in {\Cal P}^-(X))(c(Y) \in Y)\}.
\endalign
$$  
\mn
2) ${\frak C}$ is called symmetric if for every $\pi \in \text{ Per}(X)
=$ group of permutations of $X$ we have
$$
c \in {\frak C} \Rightarrow \pi * c \in {\frak C}
$$
\mn
where

$$
\pi * c(Y) = \pi^{-1}(c \pi(Y)).
$$
\mn
3)  ${\Cal P}_{\frak C} = {\Cal P}^-(X)$.
\enddefinition
\bigskip

\definition{\stag{1.4} Definition}  1) We say av is a $r$-averaging
function for ${\frak C}$ \ub{if}
\mr
\item "{$(a)$}"  av is a function written av$_Y(a_1,\dotsc,a_r)$
\sn
\item "{$(b)$}"  for any $c_1,\dotsc,c_r \in {\frak C}$, there is $c \in
{\frak C}$ such that \nl
$(\forall Y \in {\Cal P}^-(X))(c(Y)) = \text{ av}_Y(c_1(Y),\dotsc,c_r(Y))$
\sn
\item "{$(c)$}"  if $a \in Y \in {\Cal P}^-(X)$ \ub{then}
av$_Y(a,\dotsc,a) = a$.
\ermn
2) av is simple if av$_Y(a_1,\dotsc,a_r)$ does not depend on $Y$ so we
may omit $Y$. \nl
3) AV$_r({\frak C}) = \{\text{av}:\text{av}$ is 
an $r$-averaging funtion for ${\frak
C}\}$, similarly AV$^s_r({\frak C}) = \{\text{av}:\text{av is a}$
simple $r$-averaging function for ${\frak C}\}$. \nl
4) AV$({\frak C}) = \dbcu_r AV_r({\frak C})$ and $AV^s({\frak C}) = 
\dbcu_r AV^s_r({\frak C})$.
\enddefinition
\bigskip

\definition{\stag{1.3} Definition}  1) ${\frak C}$ has the simple $r$-Arrow
property \ub{if} \nl
av $\in AV^s_r({\frak C}) \Rightarrow \dsize \bigvee^r_{t=1}
(\forall a_1,\dotsc,a_2)$ av$(a_1,\dotsc,a_2) = a_t$ \nl
such av is called monarchial. \nl
2) Similarly without simple.
\enddefinition
\bn
\ub{\stag{1.3a} Question}:  1) When does ${\frak C}$ have the Arrow
property? \nl
2) Does $|{\frak C}| \le \text{ poly}(|X|) \Rightarrow r$-Arrow property?
\bigskip

\remark{\stag{1.5} Remark}  The question was asked with ${\frak C}_{(X)}$ defined for
every $X$; but in the treatment this does not influence.
\endremark
\bn
We actually deal with:
\definition{\stag{1.6} Definition}  If $1 \le k \le |X|-1$ and we
replace ${\Cal P}^-(X)$ by $\binom Xk =: \{Y:Y \subseteq X,|Y| = k\}$,
\ub{then} ${\frak C}$ is called $(X,k)$ - FCF, ${\Cal P}_{\frak C} =
\binom X k,k = k({\frak C})$, av is [simple]
$r$-averaging function for ${\frak C}$; let $k({\frak C}) = \infty$ if
${\Cal P}_{\frak C} = {\Cal P}^-(X)$; let ${\Cal F}({\frak C}) = 
Av^s({\frak C})$ let ${\Cal F}_{[r]} = \{f \in {\Cal F}:f$ is $r$-place$\}$.
\enddefinition
\bn
\ub{\stag{1.7} Discussion}:  This is justified: \nl
1) For simple averaging function, $k \ge r$ the restriction 
to $\binom Xk$ implies the full result.  [In the end, draw the conclusion.] \nl
2) For the nonsimple case there is a little connection between the
various ${\frak C} \restriction \binom Xk$ (exercise).  Our aim is (but
we shall first prove the simple case):
\bigskip

\proclaim{\stag{1.8} Main Theorem}  There are $r^*_1,r^*_2 < \omega$.
(We shall be able to give explicit values, e.g. $r^*_1 = r^*_2 = 7$ are O.K.) such
that:
\mr
\item "{$\circledast$}"  if $X$ is finite, $r^*_1 \le k,|X| - r^*_2 \ge
k$ and ${\frak C}$ is an $(X,k)$-FDF and some {\text{\rm av\/}} 
$\in \text{ AV}_r(\Bbb C)$ is not
monarchial, \ub{then} every choice function for $\binom Xk$ belongs to
${\frak C}$ ($= {\frak C}$ is full).
\endroster
\endproclaim
\newpage

\noindent
\ub{Part A}:  Simple case.
\bn
\head {\S2 Context and on nice $f$'s} \endhead  \resetall \sectno=2
\bn
\ub{Note}:  Sometimes Part B gives alternative ways.
\bigskip

\demo{\stag{2.0} Hypothesis (for part A)}
\mr
\item "{$(a)$}"  $X$ a finite set
\sn
\item "{$(b)$}"  $5 < k < |X| - 5$
\sn
\item "{$(c)$}"  ${\frak C}$ a symmetric $(X,k)$-FCF and ${\frak C} \ne
\emptyset$
\sn
\item "{$(d)$}"  ${\Cal F}_{[r]} = \{f:f \text{ an } r$-place function
from $X$ to $X$ such that ${\frak C}$ is closed under $f$ that is $f
\in \text{ AV}^s_r({\frak C})\}$
\sn
\item "{$(e)$}"  ${\Cal F} = \cup\{{\Cal F}_{[r]}:r < \omega\}$
\sn
\item "{$(f)$}"  $r({\Cal F}) = \text{ Min}\{r:\text{ some } f \in
{\Cal F}_{[r]}$ is not a monarchy$\}$.
\endroster
\enddemo
\bigskip

\demo{\stag{2.0a} Fact}  ${\Cal F}$ is a clone on $X$ (see \scite{2.1} below).
\enddemo
\bigskip

\definition{\stag{2.1} Definition}  1) $f$ is monarchial = projection
if $f$ is an $r$-place function (from $X$ to $X$) and for some
$t,(\forall x_1,\dotsc,x_n)f(x_1,\dotsc,x_r) = x_t$. \nl
2) ${\Cal F}$ is a clone on $X$ if it is a family of functions from $X$ to $X$ (for all
arities) including the projections and closed under composition.
\enddefinition
\bigskip

\definition{\stag{2.3} Definition}  For ${\frak C},{\Cal F}$ as in \scite{2.0}:

$$
r({\frak C}) = r_m({\Cal F}) =: \text{ Min}\{r:\text{ some } f 
\in {\frak C}_r \text{ is not monarchial}\}
$$
\mn
(let $r({\Cal F}) = \infty$ if ${\frak C}$ is monarchial).
\enddefinition
\bigskip

\proclaim{\stag{2.5} Claim}  Assume
\mr
\item "{$(a)$}"  $f \in {\Cal F}_{[r]}$
\sn
\item "{$(b)$}"  $4 \le r = r_m({\Cal F}) = \text{ Min}\{r: \text{
some } f \in {\Cal F}$ is not a monarchy$\}$.
\ermn
\ub{Then} \nl
1) for some $\ell \in \{1,\dotsc,r\}$ we have $f(x_1,\dotsc,x_r) =
x_\ell \text{ if } x_1,\dotsc,x_r$ has some repetition. \nl
2) $r \le k$. 
\endproclaim
\bigskip

\demo{Proof}  1) Clearly 
there is a two-place function $h$ from $\{1,\dotsc,r\}$ to
$\{1,\dotsc,r\}$ such that: if $y_\ell = y_k \wedge \ell \ne k
\Rightarrow f(y_1,\dotsc,y_r) = y_{h(\ell,k)}$; we have some freedom
so let \wilog
\mr
\item "{$\boxtimes$}"  $\ell \ne k \Rightarrow h(\ell,k) \ne k$.
\ermn
Assume toward contradiction that 
\mr
\item "{$\circledast$}"   (1)'s conclusion fails,
i.e. $h \restriction h(\ell,k):1 \le \ell < k \le r\}$ is not
constant.
\endroster
\enddemo
\bn
\ub{Case 1}:  For some $\bar x \in {}^r X$ and $\ell_1 \ne k_1 \in
\{1,\dotsc,r\}$ we have 

$$
x_{\ell_1} = x_{k_1}
$$

$$
f(\bar x) \ne x_{\ell_1}
$$
\mn
equivalently:  $h\{\ell_1,k_1\} \notin \{\ell_1,k_1\}$.
\sn
Without loss of generality $\ell_1 = r-1,k_1 =r,f(\bar x) = x_1$ (as
for a permutation $\sigma$ of $\{1,\dotsc,r\}$ we can replace $f$ by
$f_\sigma,f_\sigma(x_1,\dotsc,x_r) =
f(x_{\sigma(1)},\dotsc,x_{\sigma(r)}))$.

We can choose $x \ne y$ in $X$ so $h(x,y,\dotsc,y) = x$ hence $\ell
\ne k \in \{2,\dotsc,r\} \Rightarrow h(\ell,k) = 1$. \nl
Now for $\ell \in \{2,\dotsc,r\}$ we have agreed $h(1,\ell) \ne \ell$, (see
$\boxtimes$) so as $h \restriction \{(\ell,k):\ell < k\}$
is not constantly 1 (by $\circledast$)
\wilog \, $h(1,2) = 3$.  But as $r \ge 4$ letting $x \ne y \in X$ we
have $f(x,x,y,y)$ is $y$ as $h(1,2) = 3$ and is $x$ as $h(3,4) = 1$,
contradiction. 
\bn
\ub{Case 2}:  Not Case 1.

Let $x \ne y$ compute $f(x,x,y,y,\ldots)$ it is $x$ as $h(1,2) \in
\{1,2\}$ and it is $y$ as $h(3,4) \in \{3,4\}$, contradiction. \nl
2) Follows as for $r > k$ we always have a repetition.
\hfill$\square_{\scite{2.5}}$\margincite{2.5}
\bigskip

\definition{\stag{2.6} Definition}  $f_{r;\ell,k} = f_{r,\ell,k}$ is the $r$-place
function on $X$ defined by

$$
f_{r;\ell,k}(\bar x) = \cases
x_\ell \qquad &\bar x \text{ is with repetition} \\
x_k \qquad &\text{ otherwise}
\endcases
$$
\enddefinition 
\bigskip

\proclaim{\stag{2.7} Claim}  1) If $f_{r,1,2} \in {\Cal F}$ then
$f_{r,\ell,k} \in {\frak C}$ for $\ell \ne k \in \{1,\dotsc,r\}$. \nl
2) If $f_{r,1,2} \in {\Cal F}$ and $r = r[{\Cal F}] \ge 3$ 
then $f_{r+1,1,2} \in {\Cal F}$.
\endproclaim
\bigskip

\demo{Proof}  1) Trivial. \nl
2) \ub{First assume $r \ge 5$}. 
Let $g(x_1,\dotsc,x_{r+1}) = 
f_{r,1,2}(x_1,x_2,\tau_3 \ldots \tau_r)$ where \nl
$\tau_m \equiv f_{r,1,m}(x_1,\dotsc,x_m,x_{m+2},\dotsc,x_{r+1})$; (that is
$x_{m+1}$ is omitted).
\sn
So for any $\bar a$:

\ub{if} $\bar a$ with no repetitions:

$$
\tau_3(\bar a) = a_3,\dotsc,\tau_m(\bar a) = a_m
$$

$$
g(\bar a) = f(a_1,a_2,a_3,\dotsc,a_r) = a_2
$$  
\mn
\ub{if} $\bar a$ has repetitions say $a_\ell = a_k$ then there is $m \in
\{3,\dotsc,r\} \backslash \{\ell -1,k-1\}$ hence $\langle
a_1,\dotsc,a_m,a_{m+2},\dotsc,a_{r+1} \rangle$ is with repetition so
$\tau_m(\bar a) = a_1$ so
$(a_1,a_2,\dotsc,\tau_m(\bar a),\ldots)$ has repetition so $g(\bar a)
= a_1$.
\bn
\ub{Second assume $r=4$}:

Let $(\tau_1,\dotsc,\tau_4),g$ be of arity 5 and for $\bar x =
(x_1,\dotsc,x_5)$ we let $g(\bar x) = f_{r,1,2}(\tau_1(\bar
x),\dotsc,\tau_4(\bar x))$ where
\mr
\sn
\item "{$(*)_1$}"  $\tau_1(\bar x) = x_1$
\sn
\item "{$(*)_2$}"  $\tau_2(\bar x) = f_{r,1,2}(x_1,x_2,x_3,x_4)$
\sn
\item "{$(*)_3$}"  $\tau_3(\bar x) = 
f_{r,1,3}(x_1,x_2,x_3,x_5),\tau_4(\bar x) =
f_{r,1,3}(x_1,x_2,x_4,x_5)$. \nl
Note that
\sn
\item "{$(*)_4$}"   for $\bar x$ with no repetition
$\tau_\ell(\bar x) = x_\ell$.
\ermn
Now check that $g$, i.e. is as required.
\bn
\ub{Third assume $r=3$}:

Let $g'(x_1,x_2,x_3,x_4) = f_{r,1,2}(\tau_1,\tau_2,\tau_3)$ where

$$
\tau_1 = x_1
$$

$$
\tau_2 = f_{r,1,2}(x_1,x_2,x_4)
$$

$$
\tau_3 = f_{r,1,2}(x_2,x_3,x_4).
$$
\mn
Now check (or see \scite{13.7}'s proof).  \hfill$\square_{\scite{2.7}}$\margincite{2.7}
\enddemo
\bigskip

\proclaim{\stag{7.8} Claim}  Assume
\mr
\item "{$(\alpha)$}"  ${\Cal F}$ is a symmetric clone,
\sn
\item "{$(\beta)$}"  every $f \in {\Cal F}_{[2]}$ is a monarchy,
$r=r[{\Cal F}] = 3$
\sn
\item "{$(\gamma)$}"  $f^* \in {\Cal F}_{[3]}$ and for no $i \in
\{1,2,3\}$ do we have $(\forall \bar b \in {}^3 X)(\bar b$ not
one-to-one $\Rightarrow f^*(\bar b) = b_i)$. 
\ermn
\ub{Then} for some $g \in {\Cal F}_{[3]}$ not a monarchy we have:
(a) and (b) where
\mr
\item "{$(a)$}"  for $\bar b \in {}^3 X$ which is not one-to-one
$g(\bar b) = f_{r;1,2}(\bar b)$, i.e. $= b_1$
\sn 
\item "{$(b)$}"  for $\bar b \in {}^3 X$ which is not one-to-one
$g(\bar b) = g_{r;1,2}(\bar b)$, see below
\endroster
\endproclaim
\bn
Where
\definition{\stag{2.9} Definition}  $g_{r;1,2}$ is the following
function \footnote{this is the majority function for $r=3$} from $X$ to $X$. \nl
$g_{r;1,2,}(x_1,x_2,\dotsc,x_2) =
\cases x_2 \quad &\text{ if } x_2 = x_3 = \ldots = x_r \\
x_1 &\text{ if otherwise}\endcases$.
\enddefinition
\bigskip

\demo{Proof of \scite{7.8}}  The same as the proof of the next claim ignoring the
one-to-one sequences (i.e. $f(a_1,a_2,a_3)$), see more later.
\enddemo
\bigskip

\proclaim{\stag{2.9A} Claim}  Assume ${\Cal F}$ is a symmetric clone,
$r = r({\Cal F}) = 3,f^* \in {\Cal F},f^*$ is a 3-place function and
not a monarchy and $\bar a \in {}^3 X$ is with no
repetition such that: if $\bar a' = (a'_1,a'_2,a'_3)$ is a permutation of $\bar
a$ then $f^*(\bar a') = a'_1$; but $\neg(\forall \bar b \in {}^3 X)(\bar
b \text{ not one-to-one } \rightarrow f(\bar b) = b_1))$. \nl
\ub{Then} for some $g \in {\Cal F}_3$ we have (a) or we have (b) where:
\mr
\item "{$(a)(i)$}"   for $\bar b \in {}^3 X$ with repetition, 
$g(\bar b) = f_{r;1,2}(\bar b)$, i.e. $g(\bar b) = b_1$
\sn
\item "{${{}}(ii)$}"   $g(\bar a') = a'_2$ for any permutation $\bar
a'$ of $\bar a,\bar b$ 
\sn
\item "{$(b)(i)$}"   for $\bar b \in {}^3 X$ with repetition, $g(\bar
b) = g_{r;1,2}(\bar b)$
\sn
\item "{${{}}(ii)$}"   $g(\bar a') = a'_1$ for any permutation $\bar
a'$ of $\bar a$ (see on $g_{r;1,2}$ in \scite{2.9}).
\endroster 
\endproclaim
\bigskip

\demo{Proof}  Let $\bar a = (a_1,a_2,a_3);(a,b,c)$ denote any
permutation of $\bar a$.

Let $W = \{\bar b:\bar b \in {}^3 X$ and $[\bar b$ is a
permutation of $\bar a$ or $\bar b$ not one-to-one$]\}$. \nl
Let ${\Cal F}^- = \{f \restriction W:f \in {\Cal F}\},f = f^*
\restriction W$; recall that $f(a_1,a_2,a_3) = a_1$. \nl
Let for $\eta \in {}^3\{1,2\},f_\eta$ be the 3-place function with
domain $W$, such that \nl
[remember that $f(x,y,y),f(x,y,x),f(x,x,y)$ are monarchies and]:
\mr
\item "{$\boxtimes_0$}"  $f_\eta(a_{\sigma(1)},a_{\sigma(2)},a_{\sigma(3)}) =
a_{\sigma(1)}\text{ for } \sigma \in \text{ Per}\{1,2,3\}$
\sn
\item "{$\boxtimes_1$}"   $f_\eta(a_1,a_2,a_2) = a_{\eta(1)}$
\sn
\item "{$\boxtimes_2$}"  $f_\eta(a_1,a_2,a_1) = a_{\eta(2)}$
\sn
\item "{$\boxtimes_3$}"  $f_\eta(a_1,a_1,a_2) = a_{\eta(3)}$
\ermn
Now
\mr
\item "{$(*)_0$}"  $f \in \{f_\eta:\eta \in {}^3 2\}$ \nl
[why?  just think: by the assumption on $f^*$ and as $r({\Cal F}) = 3$?]
\sn
\item "{$(*)_1$}"  if $\eta = \langle 1,1,1\rangle$ then $f_\eta$ is
$\ne f$ \nl
[why?  $f_\eta(x_1,x_2,x_3) = x_1$ on $W$, i.e. is a monarchy]
\sn
\item "{$(*)_2$}"  if $\eta,\nu \in {}^3\{1,2\},\eta(1) =
\nu(1),\eta(2) = \nu(3),\eta(3) = \nu(2)$, then $f_\eta \in {\Cal
F}^- \Leftrightarrow f_\nu \in {\Cal F}^-$
\sn
[Why?  In $f(x,y,z)$ we just exchange $y$ and $z$]
\sn
\item  "{$(*)_3$}"  if $f_{<2,2,2,>} \in {\Cal F}^-$ then $f_{<1,2,2>}
\in {\Cal F}^-$
\nl
[Why?  Define $g$ by $g(x,y,z) =
f_{<2,2,2>}(x,f_{<2,2,2>}(y,x,z),f_{<2,2,2>}(z,x,y))$
\nl
(so $g \in {\Cal F}^-$) hence

$$
g(a,b,c) = f_{<2,2,2>}(a,b,c), \text{ generally } g \text{ satisfies } \boxtimes_0
$$

$$
\align
g(a,b,b) = &f_{<2,2,2>}(a,f_{<2,2,2>}(b,a,b),f_{<2,2,2>}(b,a,b))) \\
  &= f_{<2,2,2>}(a,a,a) = a
\endalign
$$

$$
\align
g(a,b,a) = &f_{<2,2,2>}(a,f_{<2,2,2>}(b,a,a),f_{<222>}(a,a,b)) \\
  &= f_{<2,2,2>}(a,a,b) = b
\endalign
$$

$$
\align
g(a,a,b) = &f_{<2,2,2>}(a,f_{<2,2,2>}(a,a,b),f_{<2,2,2>}(b,a,a)) \\
  &= f_{<2,2,2>}(a,b,a) = b.
\endalign
$$
\nl
So $g = f_{<1,2,2>}$
\item "{$(*)_4$}"  $f_{<1,2,2>} \in {\Cal F}^- \Rightarrow f_{<2,1,2>}
\in {\Cal F}^-$ \nl
[Why?  Let

$$
g(x,y,z) = f_{<1,2,2>}(x,y,f_{<1,2,2>}(z,x,y))
$$
So

$$
g(a,b,c) = a \text{ and } g \text{ satisfies } \boxtimes_0
$$
and

$$
\align
g(a,b,b) = &f_{<1,2,2>}(a,b,f_{<1,2,2>}(b,a,b)) \\
  &= f_{<1,2,2>}(a,b,a) = b
\endalign
$$

$$
\align
g(a,b,a) = &f_{<1,2,2>}(a,b,f_{<1,2,2>}(a,a,b)) \\
  &= f_{<1,2,2>}(a,b,b) = a
\endalign
$$

$$
\align
g(a,a,b) = &f_{<1,2,2>}(a,a,f_{<1,2,2>}(b,a,a)) \\
  &= f_{<1,2,2>}(a,a,b) = b.
\endalign
$$
So $g = f_{<2,1,2>}$ so $f_{<2,1,2>} \in {\Cal F}^-$ as promised
\sn
\item "{$(*)_5$}"  $f_{<2,1,2>} = f_{3;3,1}$, i.e.

$$
f_{<2,1,2>}(x_1,x_2,x_3) = \cases x_1 \quad &\text{ \ub{if} } \quad |\{x_1,x_2,x_3\}|
= 3 (\text{i.e. } (\{x_1,x_2,x_3\} = \bar a) \\
  x_3 &\text{ \ub{ if } } \quad |\{x_1,x_2,x_3\}| \le 2
\endcases
$$
\sn
when $(x_1,x_2,x_3) \in W$ \nl
[Why?  Check.]
\sn
\item "{$(*)_6$}"  $f_{<2,2,1>}(x_1,x_2,x_3) = x_2$ if $2 \ge
|\{x_1,x_2,x_3\}|$  \nl
[Why?  Check.]
\sn
\item "{$(*)_7$}"   $f_{<2,1,2>} \in {\Cal F}^- \Leftrightarrow
f_{<2,2,1>} \in {\Cal F}^-$ \nl
[Why?  See $(*)_2$ in the beginning.]
\sn
\item "{$(*)_8$}"  $f_{<1,2,1>} \in {\Cal F}^- \Leftrightarrow
f_{<1,1,2>} \in {\Cal F}^-$ \nl
[Why?  By $(*)_2$ in the beginning.]
\sn
\item  "{$(*)_9$}"   $f_{<1,2,1>} \in {\Cal F}^- \Rightarrow f_{<2,2,1>}
\in {\Cal F}^-$.
\nl
Let $g(x,y,z) = f_{<1,2,1>}(x,f_{<1,2,1>}(y,z,x),f_{<1,2,1>}(z,x,y))$

$$
\align
g(a,b,c) = &f_{<1,2,1>}(a,f_{<1,2,1>}(b,c,a),f_{<1,2,1>}(c,a,b)) \\
  &= f_{<1,2,1>}(a,b,c) = a \\
  &\text{and generally } g \text{ satisfies } \boxtimes_0
\endalign
$$

$$
g(a,b,b) = f_{<1,2,1>}(a,f_{<1,2,1>}(b,b,a),f_{<1,2,1>}(b,a,b)) = f_{<1,2,1>}(a,b,a)) = b
$$

$$
g(a,b,a) = f_{<1,2,1>}(a,f_{<1,2,1>}(b,a,a),f_{<1,2,1>}(a,a,b)) = 
f_{<1,2,1>}(a,b,a) = b
$$

$$
g(a,a,b) = f_{<1,2,1>}(a,f_{<1,2,1>}(a,b,a),f_{<1,2,1>}(b,a,a)) = 
f_{<1,2,1>}(a,b,b) = a.
$$
\mn

So $g = f_{<2,2,1>}$ hence $f_{<2,2,1>} \in {\Cal F}^-$.
\endroster
\enddemo
\newpage

\head {Diagram} \endhead  \resetall 
\bn
Diagram (arrows mean belonging to ${\Cal F}^-$ follows)

$$
f_{<2,2,2>} \in {\Cal F}^-
$$
\medskip

\hskip60pt $\downarrow (*)_3$

$$
f_{<1,2,2>} \in {\Cal F}^- \qquad \qquad f_{<1,2,1>} \in {\Cal
F}^-_{(*)_8} \underset (*)_8 {}\to 
\Leftrightarrow {\Cal F}_{<1,1,2>} \in {\Cal F}^-
$$
\medskip

\hskip60pt $\downarrow (*)_4 \qquad \quad \qquad  (*)_9 \swarrow$

$$
f_{<2,1,2>} \in {\Cal F} \underset (*)_7 {}\to
\leftrightarrow f_{<2,2,1>} \in {\Cal F}^-
$$
\mn
among the $2^3$ function $f_\eta$ one is discarded being a monarchy, see
$(*)_1$, six appear in the diagram above and implies 
$f_{r;2,1} \in {\Cal F}^-$ by $(*)_5$ hence clause (a) holds, 
and one is $g_{r;1,2}$ because
\mr
\item "{$(*)$}"  $g_{r;1,2} = f_{<2,1,1>}$ on $W$.
\ermn
(Why check?), so clause (b) holds.  \hfill$\square_{\scite{2.9A}}$\margincite{2.9A}
\bn
\ub{Continuation of the proof of \sciteu{2.8}}:  As $r({\Cal F}) = 3$ 
for some $\eta \in {}^3 2,f^*$ agrees with 
$f_\eta$ for all not one-to-one triples $\bar b$.
If $\eta = \langle 1,1,1 \rangle$ we contradict an assumption and if
$\eta = \langle 2,1,1 \rangle$, possibility (b) of \scite{7.8} holds.
If $\eta = \langle 2,1,2 \rangle$ then $f^*(\bar b) = b_2$ for $\bar b
\in {}^3 X$ not one-to-one (see $(*)_5$) and as $f^*$ is not a monarchy we get
possibility (a), possibly permutating of the variables; similarly
$\eta = \langle 2,2,1 \rangle$.  In the remaining case (see the
diagram), there is $f \in {\Cal F}$ agreeing on $\{\bar b \in {}^3
X:\bar b$ is not one-to-one$\}$ with $f_\eta$ for $\eta = \langle
1,2,2 \rangle$ or $\eta = \langle 1,2,1 \rangle$, \wilog \, $f^* = f$.

If $\eta = \langle 1,2,2 \rangle$ define $g$ as in $(*)_4$,
i.e. $g(x,y,z) = f^*(x,y,f^*(z,x,y))$ so for a non one-to-one sequence
$\bar b \in {}^3 X$ we have $g(\bar b) = f_{<2,1,2>}(\bar b) = b_3$.
If for some one-to-one $\bar a \in {}^3 X$ we have $f^*(a_3,a_1,a_2) \ne
a_3$ then $g(a_1,a_2,a_3) = f^*(a_1,a_2,f^*(a_3,a_1,a_2)) \in
\{a_1,a_2\}$ so permuting the variables we get possibility (a).  So we
are left with the case $\bar a \in {}^3 X$ is one-to-one $\Rightarrow
f^*(\bar a) = a_1$. 
\nl
Let us define $g \in {\Cal F}_{[3]}$ by 
$g(x_1,x_2,x_3) = f^*(f(x_2,x_3,x_1),x_3,x_2)$.  Let $\bar b \in {}^3
X$; if $\bar b$ is with no repetitions then $g(\bar b) =
f^*(b_2,b_3,b_2) = b_3$.
If $\bar b = (a,b,b)$ then $g(\bar b) = f^*(f^*(b,b,a),b,b) =
f^*(a,b,b,) = a = b_1$ and if $\bar b = (a,b,a)$ then $g(\bar b) =
f^*(f^*(b,a,a),a,b) = f^*(b,a,b) = a = b_1$ and if 
$\bar b = (a,a,b)$ then $g(\bar b) = f^*(f^*(a,a,b),b,a) = 
f^*(b,b,a) = a = b_1$, so for $\bar b$ not one to one, $g(\bar b) =
b_1$.  So $g$ is as required in clause (a).

Lastly, let $\eta = \langle 1,2,1 \rangle$ and let $g(x,y,z) =
f^*(x,f^*(y,z,x),f^*(z,x,y))$ by $(*)_9$,  
easily [$\bar b$ is not one-to-one $\Rightarrow g(\bar b)
= f_{<2,2,1>}(\bar b) = b_2$].  Now $(a_1,a_2,a_3)$ is with no
repetitions and $f^*(a_2,a_3,a_1) = a_1$ then $g(a_1,a_2,a_3) = a_1$
and possibility $(a)$ holds.  Otherwise for $[\bar b \in{}^3 X$ is
one-to-one $\Rightarrow f^*(\bar b) \in \{b_1,b_2\}]$; so if
$(a_1,a_2,a_3) \in {}^3 X$ is one-to-one and $f^*(a_2,a_3,a_1) \ne
a_2$ then $g(a_1,a_2,a_3) \ne a_2$ (as $f^*(a_3,a_1,a_2) \ne a_2$) so
possibility $(a)$ holds.  Hence $[\bar b \in {}^3 X$ is one-to-one
$\Rightarrow f^*(\bar b) = b_2]$.
Let $g^*$ be $g^*(x,y,z) = f^*(f^*(x,y,z),f^*(x,z,y),x)$.  Now if
$\bar b$ is one to one then $g^*(\bar b) = f^*(b_2,b_3,b_1) = b_3$.
Also if $\bar b = (a,b,b)$ then $g^*(\bar b) =
f^*(f^*(a,b,b),f^*(a,b,b),a) = f^*(a,a,a) = a$, and if $\bar b =
(a,b,a)$ then $g^*(\bar b) = f^*(f^*(a,b,a),f^*(a,a,b),a) = f^*(b,a,a)
= b$ and if $\bar b = (a,a,b)$ then $g^*(\bar b) =
f^*(f^*(a,a,b),f^*(a,b,a),a) = f^*(a,b,a) = b$.  So $g^*$ is as
required in the case $\eta = \langle 1,2,2 \rangle$.
\hfill$\square_{\scite{7.8}}$\margincite{7.8}
\bigskip

\proclaim{\stag{7.9B} Claim}  Assume
\mr
\item "{$(\alpha)$}"  ${\Cal F}$ is a symmetric clone of $X$
\sn
\item "{$(\beta)$}"  every $f \in {\Cal F}_{[2]}$ is monarchial
\sn
\item "{$(\gamma)$}"  $f^* \in {\Cal F}_{[3]}$ is not monarchial.
\ermn
\ub{Then} one of the following holds
\mr
\item "{$(a)$}"  for every one-to-one $\bar a \in {}^3 X$ for some $f
= f_{\bar a}$ we have
{\roster
\itemitem{ $(i)$ }  $f_{\bar a}(\bar a) = a_2$ (even $f_{\bar a}(\bar
a') = a'_2$ for $\bar a'$ a permutation of $\bar a$)
\sn
\itemitem{ $(ii)$ }  if $\bar b \in {}^3 X$ is not one-to-one then
$f_{\bar a}(\bar b) = b_1$
\endroster}
\item "{$(b)$}"  for every one-to-one $\bar a \in {}^3 X$, for some $f
= f_{\bar a} \in {\Cal F}_{[3]}$ we have:
{\roster
\itemitem{ $(i)$ }  if $\bar b$ is a permutation of $\bar a$ then
$f_{\bar a}(\bar b) = b_1$
\sn
\itemitem{ $(ii)$ }  if $\bar b \in {}^3 X$ is not one-to-one then
$f_{\bar a}(\bar b) = g_{r;1,2}(\bar b)$
\endroster}
\sn
\item "{$(c)$}"  for every one-to-one $\bar a \in {}^3 X$, for some
$f = f_{\bar a} \in {\Cal F}_{[3]}$ we have
{\roster
\itemitem{ $(i)$ }  if $\bar b$ is a permutation of $\bar a$ then
$f_{\bar a}(\bar b) = a_1$
\sn
\itemitem{ $(ii)$ }  if $\bar b \in {}^3 X$ is not one-to-one then
$f_{\bar a}(\bar b) = g_{r;1,2}(\bar b)$
\endroster}
\endroster
\endproclaim
\bigskip

\demo{Proof}  As ${\Cal F}$ is symmetric, it suffices to prove ``for
some $\bar a$" instead of ``for every $\bar a$".
\mn
\ub{Case 1}:  For some $\ell(*)$ if $\bar b \in {}^3 X$ is not one to
one then $f^*(\bar b) = b_{\ell(*)}$.

As $f^*$ is not monarchial for some one-to-one $\bar a \in {}^3
X,f^*(\bar a) \ne a_{\ell(*)}$ say $f^*(\bar a) = a_{k(*)},k(*) \ne
\ell(*)$.  As ${\Cal F}$ is symmetrical \wilog \, $\ell(*) = 1,k(*) =
2$.  So possibility (a) holds.
\mn
\ub{Case 2}:  Not case 1.

By \scite{7.8}, \wilog $f^*$ satisfies (a),(b) or (c) of \scite{7.8} with
$f^*$ instead of $g$.  But clause (a) of \scite{7.8} is case 1 above.
So we can assume that case (b) of \scite{7.8} holds, i.e.
\mr
\item "{$(*)$}"  if $\bar b \in {}^3 X$ is not one-to-one then
$f^*(\bar b) = g_{r;1,2}$, i.e.

$$
f^*(\bar b) = \cases b_2 &\text{ if } b_2 = b_3 \\
b_1 &\text{ if } b_2 \ne b_3
\endcases
$$
\ermn
If \scite{2.9A} applies we are done as (a) or (b) of \scite{2.9A}
holds hence (a) or (b) of \scite{7.9B} respectively holds, so assume
\scite{2.9A} does not apply.  So consider a one-to-one 
sequence $\bar a \in {}^3 X$ and it follows that we have sequences
$\bar a^1,\bar a^2$, both a permutation of $\bar a$ such that
$\dsize \bigwedge_i [(f^*(\bar a^1) = a^1_i) \equiv (f^*(\bar a^2) \ne
a^2_i)]$.

Using closure under composition of ${\Cal F}$ and its being symmetric
for every permutation $\sigma$ of $\{1,2,3\}$
(and as $g_{r;1,2}$ is preserved by permuting the variables)
 there is $f_\sigma \in
{\Cal F}_{[3]}$ such that:
\mr
\item "{$(i)$}"  $f_\sigma(a_{\sigma(1)},a_{\sigma(2)},a_{\sigma(3)})
= a_1$
\sn
\item "{$(ii)$}"  if $\bar b \in {}^3 X$ not one-to-one then $f(\bar
b) = g_{r;1,2}(\bar b)$.
\ermn
Let $\langle \sigma_\rho:\rho \in {}^3 2 \rangle$ list the
permutations of $\{1,2,3\}$, necessarily with repetitions.  Now we
define by downward induction of $k \le 3,f_\rho \in {\Cal F}$ for
$\rho \in {}^k 2$ (sequences of zeroes and ones of length $k$) as follows:

$$
\ell g(\rho) = 3 \Rightarrow f_\rho = f_{\sigma_\rho}
$$

$$
\ell g(\rho) < 3 \Rightarrow f_\rho(x_1,x_2,x_3) = f_\rho(x_1,f_{\rho
\char 94 <0>} (x_1,x_2,x_3),f_{\rho \char 94 <1>} (x_1,x_2,x_3)),
$$
\mn
Easily (by downward induction):
\mr
\item "{$(*)_1$}"  if $\bar b \in {}^3 X$ is with repetitions then
$f_\rho(\bar b) = g_{r;1,2}(\bar b)$ (as $g_{r;1,2}$ act as majority). \nl
Now we prove by downward induction on $k \le 3$
\sn
\item "{$(*)_2$}"  if $\bar b$ is a permutation of $\bar a,\rho \in
{}^k 2,\rho \triangleleft \nu \in {}^3 2$ and $f_\nu(\bar b) = a_1$
then $f_\rho(\bar b) = a_1$.
\ermn
This is straight and so $f_{<>}$ is as required in clause (c).
\hfill$\square_{\scite{7.9B}}$\margincite{7.9B}
\enddemo
\bn
Similarly
\proclaim{\stag{7.10} Claim}  1) If $g_{r;\ell,k} \in {\Cal F}$ \ub{then}
\mr
\item "{$(a)$}"  $g_{r;\ell_1,k_1} 
\in {\Cal F}$ when $\ell_1 \ne k_1\in \{1,\dotsc,r\}$
\sn
\item "{$(b)$}"  $g_{s;1,2} \in {\Cal F}$ for $s \ge r$ [seemingly not
used].
\endroster
\endproclaim 
\bigskip

\demo{Proof}  1)(a)  Trivial. \nl
(b) By induction on $s$ so \wilog \, $s = r+1$.  Let
$g(x_1,\dotsc,x_{r+1}) = g_{r;1,2}(x_1,g_{r;1,2}
(x_1,\dotsc,x_r),x_3,\dotsc,x_{r+1})$ so $g \in {\frak C}$ and \nl
$g(a,b,b,\dotsc,b) =
g_{r;1,2}(a,g_{r;1,2}(a,b,\dotsc,b),b,\dotsc,b) = g_{r;1,2}
(a,b,b,\ldots) = b$. \nl
Now if $g_{r+1;1,2}(a_1,\dotsc,a_{r+1}) = a_1$ then
\sn
\ub{Case 1}:  $g_{r;1,2}(a_1,a_2,\dotsc,a_r) = a_1$ - trivial.
\sn
\ub{Case 2}:  Not 1; so $(a_1,\dotsc,a_{r+1}) = (a,b,\dotsc,b,c)$
where $b \ne c$ so easy to check.
\enddemo
\newpage

\head {\S3 Getting ${\frak C}$ is full} \endhead  \resetall \sectno=3
\bigskip

\demo{\stag{12.1} Lemma}:  Assume
\mr
\item "{$(a)$}"  $r \ge 3$, ${\Cal F}$ is a clone on $X$
{\roster
\itemitem{ $(*)$ }  $f_{r;1,2} \in {\Cal F}$ or just
\sn
\itemitem{ \,\,\,$(*)^-$ }  if $\bar a \in {}^r X$ one-to-one for some $f =
f_{\bar a} \in {\Cal F},f_{\bar a}(\bar a) = a_2$ and $[\bar b \in {}^r X$
not one-to-one $\Rightarrow f_{\bar a}(\bar b) = b_1]$
\endroster}
\item "{$(b)$}"  ${\frak C}$ is a (non empty) family of choice functions for
$\binom Xk = \{Y \subseteq X:|Y| = k\}$
\sn
\item "{$(c)$}"  ${\frak C}$ is closed under every $f \in {\Cal F}$
\sn
\item "{$(d)$}"  ${\frak C}$ is symmetric
\sn
\item "{$(e)$}"  $k \ge r > 2,k \ge 5,|X| - k \ge 7,r$.
\ermn
\ub{Then} ${\frak C}$ is full (i.e. every choice function is in).
\enddemo
\bigskip

\demo{Proof}  Without loss of generality $r \ge 4$ (if $r = 3$ then
clause (e) is fine also for $r=4$, 
the case $(*)$ holds is O.K. by \scite{2.7}, and if $(*)^-$ we
repeat the proof of \scite{2.7} for the case $r=3$, only
$g(x_1,x_2,x_3,x_4) = f_{<a_1,a_2,a_3>}(x_1,\tau_2,\tau_3)$ where $\tau_2 =
f_{<a_1,a_2,a_4>}(x_1,x_2,x_4),\tau_3 =
f_{<a_1,a_3,a_4>}(x_1,x_3,x_4)$ where for one-to-one $\bar a \in
{}^3 X,f_{\bar a}$ is defined by the symmetry; this is the proof of \scite{13.7}). 
Assume
\mr
\item "{$\boxtimes$}"  $c^*_1 \in {\frak C},Y^* \in \binom Xk,c^*_1(Y^*)
= a^*_1$ and $a^*_2 \in Y^* \backslash \{a^*_1\}$.
\ermn
\ub{Question}:  Is there $c \in {\frak C}$ such that $c(Y^*) =
a^*_2$ and $(\forall Y \in \binom Xk)(Y \ne Y^* \Rightarrow c(Y) =
c^*_1(Y))$? \nl
Choose $c^*_2 \in {\frak C}$ such that
\mr
\item "{$(a)$}"  $c^*_2(Y^*) = a^*_2$
\sn
\item "{$(b)$}"  $n(c^*_2) = |\{Y \in \binom Xk:c^*_2(Y) =
c^*_1(Y)\}|$ is maximal under (a).
\ermn
Easily ${\frak C}$ is not a singleton so $n(c^*_2)$ is well defined.
\enddemo
\bn
\ub{\stag{12.1A} 
Subfact}:  A positive answer to the question implies that ${\frak C}$
is full. \nl
[Why?  Easy.]  

Hence if $n(c^*_2) = \binom {|X|}k -1$ we are done so assume not
and let $Z \in \binom Xk,Z \ne Y^*,c^*_1(Z) \ne c^*_2(Z)$.
\bn
\ub{Case 1}:  For some $Z$ as above and $c^*_3 \in {\frak C}$ we have

$$
c^*_3(Y^*) \notin \{a^*_1,a^*_2\}
$$

$$
c^*_3(Z) \in \{c^*_1(Z),c^*_2(Z)).
$$
\mn
If so, let $a^*_3 = c^*_3(Y^*)$ and $a^*_4 \in Y^* \backslash
\{a^*_1,a^*_2,a^*_3\}$, etc. so $\langle a^*_1,\dotsc,a^*_r \rangle$
is one-to-one, $a^*_\ell \in Y^*$. \nl
Let $c^*_\ell \in {\frak C}$ for $\ell = 4,\ldots$ 
be such that $c^*_\ell(Y^*) = a_\ell$ exists as ${\frak
C}$ is symmetric. \nl
By assumption (a) we can choose $f \in {\Cal F}_{[r]}$ such that

$$
{{}} \cases f(a^*_1,\dotsc,a^*_r) = a^*_2, \\
\bar a \in {}^r X \text{ is with repetitions } \rightarrow f(\bar a) =
a_1
\endcases
$$
\mn
Let $c = f(c^*_1,c^*_2,\dotsc,c^*_r)$ so $c \in {\frak C}$ and:

$$
c(Y^*) = f(a^*_1,a^*_2,\dotsc,a^*_r) = a^*_2
$$

$$
\align
Y \in \binom Xk \and c^*_1(Y) = c^*_2(Y) \Rightarrow c(Y) &=
f(c^*_1(Y),c^*_2(Y),\ldots) \\
  &= f(c^*_1(Y),c^*_1(Y),\ldots) = c^*_1(Y)
\endalign
$$

$$
\align
c(Z) = &f(c^*_1(Z),c^*_2(Z),c^*_3(Z),\ldots) \\
  &= c^*_1(Z) \text{ (as } |\{c^*_1(Z),c^*_2(Z),c^*_3(Z)\}| \le 2).
\endalign
$$
\mn
So $c$ contradicts the choice of $c^*_2$.
\bn
\ub{Case 2}:  There are $c^*_3,c^*_4 \in {\frak C}$ such that $c^*_3(Y^*) \ne
c^*_4(Y^*)$ are $\ne a^*_1,a^*_2$ but $c^*_3(Z) = c^*_4(Z)$ or at
least $|\{c^*_1(Z),c^*_2(Z),c^*_3(Z),c^*_4(Z)\}| < 4$.
\bigskip

\demo{Proof}  Similar.
\enddemo
\bn
\ub{Case 3}:  Not case 1 nor 2.

Let ${\Cal P} = \{Z:Z \subseteq X,|Z| = k$ and $c^*_1(Z) \ne c^*_2(Z)\}$ so 
\mr
\item "{$(*)_1$}"  ${\Cal P} \ne \emptyset,\{Y^*\}$. \nl
[Why?  ${\Cal P} \ne \{Y^*\}$ by the subfact above.  Also we can find
$Z \in \binom X k$ such that $|Y^* \cap Z| = 2,c^*_1(Y^*) \notin Z$.
Let $\pi \in \text{ Per}(X)$ be the identity $Z,\pi(c^*_1(Y^*)) \ne
c^*_1(Y^*)$.  So conjugating $c^*_1$ by $\pi$ we get $c^*_2,n(c^*_2) > 0$.]
\sn
\item "{$(*)_2$}"  if $Z \in {\Cal P},c \in {\frak C}$ and $c(Z) \in
\{c^*_1(Z),c^*_2(Z)\}$ then $c(Y^*) \in \{c^*_1(Y^*),c^*_2(Y^*)\}$. \nl
[Why?  By not case 1.]
\endroster
\bn
\ub{Subcase 3a}:  For some $Z$ we have

$$
{{}} \cases Z \in {\Cal P} &\text{ and}: \\
|Y^* \backslash Z| \ge 4 &\text{ or just } |Y^* \backslash Z
\backslash \{a^*_1,a^*_2\}| \ge 2 \text{ and } |Y^* \backslash Z| \ge 3 \endcases
$$
\mn
Let $b_1,b_2,b_3 \in Y^* \backslash Z$ be pairwise distinct.  As
${\frak C}$ is symmetric there are $d_1,d_2,d_3 \in {\frak C}$ such
that $d_\ell(Y^*) = b_\ell$ for $\ell = 1,2,3$.  The number of
possible truth values of $d_\ell(Z) \in Y^*$ is 2 so \wilog \, $d_1(Z) \in Y^*
\leftrightarrow d_2(Z) \in Y^*$ and we can forget $b_3,d_3$. \nl
So for some $\pi \in \text{ Per}(X)$ we have $\pi(Y^*) = Y^*,\pi(Z) =
Z,\pi \restriction (Y^* \backslash Z) =$ identity hence $\pi(b_\ell) =
b_\ell$ for $\ell = 1,2$ and $\pi(d_1(Z)) = d_2(Z)$, note that
$d_\ell(Z) \in Z$, so this is possible; so \wilog \, $d_1(Z) = d_2(Z)$.

As $|Y^* \backslash Z \backslash \{a^*_2,a^*_2\}| \ge 2$, using another
$\pi \in \text{ Per}(X)$ without loss of generality $\{b_1,b_2\} \cap
\{a^*_1,a^*_2\} = \emptyset$.  So $d_1,d_2$ gives a contradiction by
case 2.
\bigskip

\remark{Remark}  This is enough for non polynomial $|{\frak C}|$ as
$|\{Y:|Y \backslash Z^*| \le 3\}| \le |Y|^6$.
\endremark
\bn
\ub{Case 3b}:  Not case 3a.

So $Z \in {\Cal P} \backslash \{Y^*\} \Rightarrow |Z \backslash Y^*|
\le 3$ hence (recalling $|Z \backslash Y^*| = |Y^* \backslash Z|$) we
have $Z \in {\Cal P} \backslash \{Y^*\} \Rightarrow |Z \cap Y^*| 
\ge k-3 \ge 1$. \nl
Now
\mr
\item "{$\boxtimes_0$}"  for $Z \in {\Cal P} \backslash \{Y^*\}$
there is $c^* \in {\frak C}$ such that
$c^*(Y^*) \ne c^*(Z)$ \nl
[Why? Otherwise ``by ${\frak C}$ is symmetric" for any $Z \in {\Cal P}
\backslash \{Y^*\}$ we have:
\ermn

$$
c \in {\frak C} \wedge Y',Y'' \in \binom Xk \wedge |Y' \cap Y''| = |Z \cap
Y^*| \Rightarrow c(Y') = c(Y'') \tag"{$\circledast$}"
$$

Define a graph ${\frak G} = {\frak G}_Z$: the set of nodes $\binom Xk$

$$
\text{the set of edge }\{(Y',Y''):|Y' \cap Y''| = |Y^* \cap Z|\}
$$
\mn
now this graph is connected: if ${\Cal P}_1,{\Cal P}_2$ are nonempty
disjoint set of nodes with union $\binom Xk$, then there is a cross
edge by \scite{12.3} below (why? clause $(\alpha)$ there is impossible
by $(*)_1$ and clause $(\beta)$ is impossible by 
the first line of case 3).  This gives
contradiction to $\circledast$.  So $\boxtimes_0$ holds.] \nl
\mn
We claim:
\mr
\item "{$\boxtimes_1$}"  for $Z \in {\Cal P}$ and $d \in {\frak C}$ we
have \nl
$d(Y^*) \in Z \cap Y^* \Rightarrow d(Z) = d(Y^*)$. \nl
[Why?  Assume $d,Z$ forms a counterexample; as $|Y^* \backslash Z| \le
3$ and $k \ge 7$ (see \scite{12.1}(e)) 
so if $k \ge 8$ then $|Y^* \cap Z| \ge k-3 \ge 5$ so $Y^* \cap Z
\backslash \{a^*_1,a^*_2\}$ has $\ge 3$ members looking again at
subcase 3 this always holds.  Now
for some $\pi_1,\pi_2 \in \text{ Per}(X)$ 
we have $\pi_1(Y^*) = Y^* = \pi_2(Y^*),\pi_1(Z) = Z =
\pi_2(Z),\pi_1(d(Z)) = \pi_2(d(Z)),\pi_1(d(Y^*)) \ne \pi_2(d(Y^*))$ are from
$Z \cap Y^* \backslash \{a^*_1,a^*_2\}$; recall we are assuming that
$d(Y^*) \in Z \cap Y^*$ and $d(Z) \ne d(Y^*)$.  Let $d_1,d_2$ be gotten from
$d$ by conjugating by $\pi_1,\pi_2$, so we get Case 2, contradiction.]
\sn
\item "{$\boxtimes_2$}"  if $d \in {\frak C},Y \in \binom Xk$ and
$d(Y) = a$ then \nl
$(\forall Y')(a \in Y' \in \binom Xk \rightarrow d(Y') = a)$. \nl
[Why?  By $\boxtimes_1 + ``{\frak C}$ closed under permutations of
$X$", we get: if $k^* \in \{|Z \cap Y^*|:Z \in {\Cal P} \backslash
\{Y^*\}\}$ (which is not empty) then: if $Z_1,Z_2 \in \binom X k,|Z_1
\cap Z_2| = k^*,d \in {\frak C}$ and $d(Z_1) \in Z_2$ then
$d(Z_1) = d(Z_2)$.  Clearly $k^* < k$ (by $Z \ne Y^*$) and $Zk-k^* \le
|X|$.  First assume that we can choose such $k^* > 0$.
Now for any $d \in {\frak C}$ and $a \in X$, claim
\scite{12.3} below applied to $k^*-1,k-1,X \backslash \{a\},
(\{Y \backslash \{a\}:a \in Y$ and $d(Y) = a\},\{Y \backslash \{a\}:a
\in Y$ and $d(Y) \ne a\})$.
We can demand $a \in \text{ Rang}(d)$.  Now clause $(\alpha)$ there
gives the desired conclusion (for $Y,a$ as in $\boxtimes_2$).  As we
know $k-k^* \le 3,k \ge 7$ clause $(\beta)$ is imposible too so we are done.]
\ermn
Now we get a contradiction: as said above in $\boxtimes_0$ for some $c^*
\in {\frak C}$ and $Z \in {\Cal P} \backslash \{Y^*\}$ we have
$c^*(Y^*) \ne c^*(Z)$, choose $Y \in \binom Xk$ such
that $\{c^*(Y^*),c^*(Z)\} \subseteq Y$.  So by $\boxtimes_2$ we have
$d(Y) = d(Y^*)$ and also $d(Y) = d(Z)$, contradiction. \hfill$\square_{\scite{12.1}}$\margincite{12.1}  
\bigskip

\proclaim{\stag{12.2} Claim}  In \scite{12.1} we can replace $(a)$ by
\mr
\item "{$(a)^*(i)$}"  ${\Cal F}$ a clone on $X$,
\sn
\item "{${{}}(ii)$}"  $g^* \in {\Cal F}_{[3]}$ where (note $g^* =
g_{3;1,2}$)

$$
g^*(x_1,x_2,x_3) = \cases x_2  &x_2 = x_3 \\
x_1 \quad &\text{ otherwise}
\endcases
$$
or just
\sn
\item "{${{}}(ii)^-$}"  for $\bar a^* \in {}^r X$ with no repetitions
for some $g = g_{\bar a^*},g(\bar a^*) = a^*_1$ and if $\bar a \in
{}^r X$ is with repetitions then $g_{{\bar a}^*}(\bar a) = g^*(\bar
a)$.
\endroster 
\endproclaim
\bigskip

\demo{Proof}   Let $c^*_1 \in {\frak C},Y^* \in \binom Xk,a^*_1 =
c^*_1(Y^*),a^*_2 \in Y^* \backslash \{a^*_1\}$,  we choose
$c^*_2$ as in the proof of \scite{12.1}.

Let ${\Cal P} = \{Y:Y \in \binom Xk,Y \ne Y^*,c^*_1(Y) \ne
c^*_2(Y)\}$; we assume ${\Cal P} \ne \emptyset$ and shall get a
contradiction, (this suffices).
\mr
\item "{$(*)_1$}"  there are no $Z \in {\Cal P}$ and $d \in {\frak C}$
such that

$$
d(Y^*) = c^*_2(Y^*)
$$

$$
d(Z) \ne c^*_2(Z).
$$ 
[Why?  If so, let $c = g(c^*_1,c^*_2,d),g$ is $g^*$ or just any
$g_{<c^*_1(Z),c^*_2(Z),d(Z)>}$ (from $(a)^*(iii)$ of the assumption.
\ermn
So $c \in {\frak C}$ and
\mr
\item "{$(A)$}"  $c(Y^*) = g(c^*_1(Y^*),c^*_2(Y^*),d(Y^*)) =
g(c^*_1(Y),c^*_2(Y^*),c^*_2(Y^*)) = c^*_2(Y^*)$
\sn
\item "{$(B)$}"  $c(Z) = g(c^*_1(Z),c^*_2(Z),d(Z)) = c^*_1(Z)$ as
$d(Z) \ne c^*_2(Z)$ \nl
(two cases: if $\langle c^*_1(Z),c^*_2(Z),d(Z)\rangle$ with no
repetitions - by the choice of $g$, otherwise it is equal to
$g^*(c^*_1(Z),c^*_2(Z),d(Z))$
\sn
\item "{$(C)$}"  $Y \in \binom Yk,Y \ne Y^*,Z \notin {\Cal P}
\Rightarrow c^*_2(Y) = c^*_1(Y) \Rightarrow c(Y) =
g(c^*_1(Y),c^*_2(Y),d(Y)) = g^*(c^*_1(Y),c^*_1(Y),d(Y)) = c^*_1(Y)$.  
\hfill$\square_{(*)_1}$
\ermn
So $(*)_1$ holds by $c^*_2$-s choice.]
\mr
\item "{$(*)_2$}"  if $\pi \in \text{ Per}(X);\pi(Y^*) =
Y^*,\pi(c^*_2(Y^*)) = c^*_2(Y^*)$ then
{\roster
\itemitem{ $(\alpha)$ }  $Y \in {\Cal P} \and \pi(Y) = Y \Rightarrow
\pi(c^*_2(Y)) = c^*_2(Y)$ 
\sn
\itemitem{ $(\beta)$ } $Y \in {\Cal P} \Rightarrow c^*_2(\pi(Y)) =
\pi(c^*_2(Y))$. \nl
[Why?  Otherwise may ``conjugate" $c^*_2$ by $\pi^{-1}$ getting $d \in
{\frak C}$ which gives a contradiction to $(*)_1$.]
\endroster}
\item "{$(*)_3$}"  let $Z \in {\Cal P}$ then there are no $d_1,d_2
\in {\frak C}$ such that \nl
$d_1(Z) = d_2(Z) \ne c^*_2(Z)$ \nl
$d_1(Y^*) \ne d_2(Y^*)$. 
\nl
[Why?  Let $g = g_{<c^*_2(Y^*),d_1(Y^*),d_2(Y^*)>}$ be as in the proof
of $(*)_1$.
If the conclusion fails we let $c = g(c^*_2,d_1,d_2)$ so \nl
$c(Y^*) = g(c^*_2(Y^*),d_1(Y^*),d_2(Y^*)) = c^*_2(Y^*)$ as $d_1(Y^*)
\ne d_2(Y^*) +$ choice of $g$ and $c(Z) = g(c^*_2(Z),d_1(Z),d_2(Z)) =
d_1(Z) \ne c^*_2(Z)$ as $d_1(Z) = d_2(Z) \ne c^*_2(Z)$. \nl
So $c$ contradicts $(*)_1$.]
\sn
\item "{$(*)_4$}"   for $Z \in {\Cal P}$, there are no $d_1,d_2 \in
{\frak C}$ such that $d_1(Z) = d_2(Z),d_1(Y^*) \ne d_2(Y^*)$ except
possibly when $\{d_1(Z)\} = \{c^*_2(Z)\} \in \{Z \cap Y^*,Z \backslash
Y^*\}$. \nl
[Why?  If $d_1(Z) \ne c^*_2(Z)$ use $(*)_3$, if $d_1(Z) = c^*_2(Z)$
and there is $\pi \in \text{ Per}(X)$ satisfying $\pi(Y^*) = Y^*,\pi(Z) =
Z,\pi(c^*_2(Z)) \ne c^*_2(Z)$ we use it to get the situation in
$(*)_3$.]
\ermn
Let

$$
K = \{(m):\text{ for some } Z \in {\Cal P} \text{ we have }
|Z \cap Y^*| = m\}
$$
\mn
we are assuming $K \ne \emptyset$. By $(*)_4 +$ symmetry we know
\mr
\item "{$(*)_5$}"   if $(m) \in K,1 \ne m < k-1,c_1,c_2 \in {\frak
C}$ and $Z_1,Z_2 \in \binom Xk$ satisfying $c_1(Z_1) = c_2(Z_1)$ and $|Z_1 \cap Z_2| = m$,
\ub{then} $c_1(Z_2) = c_2(Z_2)$.
\ermn
\ub{Case 1}:  There is $(m) \in K$ such that $1 \ne m < k-1$, let
${\Cal P}' = {\Cal P} \cup \{Y^*\}$.

For any $c_1,c_2 \in {\frak C}$ let ${\Cal P}_{c_1,c_2} = \{Y \in
\binom Yk:c_1(Y) = c_2(Y)\}$. \nl
By $(*)_5$ we have $[Y_1,Y_2 \in \binom Xk \wedge |Y_1 \cap Y_2| = m
\Rightarrow [Y_1 \in {\Cal P}_{c_1,c_2} \equiv Y_2 \in {\Cal
P}_{c_1,c_2}]]$.

Let $Y_1 \in \binom Xk,c_1 \in {\frak C}$, let $a = c_1(Y_1)$ let $Y_2
\in \binom Xk$ be such that $\{a,b\} = Y_1 \backslash Y_2$ for some $b
\ne a$.  By
conjugation there is $c_2 \in {\frak C}$ such that $c_2(Y_1) = a =
c_1(Y_1) \and c_1(Y_2) \ne c_2(Y_2)$.  To ${\Cal P}_{c_1,c_2}$ apply 
\scite{12.3} below; so necessarily $|X| = 2k,m=0$.  But as $m=0,(m)
\in K$ there is $Y \in {\Cal P},|Y \cap Y^*| = m=0$ hence $Y = X
\backslash Y^*$, and by $(*)_2(\alpha)$ we get a contradiction,
i.e. we can find a $\pi$ contradicting it. 
\bn
\ub{Case 2}:  $(m) \in K,m = k-1$ and not Case 1.

Let $Z \in {\Cal P},|Z \cap Y^*| = k-1$ so by $(*)_4$ and 
${\frak C}$ being symmetric
\mr
\item "{$(*)_6$}"  if $Z_1,Z_2 \in \binom Xk,|Z_1 \cap Z_2| =
k-1,d_1,d_2 \in {\frak C},d_1(Z_1) = d_2(Z_1),d_1(Z_2) \ne d_2(Z_2)$
then $\{d_1(Z_1)\} = Z_1 \backslash Z_2$. \nl
Also
\sn
\item "{$(*)_7$}"  if $Z_1,Z_2 \in \binom Xk,|Z_1 \cap Z_2| =
k-1$ then for no $d \in {\frak C}$ do we have 
$d(Z_1) \ne d(Z_2) \and \{d(Z_1),d(Z_2)\} \subseteq Z_1 \cap Z_2$.
\nl
[Why?  Applying $\pi \in \text{ Per}(X)$ we get a contradiction to
$(*)_6$.]
\ermn
This case (2) is finished by the following claim (and then we shall continue).
\enddemo
\bigskip

\proclaim{\stag{12.2a} Claim}  Assume $(a)^*(i)$ 
of \scite{12.2} and (b),(c) of \scite{12.1}
and $(*)_7$ above (on ${\frak C}$).  \ub{Then} ${\frak C}$ is full.
\endproclaim
\bigskip

\demo{Proof}  Now
\mr
\item "{$(*)_8$}"  for every $Z_1,Z_2 \in \binom Xk,|Z_1 \cap Z_2| =
k-1$ and $a \in Z_1 \cap Z_2$ there is no $d \in {\frak C}$ such that
$d(Z_1) = d(Z_2) = a$. 
\ermn
Why?  Otherwise we can find $Z_1,Z_2$ such that $|Z_1 \cap Z_2| =
k-1,d(Z_1) = d(Z_2) = a$ 
hence for every $Z_1,Z_2 \in \binom Xk$ such
that $|Z_1 \cap Z_2| = k-1$ and $a \in Z_1 \cap Z_2$ there is such a
$d$ (using $\pi \in \text{ Per}(X)$).
\sn
Let $Z_1,Z_2 \in \binom Xk$ such that $|Z_1 \cap Z_2| = k-1$.  Let $x
\ne y \in Z_1 \cap Z_2$.  Choose $d_1 \in {\frak C}$ such that
$d_1(Z_1) = d_1(Z_2) = x$. \nl
Choose $d_2 \in {\frak C}$ such that $d_2(Z_1) = d_2(Z_2) = y$. \nl
Choose $d_3 \in {\frak C}$ such that $d_3(Z_1) = y,d_3(Z_2) \in Z_2 \backslash Z_1$.
\sn
\ub{Why is it possible to choose $d_3$}?  (Using $\pi \in \text{ Per}(X)$), otherwise (using
$(*)_7$) we have
\mr
\item "{$\bigotimes$}"  if $Y_1,Y_2 \in \binom Xk,|Y_1 \cap Y_2| = k-1$
\nl
$d \in {\frak C},d(Y_1) \in Y_1 \cap Y_2$ then $d(Y_2) \in Y_1 \cap
Y_2$ \nl
hence by $(*)_7,d(Y_2) = d(Y_1)$ \nl
so for $d \in {\frak C}$ we have (by a chain of $Y$'s)

$$
Y_1,Y_2 \in \binom Xk,d(Y_1) \in Y_1 \cap Y_2 \Rightarrow d(Y_2) =
d(Y_1).
$$
\ermn
Let $c \in {\frak C},Y_1 \in \binom Xk,x_1 = c(Y_1)$.  Let $x_2 \in X
\backslash Y_1,Y_2 = Y_1 \cup \{x_2\} \backslash \{x_1\}$, so if
$c(Y_2) \in Y_1 \cap Y_2$ we get a contradiction, so $d(Y_2) = x_2$. \nl
Let $x_3 \in Y_1 \cap Y_2,Y_3 = Y_1 \cup Y_2 \backslash \{x_3\}$ so
$Y_3 \in \binom Xk,|Y_3 \cap Y_1| = k-1 = |Y_3 \cap Y_2|$ and clearly
$c(Y_1),c(Y_2) \in Y_3$. \nl
If $c(Y_3) \notin Y_1$ then $Y_3,Y_1$ contradict $\bigotimes$.  If
$c(Y_3) \notin Y_2$ then $Y_3,Y_2$ contradict $\bigotimes$.  But $c(Y_3)
\in Y_3 \subseteq Y_1 \cup Y_2$ contradiction.  So $d_3$ exists.

We shall use $d_1,d_2,d_3,Z_1,Z_2$ to get a contradiction (thus
proving $(*)_8$). \nl
Let $\{z\} = Z_2 \backslash Z_1$ so $\langle x,y,z \rangle$ is with no repetitions. \nl
Let $d = g(d_1,d_2,d_3)$ so with $g=g^*$ or $g=g_{<x,y,z>}$

$$
\align
d(Z_1) = &g(d_1(Z_1),d_2(Z_1),d_3(Z_1)) = g(x,y,y) = y \\
  &\text{ (see Definition of } g)
\endalign
$$

$$
\align
d(Z_2) = &g(d_1(Z_2),d_2(Z_2),d_3(Z_2)) \\
  &= g(x,y,z) = x \\
  &\text{ by Definition of } g \text{ as } y \ne z \text{ because } y \in Z_1,z \notin Z_1.
\endalign
$$
\mn
So $Z_1,Z_2,d$ contradicts $(*)_7$. \nl
So we have proved $(*)_8$.
\mr
\item "{$(*)_9$}"  if $|Z_1 \cap Z_1| = k-1,Z_1,Z_2 \in \binom Xk,d
\in {\frak C},d(Z_1) \in Z_1 \cap Z_2$, \ub{then} $d(Z_2) \in Z_2
\backslash Z_1$. \nl
[Why?  By $(*)_7, d(Z_2) \notin Z_1 \cap Z_2 \backslash \{d(Z_1)\}$
and by $(*)_8, d(Z_2) \notin \{d(Z_1)\}$.]
\ermn
Now we shall finish proving the claim \scite{12.2a}.
\enddemo
\bn
Let $c \in {\frak C}$. \nl

Now let $x_1,x_2 \in X$ be distinct and $Y \subseteq X
\backslash \{x_1,x_2\},|Y| = k$.  Let $x_3 = c(Y),x_4 \in Y \backslash
\{x_3\}$ and $x_5 \in Y \backslash \{x_3,x_4\}$. \nl
So $Y_1 = Y \cup \{x_1\} \backslash \{x_4\}$ belong to $\binom Xk$
satisfies $|Y_1 \cap Y| = k-1$ and $c(Y) = x_3 \in Y_1 \cap Y$ hence by $(*)_9$, we have
$c(Y_1) = x_1$. \nl
Let $Y_2 = Y \cup \{x_2\} \backslash \{x_4\}$ so similarly $c(Y_2) =
x_2$. \nl
Let $Y_3 = Y \cup \{x_1,x_2\} \backslash \{x_4,x_5\}$, so $Y_3 \in
\binom Xk$ and $Y_3 \backslash Y_1 = \{x_2\}$ and $Y_3 \backslash Y_2
= \{x_1\}$. \nl
\ub{If} $c(Y_3) \in Y$, \ub{then}

$$
c(Y_3) \in Y_3 \cap Y = Y \backslash \{x_4,x_5\} \subseteq Y_1 \text{
hence } c(Y_3) \in Y_3 \cap Y_1
$$
\mn
recall

$$
c(Y_1) = x_1 \in Y_3 \cap Y_1
$$
\mn
and $c(Y_3) \ne x_1$ as $x_1 \notin Y$ so $(Y_3,Y_1,c)$ contradicts $(*)_7$. \nl
\ub{If} $c(Y_3) = x_1$, \ub{then} $c,Y_3,Y_1$ contradicts $(*)_8$. \nl
If $c(Y_3) = x_2$, \ub{then} $c,Y_3,Y_2$ contradicts $(*)_8$. \nl
Together contradiction, so we have finished proving \scite{12.2a} hence
Case 2 in the proof of \scite{12.1}.  \hfill$\square_{\scite{12.2a}}$\margincite{12.2a}
\bn
\ub{Continuation of the proof of \scite{12.1}}:
\mn
\ub{Case 3}:  Neither case 2 nor case 3.

As ${\Cal P} \ne \emptyset$ (otherwise we are done) clearly $K =
\{(1)\}$.  So easily (clearly $2k-1 \le |X|$ as $(1) \in K$) and
\mr
\item "{$\boxtimes_1$}"  if $|Y_1 \cap Y_2| = 1,Y_1 \in \binom Xk,Y_2
\in \binom Xk$ and $d \in {\frak C}$ then $d(Y_1) \in Y_1 \cap Y_2
\vee d(Y_2) \in Y_1 \cap Y_2$. \nl
[Why?  Otherwise by conjugation we can get a contradiction to $(*)_4$ above.]
\sn
\item "{$\boxtimes_2$}"   $Y_1,Y_2 \in \binom Xk,|Y_1 \cap Y_2| = k-1,d \in {\frak
C},d(Y_1),d(Y_2) \in Y_1 \cap Y_2$ is impossible. \nl
[Why?  Let $x \in Y_1 \backslash Y_2$ 
and $y \in Y_2 \backslash Y_1$, we can find $Y_3
\in \binom Xk$ such that $Y_3 \cap (Y_1 \cup Y_2) = \{x,y\}$ so $Y_3
\cap Y_1 = \{x\},Y_3 \cap Y_2 = \{y\}$; this is possible as $|X| \ge 2k-1$.  
Apply $\boxtimes_1$ to $Y_3,Y_1,d$ and as 
$d(Y_1) \ne x$ (as $d(x_1) \in Y_Z$) we have $c(Y_3) = x$. \nl
Apply $\boxtimes_1$ to $Y_3,Y_2,d$ and as $d(Y_2) \ne y$ (as $d(Y) \in
Y_1$) we get $d(Y_3) = y$. \nl
But $x \ne y$, contradiction.]
\ermn
By $\boxtimes_2$ we can use the proof of case 2 from $(*)_7$,
i.e. Claim \scite{12.2a} to get
contradiction.  \hfill$\square_{\scite{12.2}}$\margincite{12.2}
\bigskip

\proclaim{\stag{12.3} Claim}  Assume
\mr
\item "{$(a)$}"   $k^* < k < |X| < \aleph_0$,
\sn
\item "{$(b)$}"   ${\Cal P} \subseteq \binom Xk$
\sn
\item "{$(c)$}"  if $Z,Y \in \binom Xk,|Z \cap Y| = k^*$ then $Z \in
{\Cal P} \Leftrightarrow Y \in {\Cal P}$. 
\sn
\item "{$(d)$}"  $2k-k^* \le |X|$ (this is equivalent to clause (c)
being non empty).
\ermn
\ub{Then}
\mr
\item "{$(\alpha)$}"  ${\Cal P} = \emptyset \vee {\Cal P} = \binom
Xk$ \ub{or}
\sn
\item "{$(\beta)$}"  $|X| = 2k,k^* = 0$ and so $E = E_{X,k} =:
\{(Y_1,Y_2):Y_1 \in \binom Xk,Y_2 \in \binom Xk,(Y_1 \cup Y_2 = X)\}$ 
is an equivalence relation on $X$, with each equivalence class is a
doubleton and ${\Cal P}$ is a union of a set of $E$-equivalence classes.
\endroster
\endproclaim
\bigskip

\demo{Proof}  If not clause $(\alpha)$, then for some $Z_1 \in {\Cal P},Z_2 \in \binom
Xk \backslash {\Cal P}$ we have $|Z_1 \backslash Z_2| = 1$. \nl
Let $Z_1 \backslash Z_2 = \{a^*\},Z_2 \backslash Z_1 = \{b^*\}$.
\enddemo
\bn
\ub{Case 1}:  $2k - k^* < |X|$.

We can find a set $Y^+ \subseteq X \backslash (Z_1 \cup Z_1)$ with
$k-k^*$ members (use $|Z_1 \cup Z_2| = k+1,|X \backslash (Z_1 \cup
Z_2)| = |X| - (k+1) \ge (2k-k^* +1) - (k+1) = k-k^*$).
\sn
Let $Y^- \subseteq Z_1 \cap Z_2$ be such that $|Y^-| = k^*$. \nl
Let $Z = Y^- \cup Y^+$ so $Z \in \binom Xk,|Z \cap Z_1| = |Y^-| =
k^*,|Z \cap Z_2| = |Y^-| = k^*$ hence $Z_1 \in {\Cal P}
\leftrightarrow Z \in {\Cal P} \leftrightarrow Z_2 \in {\Cal P}$,
contradiction.
\bn
\ub{Case 2}:  $2k - k^* = |X|$ and $k^* > 0$.

Let $Y^+ = X \backslash (Z_1 \cup Z_2)$ so

$$
|Y^+| = (2k - k^*) - |k+1| = k - k^* - 1.
$$
\mn
Let $Y^- \subseteq Z_1 \cap Z_2$ be such that $|Y^-| = k^*-1$ (O.K. as
$|Z_1 \cap Z_2| = k-1 \ge k^*)$. \nl
Let $Z = Y^+ \cup Y^- \cup \{a^*,b^*\}$.  So $|Z| = (k - k^*-1) +
(k^*-1) + 2 = k,|Z_1 \cap Z| = |Y^- \cup \{a^*\}| = k^*,|Z_2 \cap Z| =
|Y^- \cup \{b^*\}| = k^*$ and as in case 1 we are done.  \hfill$\square_{\scite{12.3}}$\margincite{12.3}
\bigskip

\proclaim{\stag{12.4} Claim}  If $r({\Cal F}) < \infty$ \ub{then}
\scite{12.1} or \scite{12.2} apply so ${\frak C}$ is full.
\endproclaim
\bigskip

\remark{Remark}  Recall $r({\Cal F}) = \text{ Inf}\{r:\text{some } f
\in {\Cal F}_{[r]}$ is not a monarchy$\}$, see Definition \scite{2.3}.
\endremark
\bigskip

\demo{Proof}
\sn
\ub{Case 1}:  $r({\Cal F}) \ge 4$. \nl

Let $f \in {\Cal F}_{[r]}$ exemplify it, so by \scite{2.5} for some $\ell(*)$:

$$
\bar a \in {}^r X \text{ with repetitions } \Rightarrow f(\bar a) =
a_{\ell(*)}.
$$
\mn
As $f$ is not a monarchy for some $k(*) \in \{1,\dotsc,r\}$ and 
$\bar a^* \in {}^r X,f(\bar a^*) = a_{k(*)} \ne a_{\ell(*)}$.

Without loss of generality $\ell(*) = 1,k(*) = 2$ and \scite{12.1}
apply.
\bn
\ub{Case 2}:  $r({\Cal F}) = 3$. \nl

Let $f \in {\Cal F}_{[r]}$ exemplify it.
Now apply \scite{7.9B}; if $(a)$ there holds, apply \scite{12.1}, if
$(b)$ or $(c)$ there holds, apply \scite{12.2}.
\bn
\ub{Case 3}:  $r({\Cal F}) = 2$.

By \scite{13.7} below, clause $(a)$ of \scite{12.1} holds so we are
done.  \hfill$\square_{\scite{12.4}}$\margincite{12.4} 
\enddemo
\newpage

\head {\S4 The case $r = 2$} \endhead  \resetall \sectno=4
\bigskip

This is revisted in \S6 (non simple case), and we can make
presentation simpler (e.g. \scite{15.3}).
\demo{\stag{13.1} Hypothesis}
\mr
\item "{$(a)$}"  $r({\Cal F}) = 2$
\sn
\item "{$(b)$}"  $|X| \ge 5$ (have not looked at 4).
\endroster
\enddemo
\bigskip

\proclaim{\stag{13.2} Claim}  Choose $\bar a^* = \langle a^*_1,a^*_2
\rangle,a^*_1 \ne a^*_2 \in X$.
\endproclaim
\bigskip

\proclaim{\stag{13.3} Claim}  For some $f \in {\Cal F}$ and $\bar b \in
{}^2 X$ we have
\mr
\item "{$(a)$}"  $f(\bar a^*) = a^*_2$
\sn
\item "{$(b)$}"  $\bar a^* \char 94 \bar b$ is with no repetition
\sn
\item "{$(c)$}"  $f(\bar b) = b_1 \ne b_2$.
\endroster
\endproclaim
\bigskip

\demo{Proof}  There is $f \in {\Cal F}$ not monarchial so for some
$\bar b,\bar c \in {}^2 X$

$$
f(\bar b) = b_1 \ne b_2,f(\bar c) = c_2 \ne c_1.
$$
\ub{If} Rang$(\bar b) \cap \text{ Rang}(\bar c) = \emptyset$ we can
conjugate $\bar c$ to $\bar a^*,f$ to $f'$ which is as required. \nl
\ub{If} not, find $\bar d \in {}^2 X,d_1 \ne d_2$ satisfying $\text{Rang}(\bar d)
\cap \text{ (Rang}(\bar a) \cup \text{ Rang}(\bar b)) = \emptyset$ so
$\bar d,\bar b$ \ub{or} $\bar d,\bar c$ are like $\bar c,\bar b$ or
$\bar b,\bar c$ respectively.  \hfill$\square_{\scite{13.3}}$\margincite{13.3}
\enddemo
\bigskip

\proclaim{\stag{13.4} Claim}  There is $f^* \in {\Cal F}_{[2]}$ such
that
\mr
\item "{$(a)$}"  $f^*(\bar a^*) = \bar a^*_2$
\sn
\item "{$(b)$}"  $b_1 \ne b_2 \in X,\{b_1,b_2\} \subseteq
\{a^*_1,a^*_2\} \Rightarrow f(b_1,b_2) = b_2$
\sn
\item "{$(c)$}"  $b_1 \ne b_2,\{b_1,b_2\} \nsubseteq \{a^*_1,a^*_2\}
\Rightarrow f(b_1,b_2) = b_1$.
\endroster
\endproclaim
\bigskip

\demo{Proof}  Choose $f$ such that
\mr
\item "{$(a)$}"  $f \in {\Cal F}_{[2]}$
\sn
\item "{$(b)$}"  $f(\bar a^*) = a^*_2$
\sn
\item "{$(c)$}"  $n(f) = |\{\bar b \in {}^2 X:f(\bar b) = b_1\}|$ is
maximal under $(a) + (b)$.
\ermn
Let ${\Cal P} = \{\bar b \in {}^2 X:f(\bar b) = b_1\}$.  In each case
we can assume that the previous cases do not hold for any $f$
satisfying (a), (b), (c).
\enddemo
\bn
\ub{Case 1}:  There is $\bar b \in {}^2(X \backslash \{a^*_1,a^*_2\})$
such that $f(\bar b) = b_2 \ne b_1$. \nl
There is $g \in {\Cal F}_{[2]},g(\bar a^*) = a^*_2,g(\bar b) = b_1$ (by \scite{13.3} +
conjugation).
Let $f^+(x,y) = f(x,g(x,y))$. \nl
So
\mr
\item "{$(A)$}"  $f^+(\bar a^*) = f(a^*_1,g(\bar a^*)) =
f(a^*_1,a^*_2) = a^*_2$
\sn
\item "{$(B)$}"  $f^+(\bar b) = f(b_1,g(\bar b)) = f(b_1,b_1) = b_1$
\sn
\item "{$(C)$}"  if $\bar c \in {\Cal P}$ then $f(\bar c) = c_1$.
\ermn
[Why does (C) hold?  If 
$g(\bar c) = c_1$ then $f^+(\bar c) = f(c_1,g(\bar c)) = f(c_1,c_1)
= c_1$. \nl
If $g(\bar c) = c_2$ then $f^+(\bar c) = f(c_1,g(\bar c)) = f(c_1,c_2)
= f(\bar c) = c_1$ (last equality as $\bar c \in {\Cal P}$).] \nl
By the choice of $f$ the existence of $f^+$ is impossible so
\mr
\item "{$(*)$}"  $\bar b \in {}^2(X \backslash \{a^*_1,a^*_2\})
\Rightarrow f(\bar b) = b_1 \Rightarrow \bar b \in {\Cal P}$ (if $b_1
= b_2$ - trivial).
\endroster
\bn
\ub{Case 2}:  There are $b_1 \ne b_2$ such that
$\{b_1,b_2\} \nsubseteq \{a^*_1,a^*_2\},f(b_1,b_2) = b_2$ and $b_1 \ne
a^*_1 \wedge b_2 \ne a^*_2$.

There is $g \in {\Cal F}_{[2]}$ such that $g(a^*_1,a^*_2) =
a^*_2,g(b_1,b_2) = b_1$. \nl
[Why?  There is $\pi \in \text{ Per}(X),\pi(b_1) = a^*_1,\pi(b_2) =
a^*_2,\pi^{-1}(\{b_1,b_2\})$ is disjoint to $\{a^*_1,a^*_2\}$.
Conjugate $f$ by $\pi^{-1}$, getting $g$ so
$g(a^*_1,a^*_2) = g(\pi b_1,\pi b_2) = \pi(f(b_1,b_2)) = \pi(b_2) =
a^*_2$; let $c_1,c_2$ be such that $\pi(c_1) = b_1,\pi(c_2) = b_2$ so

$$
g(b_1,b_2) = g(\pi c_1,\pi c_2) = \pi(f(c_1,c_2)) = \pi(c_1) = b_1
$$
\mn
(third equality as $c_1,c_2 \notin \{a^*_2,a^*_2\}$ by not Case 1).
So there is such $g \in {\Cal F}$.]
\sn
Let $f^+(x,y) = f(x,g(x,y))$, as before $f^+$ contradicts the choice
of $f$. 
\bn
\ub{Case 3}:  For some $b' \ne b'' \in X \backslash \{a^*_1,a^*_2\}$
we have $f(a^*_1,b') = b' \wedge f(a^*_1,b'') = a^*_1$.

As in Case 2, using $\pi \in \text{ Per}(X)$ such that $\pi(a^*_1) =
a^*_1,\pi(a^*_2) = a^*_2,\pi(b') = b''$.
\bn
\ub{Case 4}:  For some $b' \ne b'' \in X \backslash \{a^*_1,a^*_2\}$
we have $f(b',a^*_2) = a^*_2 \wedge f(b'',a^*_2) = b''$.

As in Case 3; recall that  \wilog \, Case 1,2,3,4 fails.
\bn
\ub{Case 5}:  For some $b',b'' \in X \backslash \{a^*_1,a^*_2\}$
we have $f(a^*_1,b') = b' \wedge f(b'',a^*_2) = a^*_2$.

As Cases 1,2,3,4 fail, \wilog \, $b' \ne b''$ and prove as in Case 2
conjugating by 
$\pi \in \text{ Per}(X)$ such that $\pi(b') = a^*_2,\pi(a^*_1) =
a^*_1$ and $\pi(b'') = b''$ getting $g$ which satisfies
$g(\bar a^*_1,a^*_2) = a^*_2$ and $g(b',a^*_2) = \pi(f(b'',b')) = \pi
(b'') = b'$ whereas $f(b',a^*_2) = a^*_2$; so $f^+(x,y) = f(x,g(x,y))$
contradicts the choice of $f$.  \nl
Without loss of generality, Cases 1-5 fail.
\bn
\ub{Case 6}:  For some $b \in X \backslash \{a^*_1,a^*_2\}$
we have $f(a^*_1,b) = b$.

So as Cases 1-5 fail we have
\mr
\item "{$\circledast$}"  $\forall b_1,b_2 \in X[f(b_1,b_2) \ne b_1
\leftrightarrow (b_1 = a^*_1 \and b_2 \ne a^*_1)]$
\ermn
(except in the case $f(a^*_2,a^*_1) = a^*_1$ but then let $\pi \in
\text{ Per}(X),\pi(a^*_1) = a^*_2,\pi(a^*_2) = a^*_1$ (and $\pi(a) = a$
for $a \in X \backslash \{a^*_1,a^*_2\}$), then $g = \pi f
\pi^{-1}$ satisfies $g(a^*_1,a^*_2) = a^*_2$ but for $b \in X
\backslash \{a^*_1,a^*_2\},g(a^*_1,b) = a^*_1$, easy contradiction (or as below)). \nl
Hence for every $c \in X$ there is $f_c \in {\Cal F}_{[2]}$ such that
\mr
\item "{$\circledast_{f_c}$}"  $\forall b_1 b_2 \in X[f_c(b_1,b_2) \ne b_1
\leftrightarrow (b_1 = c \and b_2 \ne c)]$.
\ermn
Let $a \ne c$ be from $X$ and define $f_{a,c} \in {\Cal F}_{[2]}$ by
$f_{a,c}(x,y) = f_a(x,f_c(y,x))$. \nl
Assume $b_1 \ne b_2$ so $f^*_{a,c}(b_1,b_2) = b_2 \ne b_1$ implies
$f_c(b_2,b_1) \in \{b_1,b_2\},f_{a,c}(b_1,b_2) =
f_a(b_1,f_c(b_2,b_1))$ and so (by the choice of $f_a$) $b_1 = a \and
f_c(b_2,b_1) = b_2$ which (by the choice of $f_c$) implies $b_1 = a
\and b_2 \ne c$.  But $b_1 = a \and b_2 \ne c \and b_1 \ne b_2$
implies $f_c(b_2,b_1) = b_2,f_{a,c}(b_1,b_2) = f_a(b_1,b_2) = b_2$.
So $f_{a,c}(b_1,b_2) = b_2 \ne b_1 \Leftrightarrow b_1 = a \and b_2 \ne
c \and b_2 \ne b_1$.
\sn
Let $a = a^*_1$.  Let $\langle c_i:i < i^* = |X|-2 \rangle$ list $X \backslash
\{a^*_1,a^*_2\}$.  We define by induction on $i \le i^*$, a
function $f_i \in {\Cal F}_{[2]}$ by

$$
f_0 (x,y) = y
$$

$$
f_{i+1}(x,y) = f_i(x,f_{a,c_i}(x,y))
$$
\mn
and let $f' = f_{i^*}$.  Now by induction on $i$ we can show
$f_i(a^*_1,a^*_2) = a^*_2$ and
$f'(b_1,b_2) = b_2 \ne b_1 \Rightarrow (\forall i <
i^*)(f_{a,c_i}(b_1,b_2) = b_2 \ne b_1)$. \nl
So $f' \in {\Cal F}_{[2]},f'(a^*_1,a^*_2) = a^*_2$ and $b_1 \ne b_2
\wedge (b_1,b_2) \ne (a^*_1,a^*_2),(a^*_2,a^*_1)$ implies
$f'(b_1,b_2) = b_1$.  By the choice of $f$ (minimal $n(f)$) we get a contradiction.
\bn
\ub{Case 7}:  For some $b \in X \backslash \{a^*_1,a^*_2\}$
we have $f(b,a^*_2) = a^*_2$.
\sn
Similar to Case 6.

That is, as there,
\wilog \, for every $a \in X$ for some $f_a \in {\Cal F}_{[2]}$ we have
\mr
\item "{$\circledast$}"  $(\forall b_1,b_2 \in X)[(f(b_1,b_2) = b_2
\ne b_1 \leftrightarrow b_2 = a \ne b_1)]$.
\ermn
Let $a \ne c \in X$ let $f_{a,c}(x,y) = f_c(f_a(y,x),y)$.

So for $b_1 \ne b_2 \in X$
\mr
\item "{$(i)$}"  $f_{a,c}(b_1,b_2) = b_2(\ne b_1)$ implies
$f_a(f_c(b_2,b_1),b_2) = b_1$ which implies $b_2 = c \and f_c(b_2,b_1)
= b_2$ which implies $b_2 = c \and b_1 \ne a$.
\ermn
We continue as there.
\bn
\ub{Case 8}:  Not Cases 1-7; not the conclusion. 

So for $\bar a = (a_1,a_2) = {}^2X,a_1 \ne a_2$ there is $f_{\bar a} \in {\Cal F}$
such that

$$
\{b_1,b_2\} \nsubseteq \{a_1,a_2\} \Rightarrow f_{\bar a}(b_1,b_2) = b_1
$$

$$
f_{\bar a}(a_1,a_2) = a_2
$$
\mn
and (as ``not the conclusion")

$$
f_{\bar a}(a_2,a_1) = a_2
$$
\mn
Let $\langle \bar b^i:i < i^* = |X|^2 - |X| - 2 \rangle$ list the pairs
$\bar b = (b_1,b_2) \in {}^2 X$ such that $b_1 \ne b_2,\{b_1,b_2\} \ne \{a^*_1,a^*_2\}$. \nl
Define $g_i \in {\Cal F}_{[2]}$ by induction on $i$. \nl
Let $g_0(x,y) = x$. \nl
Let $g_{i+1}(x,y) = f_{{\bar b}^i}(g_i(x,y),y)$. \nl
We can prove by induction on $i \le i^*$ that: $g_i(a^*_1,a^*_2) =
a^*_1,g_i(a^*_2,a^*_1) = a^*_2,j < i \Rightarrow g_i(\bar
b^j) = b^j_2$.  So $g_{i^*}$ is as required interchanging 1 and 2 that is
$g(x,y) =: g_{i^*}(y,x)$ is as required.  \hfill$\square_{\scite{13.4}}$\margincite{13.4}
\bigskip

\definition{\stag{13.5} Definition/Choice}  For $b \ne c \in X$ let
$f_{b,c}$ be like $f$ in \scite{13.4} with $(b,c)$ instead of
$(a^*_1,a^*_2)$, so $f_{c,b}(c,b)$ is $b$ and $f(b,c) = c$ and
$f(x_1,_2) = x_1$ if $\{x_1,x_2\} \nsubseteq \{b,c\}$.
\enddefinition
\bigskip

\proclaim{\stag{13.6} Claim}  Let $a_1,a_2,a_3 \in X$ be pairwise
distinct.

Then for some $g \in {\Cal F}_{[3]}$:
\mr
\item "{$(i)$}"  $\bar b \in {}^3 X$ with repetitions $\Rightarrow
g(\bar b) = b_1$,
\sn
\item "{$(ii)$}"  $g(a_1,a_2,a_3) = a_2$.
\endroster
\endproclaim
\bigskip

\demo{Proof}  Without loss of generality we replace $a_2$ by $a_3$ in
(ii).
\sn
Let $h_\ell$ for $\ell = 1,2,3,4$ be the three place functions

$$
h_1(\bar x) = f_{a_1,a_2}(x_1,x_2)
$$

$$
h_2(\bar x) = f_{a_1,a_3}(x_1,x_3)
$$

$$
h_3(\bar x) = f_{a_2,a_3}(h_1(\bar x),h_2,(\bar x))
$$

$$
h_4(\bar x) = f_{a_1,a_3}(x,h_3(\bar x)).
$$
\mn
Clearly $h_1,h_2,h_3,h_4 \in {\Cal F}_{[3]}$.  We shall show $h_4$ is
as required.
\nl
To prove clause (ii) note that for $\bar a = (a_1,a_2,a_3)$ we have
$h_1(\bar a) = a_2,h_2(\bar a) = a_3,h_3(\bar a) =
f_{a_2,a_3}(a_2,a_3) = a_3$ and $h_4(\bar a) = f_{a_1,a_3}(a_1,a_3) =
a_3$ as agreed above.  To prove clause (i). \nl
Let $\bar b \in {}^3X$ be without repetitions such that $\bar b \ne
\bar a$.
\bn
\ub{Case 1}:  $b_1 \ne a_1,a_3$ so

$$
h_4(\bar b) = f_{a_1,a_3}(b_1,h_3(\bar b)) = b_1 \text{ as } b_1 \ne a_1,a_3.
$$
\bn
\ub{Case 2}:  $b_1 = a_1,b_2 \ne a_2$ hence $b_1 \ne a_2,a_3$, so

$$
h_1(\bar b) = f_{a_1,a_2}(b_1,b_2) = b_1 \text{ (if } b_2 = a_1 
\text{ also O.K)}
$$

$$
h_3(\bar b) = f_{a_2,a_3}(h_1(\bar b),h_2(\bar b)) =
f_{a_2,a_3}(b_1,h_2(\bar b)) = b_1 \text{ as } b_1 \ne a_2,a_3
$$

$$
h_4(\bar b) = f_{a_1,a_3}(b_1,h_3(\bar b)) = h_{a_1,a_3}(b_1,b_1) = b_1.
$$
\bn
\ub{Case 3}:  $b_1 = a_1,b_2 = a_2,b_3 \ne a_3$, so

$$
h_1(\bar b) = f_{a_1,a_2}(b_1,b_2) = b_2
$$

$$
h_2(\bar b) = f_{a_1,a_3}(b_1,b_3) = f_{a_1,a_3}(a_1,b_3) = a_1 =
b_1\text{ as } b_3 \ne a_3 \text {(if } b_3 = a_1 \text{ fine)}
$$

$$
h_3(\bar b) = f_{a_2,a_3}(h_1(\bar b),h_2(\bar b)) =
h_{a_2,a_3}(b_2,b_1) = b_2 \text{ as } b_1 = a_1 \ne a_2,a_3.
$$

$$
h_4(\bar b) = f_{a_1,a_3}(b_1,h_3(\bar b)) =
f_{a_1,a_3}(b_1,b_2) = b_1 \text{ as } b_2 = a_2 \ne a_1,a_3.
$$
\bn
\ub{Case 4}:  $b_1 = a_3,b_3 \ne a_1$. \nl
So

$$
h_1(\bar b) = f_{a_1,a_2}(b_1,b_2) = b_1 \text{ as } b_1 = a_3 \ne a_1,a_2
$$

$$
h_2(\bar b) = f_{a_1,a_3}(b_1,b_3) = f_{a_1,a_3}(a_3,b_3) =
a_3 = b_1 \text{ as } b_3 \ne a_1
$$

$$
h_3(\bar b) = f_{a_2,a_3}(h_1(\bar b),h_2(\bar b)) =
f_{a_2,a_3}(b_1,b_1) = b_1
$$

$$
h_4(\bar b) = f_{a_1,a_3}(b_1,f_3(\bar b) =
f_{a_1,a_3}(b_1,b_1) = b_1.
$$
\bn
\ub{Case 5}:  $b_1 = a_3,b_3 = a_1$. 

$$
h_1(\bar b) = f_{a_1,a_2}(b_1,b_2) = b_1 \text{ as } b_1 = a_3 \ne a_1,a_2
$$

$$
h_2(\bar b) = f_{a_1,a_3}(b_1,b_3) = b_3 \text{ as } \{b_1,b_3\} = \{a_1,a_3\}
$$

$$
h_3(\bar b) = f_{a_2,a_3}(h_1(\bar b),h_2(\bar b)) =
f_{a_2,a_3}(b_1,b_3) \equiv b_1 \text{ as } b_3 = a_1 \ne a_2,a_3
$$

$$
h_4(\bar b) = f_{a_1,a_3}(b_1,f_3(\bar b)) =
f_{a_1,a_3}(b_1,b_1) = b_1.
$$
\mn
as required.  \hfill$\square_{\scite{13.6}}$\margincite{13.6}
\enddemo
\bigskip

\proclaim{\stag{13.7} Claim}  Let $\bar a^* = (a^*_1,a^*_2,a^*_3,a^*_4)
\in {}^4 X$ be with no repetitions.  \ub{Then} for some $g \in {\Cal
F}_{[4]}$ we have
\mr
\item "{$(i)$}"  if $\bar b \in {}^4 X$ is with repetitions then
$f(\bar b) = b_1$
\sn
\item "{$(ii)$}"  $g(\bar a^*) = a^*_2$.
\endroster
\endproclaim
\bigskip

\demo{Proof}  For 
any $\bar a \in {}^3 X$ with no repetitions let $f_{\bar a}$  be as in
\scite{13.6} for the sequence $\bar a$.  Let us define with $(\bar x =
(x_1,x_2,x_3,x_4)),g(\bar x) =
g_0(x_1,g_2(x_1,x_2,x_4),g_3(x_1,x_3,x_4))$ with $g_0 = f_{\langle
a_1,a_2,a_3 \rangle},g_2 = f_{\langle a_1,a_2,a_4 \rangle},g_3 =
f_{\langle a_1,a_3,a_4 \rangle}$.  So
\mr
\item "{$(A)$}"  $g(\bar a^*) =
g_0(a^*_1,g_2(a^*_1,a^*_2,a^*_3),g_3(a^*_1,a^*_3,a^*_4)) 
= g_0(a^*_1,a^*_2,a^*_3) = a^*_2$
\sn
\item "{$(B)$}"  if $\bar b \in {}^4 X$ and $\langle b_1,b_2,b_4
\rangle$ is with repetition then $g_2(b_1,b_2,b_4) = b_1$, hence \nl
$g(\bar b) = g_0(b_1,b_1,g_3(b_1,b_3,b_4)) = b_1$
\sn
\item "{$(C)$}"  if $\bar b \in {}^4 X$ and $\langle b_1,b_3,b_4
\rangle$ is with repetition then $g_3(a_1,a_3,a_4) = b_1$, hence \nl
$g(\bar b) = g_0(b_1,g_2(b_1,b_2,b_4),b_1) = b_1$
\sn
\item "{$(D)$}"  $\bar b \in {}^4 X$ is with repetitions, but
neither (B) nor (C) then necessarily $b_2 = b_3$ so 
$\langle b_1,b_2,b_3 \rangle$ is with repetitions, so \nl
$g(\bar b) = g_0(b_1,b_2,b_3) = b_1$. \nl
${{}}$ \hfill$\square_{\sciteu{4.7}}$\sciteuphantom{4.7}\margincite{4.7}
\endroster
\enddemo
\newpage

\head {Part B: Non simple case \\
\S5 Fullness for the non simple case} \endhead  \resetall \sectno=5
\bigskip

\demo{\stag{14.0} Context}  As in \S1: ${\frak C}$ is a $(X,k)$-FDF
${\Cal F} = \cup\{{\Cal F}_{[r]}:r < \infty\}$ and ${\Cal F} =
\{f:f \in \text{AV}({\frak C})\}$, so

$$
\align
{\Cal F}_{[r]} = \{f:&f \text{ is (not necessarily simple) function
written} \\
  &f_Y(x_1,\dotsc,x_r), \text{ for } Y \in \binom Xk,x_1,\dotsc,x_r
\in Y \text{ such that} \\
  &f_Y(x_1,\dotsc,x_r) \in \{x_1,\dotsc,x_r\} \text{ and} \\
  &{\frak C} \text{ is closed under } f, \text{ i.e. if }
c_1,\dotsc,c_r \in {\frak C} \\
  &\text{ and } c \text{ is defined by } c = f(c_1,\dotsc,c_r) \\
  &\text{i.e. } c(Y) = f_Y(c_1(Y),\dotsc,c_r(Y)) \text{ then } c \in
{\frak C}\}
\endalign
$$
\mn
and we add (otherwise use Part A; alternatively combine the proofs).
\enddemo
\bigskip

\demo{\stag{14.0a} Hypothesis}  If $f \in {\Cal F}$ is simple then it
is a monarchy.
\enddemo
\bigskip

\definition{\stag{14.1} Definition}  1) ${\Cal F}[Y] = \{f_Y:f \in {\Cal F}\}$. \nl
2) ${\Cal F}_{[r]}(Y) = \{f_Y:f \in {\Cal F}_{[r]}\}$.
\enddefinition 
\bigskip

\demo{\stag{14.1a} Observation}  If $f \in {\Cal F}_{[r]},Y \in \binom
Xk$, \ub{then} $f_Y$ is an $r$-place function from $Y$ to $Y$ and
\mr
\item "{$(*)$}"  ${\Cal F}[Y]$ is a clone on $Y$.
\endroster
\enddemo
\bigskip

\definition{\stag{14.2} Definition}  1) $r({\Cal F}) 
= \text{ Min}\{r:r \ge 2$, some $f \in {\Cal F}_{[r]}$ is not a monarchy$\}$
where \nl
2) $f$ is a monarchy \ub{if} for some $t$ we have $\forall Y \forall
x_1,\dotsc,x_r \in Y[f_Y(x_1,\dotsc,x_r) = x_t]$.
\enddefinition 
\bigskip

\proclaim{\stag{14.4} Claim}  1) For proving that ${\frak C}$ is full it is
enough to prove, for some $r \in \{3,\dotsc,k\}$
\mr
\item "{$(*)$}"  for every $Y \in \binom Xk$ and $\bar a \in {}^r Y$
which is one to one there is $f \in {\Cal F}$ such that
{\roster
\itemitem{ $(i)$ }  $f_Y(\bar a) = a_2$
\sn
\itemitem{ $(ii)$ }  if $Z \in \binom Xk,Z \ne Y,\bar b \in {}^r Z$
then $f_Z(\bar b) = b_1$.
\endroster}
\ermn
2) If $ r \ge 4$ we can weaken $f_Z(\bar b) = b_1$ in clause (ii) to
$[b_3 = b_4 \vee b_1 = b_2 \vee b_1 = b_3 \vee b_2 = b_3]
\rightarrow f_Y(\bar b) = b_1$ [not used].
\endproclaim
\bigskip

\demo{Proof}  The proof is as in the proof of \sciteu{5.8} below only
we choose $c_3,c_4,\dotsc,c_5$ such that $\bar a = \langle
c_\ell(Y):\ell = 1,2,\dotsc,r \rangle$ is without repetitions and $f =
f_{\bar a}$ from $(*)$.  \hfill$\square_{\sciteu{5.6}}$\sciteuphantom{5.6}\margincite{5.6}
\enddemo
\bigskip

\proclaim{\stag{14.7} Claim}  In \scite{14.4} we can replace $(*)$ by:
$r=3$ and
\mr
\item "{$(*)$}"  if $Y \in \binom Xk$ and $\bar a \in {}^3 Y$
one-to-one (or just $a_2 \ne a_3$), \ub{then} for some $g \in {\Cal
F}_{[r]}$
{\roster
\itemitem{ $(i)$ }  $g_Y(\bar a) = a_1$
\sn
\itemitem{ $(ii)$ }  if $Z \in \binom Yk,Z \ne Y,\bar b \in {}^3 Z$ is
not one-to-one then $g_Z(\bar b) = b$ if $b_2 = b_3,b_1$ if
otherwise (i.e. $g_{3;1,2}(\bar b))$. 
\endroster}
\endroster
\endproclaim
\bigskip

\demo{Proof}  Like \scite{12.2}.  Let $c^*_1 \in {\frak C},Y^* \in
\binom Xk,a^*_1 = c_2(Y^*),a^*_2 \in Y^* \backslash \{a^*_1\}$ 
we choose $c^*_2$ as in the proof of \scite{14.4}, i.e. \scite{12.1}. 

Let ${\Cal P} = \{Y:Y \in \binom Xk,Y \ne Y^*,c^*_1(Y) \ne
c^*_2(Y)\}$.  As there it suffices to prove that ${\Cal P} =
\emptyset$.  Now
\mr
\item "{$\boxtimes$}"  there are no $Z \in {\Cal P}$ and $d \in {\frak C}$
such that

$$
d(Y^*) = c^*_2(Y^*)
$$

$$
d(Z) \ne c^*_2(Z).
$$
[Why?  If so, let 
$c = g^*(c^*_1,c^*_2,d)$ where $g$ is from 
$(*)$ for $Z,a_1 = c^*_1(Z),a_2 = c^*_2(Z),a_3 = d(Z)$.]
\ermn
Continue as there: the $g_{\bar a}$ depends also on $Y$, and we write
$c(Y) = f_Y(c_1(Y),\dotsc,c_r(Y))$. \nl
${{}}$  \hfill$\square_{\scite{14.7}}$\margincite{14.7}
\enddemo 
\bigskip

\proclaim{\stag{14.20} Claim}  Assume $r({\Cal F}) = 2,({\frak C},{\Cal F}$ as
usual) and
\mr
\item "{$(*)$}"  for every $a_1 \ne a_2 \in Y \in \binom Xk$ for
some $f = f^Y_{\langle a_1,a_2 \rangle} \in {\Cal F}$ we have
{\roster
\itemitem{ $(i)$ }  $f_Y(\bar a) = a_2$
\sn
\itemitem{ $(ii)$ }   $Z \in \binom Yk,Z \ne Y,\bar b \in {}^2 Z
\Rightarrow f_Z(\bar b) = b_1$.
\endroster}
\ermn
\ub{Then} ${\frak C}$ is full.
\endproclaim
\bigskip

\remark{Remark}  ${\frak C}$ is full iff every choice function of
$\binom Xk$ belongs to it.
\endremark
\bigskip

\demo{Proof}  If ${\frak C}$ is not full, as ${\frak C} \ne \emptyset$
there are $c_1 \in {\frak C},c_0 \notin {\frak C},c_0$ a choice
function for $\binom Xk$.  Choose such a pair $(c_1,c_0)$ with $|{\Cal
P}|$ minimal where ${\Cal P} = \{Y \in \binom Xk:c_1(Y) \ne c_0(Y)\}$.
So clearly ${\Cal P}$ is a singleton, say $\{Y\}$.  By symmetry for
some $c_2 \in {\frak C}$ we have $c_2(Y) = c_0(Y)$.  Let $f$ be
$f^Y_{c_1(Y),c_0(Y)} = f^Y_{c_1(Y),c_2(Y)}$ from the assumption so $f
\in {\Cal F}$ and let $c = f(c_1,c_2)$ so clearly $c \in {\frak C}$
(as ${\frak C}$ is closed under every member of ${\Cal F}$).
\sn
Now
\mr
\item "{$(A)$}" $c(Y) = f_Y(c_1(Y),c_2(Y)) = c_2(Y) = c_0(Y)$
\sn
\item "{$(B)$}"  if $Z \in \binom Xk \backslash \{Y\}$ then
$$
c(Z) = f_Z(c_1(Z),c_2(Z)) = c_1(Z) = c_0(Y).
$$
\ermn
So $c = c_0$ hence $c_0 \in {\frak C}$, contradiction. \hfill$\square_{\scite{14.20}}$\margincite{14.20}
\enddemo
\bigskip

\proclaim{\stag{14.21} Claim}  Assume $r({\Cal F}) = 2$ and
$\boxtimes(f^*)$ of \scite{15.7} (see Definitions \sciteu{6.3}, \scite{15.5})
below holds.  \ub{Then} ${\frak C}$ is full.
\endproclaim
\bigskip

\demo{Proof}  We use conventions from \scite{15.5}, \scite{15.5a}, 
\scite{15.7} below.
In $\boxtimes(f^*)$ there are two possibilities.
\mn
\ub{Possibility (i)}:

This holds by \scite{14.20}.
\mn
\ub{Possibility (ii)}:

Similar to the proof of \scite{14.20}.  Again ${\Cal P} = \{Y\}$ where
${\Cal P} = \{Y \in \binom Xk:c_1(Y) \ne c_0(Y)\}$.  We choose $c_2
\in {\frak C}$ such that $c_2(Y) = c_0(Y) \and c_2 (X \backslash Y) = c_0(Y)$,
continue as before.  Why is this possible?  Let $\pi \in \text{
Per}(X)$ be such that $\pi(Y) = Y,\pi(c_1(Y)) = c_0(Y),\pi(c_1(X
\backslash Y)) = c_1(X \backslash Y)$ (and of course $\pi(X \backslash
Y) = X \backslash Y$).  Now conjugating $c_1$ by $\pi$ gives $c_2$ as
required.
\hfill$\square_{\scite{14.21}}$\margincite{14.21}
\enddemo
\bigskip

\proclaim{\stag{14.30} Claim}  If $r({\Cal F}) < \infty$ \ub{then}
${\frak C}$ is full.
\endproclaim
\bigskip

\demo{Proof}  Let $r = r({\Cal F})$.
\mn
\ub{Case 1}:  $r = 2$.

So hypothesis \scite{15.0} below holds.

If $\boxtimes(f)$ of \scite{15.7} holds for some $f \in {\Cal
F}_{[r]}$, by \scite{14.21} we know that ${\frak C}$ is full.  If
$\boxtimes(f)$ of \scite{15.7} fails for every $f \in {\Cal F}_{[r]}$
then hypothesis \scite{15.9} holds hence \scite{15.10}-\scite{15.18}
holds.  So by \scite{15.18} we know that $(*)$ of \scite{14.4} holds
(and ${\Cal P}_\pm$ is a singleton).
So by \scite{14.4}, ${\frak C}$ is full.
\bn
\ub{Case 2}:  $r \ge 4$.

So hypothesis \scite{16.1} holds.  By \scite{16.5} clearly $(*)$ of
\scite{14.4} holds hence by \scite{14.4} we know that ${\frak C}$ is
full.
\bn
\ub{Case 3}:  $r = 3$.

Let $f^* \in {\Cal F}_{[3]}$ be not a monarchy.  So for $\bar b \in
{}^3 Y$ not one-to-one, $Y \in \binom X3$, clearly $f^*_Y(\bar b)$
does not depend on $Y$, so we write $f^-(\bar b)$.  If for some
$\ell(*),f^-(\bar b) = b_{\ell(*)}$ for every such $\bar b$ then easily
\scite{14.4} apply.  If $f^-(\bar b) = g_{r;1,2}(\bar b)$, let $\bar a
\in {}^3 Y,Y \in \binom Xk,\bar a$ is one-to-one, so $f_Y(\bar b) =
a_k$ for some $k$; by permuting the variables, $f^-$ does not 
change while we have $k=1$,
so \scite{14.7} apply.  If both fail, then by repeating the proof of
\scite{7.8} of Part A, for some $f' \in {\Cal F}_{[3]}$, 
for $\bar b$ not one-to-one we have 
$f_{\langle 1,2,1 \rangle}(\bar b)$ or for $\bar b$ not one to one
$f'_Y(\bar b) = f_{\langle 1,2,2 \rangle}(\bar b)$.  By the last
paragraph of the proof of \sciteu{2.8} we can assume the second case
holds.  In this case repeat the proof of the case $\eta = \langle
1,2,2 \rangle$ in the end of the proof of \sciteu{2.8}.
\enddemo
\newpage

\head {\S6 The case $r({\Cal F}) = 2$} \endhead  \resetall \sectno=6
\bn
For this section
\demo{\stag{15.0} Hypothesis} $r=2$.
\enddemo
\bn
\ub{\stag{15.1} Discussion}:  So $(\alpha)$ or $(\beta)$ holds where
\mr
\item "{$(\alpha)$}"   there are $Y \in \binom Xk$ and $f \in
{\Cal F}_{[r]}(Y)$ which is not monarchy.  Hence by \S4, i.e. \scite{13.4} for
$a \ne b \in Y$ there is $f = f^Y_{a,b} \in {\Cal F}_2[Y]$

$$
f_Y(x,y) = \cases y \quad &\text{ if }  \{x,y\} = \{a,b\} \\
x \quad &\text{ if otherwise} \endcases
$$ 
\sn
\item "{$(\beta)$}"  every $f_Y$ is a monarchy but some $f \in {\Cal
F}_{[r]}$ is not.
\endroster
\bigskip

\definition{\stag{15.2} Definition/Choice}  Choose $f^* \in {\Cal F}_2$
such that
\mr
\item "{$(a)$}"   $\neg(\forall Y)(\forall x,y \in Y)(f(x,y)=x)$
\sn
\item "{$(b)$}"   under (a), $n(f) = |\text{dom}_1(f)|$ where
dom$_1(f) = \{(Z,a,b):f_Z(a,b) = a \ne b$ and $Z \in \binom Xk$ and
$\{a,b\} \subseteq Z$ of course$\}$, is maximal.
\endroster
\enddefinition
\bigskip

\demo{\stag{15.3} Fact}  If $f_1,f_2 \in {\Cal F}_{[2]}$ and $f$ is
$f(x,y) = f_1(x,f_2(x,y))$ (formally $f(Y,x,y) = f_1(Y,x,f_2(Y,x,y))$
but we shall be careless) then dom$_1(f) = \text{ dom}_1(f_1) 
\cup \text{ dom}_1(f_2)$.
\enddemo
\bigskip

\demo{Proof}  Easy.
\enddemo
\bigskip

\proclaim{\stag{15.4} Claim}  If $Z \in \binom Xk,f^*_Z(a^*,b^*) = b^*
\ne a$ then
\mr
\item "{$(a)$}"  $(\forall x,y \in Z)[f^*_Z(x,y) = y]$ \ub{or}
\sn
\item "{$(b)$}"  $x,y \in Z \wedge \{x,y\} \nsubseteq \{a^*,b^*\}
\Rightarrow f^*_Z(x,y) = x$.
\endroster
\endproclaim
\bigskip

\demo{Proof}  As in \scite{13.4} (+ \scite{15.2} + \scite{15.3}), recalling that
${\Cal F}(Z)$ is a clone.  \hfill$\square_{\scite{15.4}}$\margincite{15.4}
\enddemo
\bigskip

\definition{\stag{15.5} Definition}  Let
\mr
\item  ${\Cal P}_1 = {\Cal P}_1(f^*) = \{Z \in \binom Xk:(\forall a,b
\in Z)(f^*_Z(a,b) = a\}$
\sn
\item  ${\Cal P}_2 = {\Cal P}_2(f^*) = \{Z \in \binom Xk:(\forall a,b
\in Z)(f^*_Z(a,b) = b\}$
\sn
\item  ${\Cal P}_\pm = {\Cal P}_\pm(f^*) = \binom Xk \backslash {\Cal
P}_1(f^*) \backslash {\Cal P}^*_2(f^*)$.
\endroster
\enddefinition
\bigskip

\proclaim{\stag{15.5a} Claim}  For $Y \in \binom Xk$ we have \nl
1) $Y \in {\Cal P}_\pm(f^*)$ iff $Y \in \binom Xk$ and $(\exists a,b
\in Y)(f^*_Y(a,b) = a \ne b)$ and $(\exists a,b \in Y)(f^*_Y(a,b) = b
\ne a)$.
\nl
2) If $Y \in {\Cal P}_\pm$, \ub{then} there 
are $a_Y \ne b_Y \in Y$ such that $f^*_Y(a_Y,b_Y) = b_Y$ and

$$
\{a,b\} \subseteq Y \wedge \{a,b\} \nsubseteq \{a_Y,b_Y\}
\Rightarrow f^*_Y(a,b) = a.
$$
\endproclaim
\bigskip

\demo{Proof}  By \scite{15.4}.
\enddemo
\bigskip

\proclaim{\stag{15.6} Claim}  1) $\langle {\Cal P}_1,{\Cal P}_2,{\Cal
P}_\pm \rangle$ is a partition of $\binom Xk$. \nl
2) For $Y \in {\Cal P}_\pm$ the pair $(a_Y,b_Y)$ is well defined (but
maybe $(b_Y,a_Y)$ can serve as well).
\endproclaim
\bigskip

\demo{Proof}  By \scite{15.4}.
\enddemo
\bigskip

\proclaim{\stag{15.7} Claim}  If ${\Cal P}_2(f^*) \ne \emptyset$ then
\mr
\item "{$\boxtimes_{f^*}(i)$}"  ${\Cal P}_2 = {\Cal P}_2(f^*)$
is a singleton, ${\Cal P}_\pm = \emptyset$ \ub{or}
\sn
\item "{${{}}(ii)$}"  $2k = |X|,{\Cal P}_2$ is a $\{Y^*,Y^{**}\}
\subseteq \binom Xk$ where $Y^* \cup Y^{**} = X$ and 
${\Cal P}_\pm = \emptyset$.
\endroster
\endproclaim
\bigskip

\demo{Proof}  Assume ${\Cal P}_2 \ne \emptyset$, let $Y^* \in {\Cal
P}_2$.  As $f^*$ is not a monarchy
\mr
\item "{$(*)_1$}"  ${\Cal P}_1 \cup {\Cal P}_\pm \ne \emptyset$.
\ermn
By Definition \scite{15.5} and Fact \scite{15.3}, $f^* \in {\Cal F}_{[r]}$
satisfies
\mr
\item "{$(*)_2(i)$}"  $f^*_{Y^*}(a,b) = b$ for $a,b \in Y^*$
\sn
\item "{${{}}(ii)$}"  if $g \in {\Cal F}_{[r]},g_{Y^*}(a,b) = b$ for $a,b \in Y^*$
then dom$_1(f^*) \supseteq \text{ dom}_1(g)$.
\ermn
Hence
\mr
\item "{$(*)_3$}"  if $Y_1 \in {\Cal P}_2,Y_1 \notin {\Cal P}_2,k^* =
|Y_1 \cap Y_2|$ and $Y \in \binom Yk,|Y \cap Y^*| = k^*$, \ub{then} $Y
\notin {\Cal P}_2$ (even $Y \in {\Cal P}_1 \leftrightarrow Y_2 \in
{\Cal P}_1$). \nl
[Why?  By $(*)_2$ as we 
can conjugate $f^*$ by $\pi \in \text{ Per}(X)$ which maps
$Y^*$ onto $Y^*$ and $Y_1$ to $Y$.]
\ermn
So by \scite{12.3} (applied to all $k^*$)
\mr
\widestnumber\item{$(*)_4(ii)$}
\item "{$(*)_4(i)$}"   ${\Cal P}_2$ is the singleton $\{Y^*\}$ \ub{or}
\sn
\item "{${{}}(ii)$}"   ${\Cal P}_2$ is a $\{Y^*,Y^{**}\},2k = |X|$ and
$Y^{**} = X \backslash Y^*$
\endroster
\sn
\roster
\item "{$(*)_5\,\,\,\,$}"   if $Z \in {\Cal P}_\pm$, \ub{then} $(\alpha)$ or
$(\beta)$
{\roster
\itemitem{ $(\alpha)$ }   $\{a_Z,b_Z\} = Z \cap Y^*,f^*_Z(b_Z,a_Z) = a_Z$
\sn
\itemitem{ $(\beta)$ }  $\{a_Z,b_Z\} = Z \backslash Y^*,f_Z(b_Z,a_Z) =
a_Z$.
\endroster}
[Why?  Otherwise as $k \ge 3$ we can choose $\pi \in \text{
Per}(X),\pi(Y^*) = Y^*,\pi(Z) = Z$ such that $\pi''\{a_Z,b_Z\}
\nsubseteq \{a_Z,b_Z\}$ and use \scite{15.2}, \scite{15.3} on a
conjugate of $f^*$.]
\ermn
It is enough by $(*)_4$ to prove ${\Cal P}_\pm = \emptyset$. \nl
So assume toward contradiction ${\Cal P}_\pm \ne \emptyset$. \nl
By $(*)_5$ one of the followng two cases occurs.
\enddemo
\bn
\ub{Case 1}:  $Z^* \in {\Cal P}_\pm,|Z^* \cap Y^*| = k-2$.

As we are allowed to assume $k+2 < |X|$ there is Y$ \in \binom Xk$
such that $|Y \cap Y^*| = k-1$ and $Y \cap Z^* = Y^* \cap Z^*$.  Now
(by $(*)_5$) we have $Y \notin {\Cal P}_\pm$ and (by $(*)_4$) we have $Y \notin {\Cal
P}_2$ so $Y \in {\Cal P}_1$.  So there is $\pi \in \text{ Per}(X)$
such that $\pi(Y^*) = Y,\pi \restriction Z^* =$ identity, let $f =
(f^*)^\pi$ so by \scite{15.3} we get a contradiction to the choice of
$f^*$.
\bn
\ub{Case 2}:  $Z^* \in {\Cal P}_\pm,|Z^* \cap Y^*| = 2$.

A proof similar to case 1 works if $Z^* \cup Y^* \ne X$.  So assume
$Z^* \cup Y^* = X$ (hence $2k-2 = |X|$).

Let $\pi \in \text{ Per}(X)$ be the identity on $Z^* \cap Y^*$ and
interchange $Z^*,Y^*$.  Apply \scite{15.3} on $f^*,(f^*)^\pi$ so
$(a_{Z^*},b_{Z^*}) \in \text{ dom}_1(f^*) \cap \text{
dom}_1((f^*)^\pi)$, etc. easy contradiction. \nl
${{}}$  \hfill$\square_{\scite{15.7}}$\margincite{15.7}
\bigskip

\remark{\stag{15.8} Remark}  if $\boxtimes(f^*)$ of \scite{15.7} holds
for some $f^*$ then ${\frak C}$ (in the context of \S5) ${\frak C}$ is full by
\scite{14.21}.
\endremark
\bigskip

\demo{\stag{15.9} Hypothesis}  For no $f \in {\Cal F}_{[r]}$ is
$\boxtimes(f)$.
\enddemo
\bigskip

\demo{\stag{15.10} Conclusion}  1) ${\Cal P}_2(f^*) = \emptyset$. \nl
2) ${\Cal P}_\pm \ne \emptyset$. \nl
3) ${\Cal P}_1 \ne \emptyset$. \nl
4) If $Y \in {\Cal P}_\pm$ and $|Y \cap Z_1| = |Y \cap Z_2|$ and $a_Y
\in Z_1 \leftrightarrow a_Y \in Z_2$ and $b_Y \in Z_1 \leftrightarrow
b_Y \in Z_2$ where, of course, $Y,Z_1,Z_1 \in
\binom Xk$, \ub{then} $Z_1 \in {\Cal P}_\pm \leftrightarrow Z_2 \in
{\Cal P}_\pm$.
\enddemo
\bigskip

\demo{Proof}  1) By \scite{15.9} + \scite{15.7}. \nl
2) Otherwise $f^*$ is a monarchy. \nl
3) Assume not, so ${\Cal P}_\pm = \binom Xk$.
Let $Y \in {\Cal P}_\pm,Z \in \binom Xk,Z \cap \{a_Y,b_Y\} =
\emptyset$ and \footnote{I am sure that after checking we can save}
$|Z \cap Y| > 2$ and $|Z \backslash Y| > 2$, we can get a
contradiction to $n(f^*)$-s minimality. \nl
4) By \scite{15.2} and \scite{15.3} as we can find $\pi \in \text{
Per}(X)$ such that $\pi(Y) = Y,\pi(Z_1) = Z_2,\pi(a_Y) = a_Y,\pi(b_Y)
= b_Y$.  \hfill$\square_{\scite{15.10}}$\margincite{15.10}
\enddemo
\bigskip

\proclaim{\stag{15.11} Claim}  If $Y,Z \in {\Cal P}_\pm$ and $Y \ne Z$,
\ub{then} there is no $\pi \in \text{ Per}(X)$ such that:

$$
\pi(Y) = Y,\, \pi(Z) = Z
$$

$$
\pi(a_Y) = a_Y, \, \pi(b_Y) = b_Y
$$

$$
\{\pi(a_Z),\pi(b_Z)\} \nsubseteq \{a_Z,b_Z\}.
$$
\endproclaim
\bigskip

\demo{Proof}  By \scite{15.2} + \scite{15.3}.
\enddemo
\bigskip

\proclaim{\stag{15.12} Claim}  If $Y \in {\Cal P}_\pm,Z \in {\Cal
P}_\pm,2 < |Y \cap Z| < k-2$ then $\{a_Z,b_Z\} = \{a_Y,b_Y\}$.
\endproclaim
\bigskip

\demo{Proof}  By \scite{15.11}.
\enddemo
\bigskip

\proclaim{\stag{15.12a} Claim}  If $Z_0,Z_1 \in {\Cal P}_\pm$ and $|Z_1
\backslash Z_0| = 1$ then $\{a_{Z_0},b_{Z_0}\} = \{a_{Z_1},b_{Z_1}\}$.
\endproclaim
\bigskip

\demo{Proof}  We shall choose by induction $i = 0,1,2,3,4$ a set $Z_i \in
{\Cal P}_\pm$ such that $j<i \Rightarrow |Z_i \backslash Z_j| = i-j$,
this holds for $i=1,j=0$, also this is sufficient as by \scite{15.12}
we have $i = j+3 \le 4 \Rightarrow \{a_{Z_i},b_{Z_i}\} =
\{a_{Z_j},b_{Z_j}\}$ as this applies to $(j,i) = (0,4)$ and $(j,i) =
(1,4)$ we get the desired conclusion by transitivity of equality.

To choose $Z_i$, let $x_i \in X \backslash Z_0 \cup \ldots \cup
Z_{i-1}$; possible as we exclude $k+i-1$ elements and choose $y_i \in
Z_0 \cap \ldots \cap Z_{i-1} \backslash \{a_{Z_{i-1}},b_{Z_{i-1}}\}$.
Now let $Z_i = Z_{i-1} \cup \{y_i\} \backslash \{x_i\}$ easily $j<i
\Rightarrow |Z_i \backslash Z_j| = i-j$ and $Z_i \in {\Cal P}_\pm$ by
\scite{15.10}(4) with $Y,Z_1,Z_2$ there standing for
$Z_{i-1},Z_{i-2},Z_i$ here. 
\enddemo
\bn
\ub{\stag{15.15} Choice}:  $Y^* \in {\Cal P}_\pm$.
\bigskip

\demo{\stag{15.16} Conclusion}
\mr
\item "{$(a)$}"  $Y^* \in {\Cal P}_\pm$
\sn
\item "{$(b)$}"  if $Y \in {\Cal P}_\pm$ then $(\{a_Y,b_Y\} = \{a_{Y^*},b_{Y^*}\}$ 
\sn
\item "{$(c)$}"  one of the following possibilities holds
{\roster
\itemitem{ $(\alpha)$ }  ${\Cal P}_\pm = \{Y^*\}$
\sn
\itemitem{ $(\beta)$ }  ${\Cal P}_\pm = \{Y \in \binom
Xk:\{a_{Y^*},b_{Y^*}\} \subseteq Y\}$
\sn
\itemitem{ $(\gamma)$ }  ${\Cal P}_\pm = \{Y^*,Y^{**}\}$ where $Y^{**}
= (X \backslash Y^*) \cup
\{a_{Y^*},b_{Y^*}\}$ and $|X| = 2k-2$.
\endroster}
\endroster
\enddemo
\bigskip

\demo{Proof of \scite{15.16}}  Note that
\mr
\item "{$(*)$}"  if $Y_1,Y_2,Y_3 \in {\Cal P}_\pm,|Y_1 \backslash Y_2|
= 1$ and $Y_3 \in {\Cal P}_\pm,Y_4 \in \binom Xk,|Y_3 \backslash Y_4|
= 1$ and $\{a_{Y_3},b_{Y_3}\} = \{a_{Y_1},b_{Y_1}\} \subseteq Y_4$
then $Y_4 \in P_\pm$ (hence $\{a_{Y_4},b_{Y_4}\} = \{a_{Y_3},b_{Y_3}\}
= \{a_{Y_1},b_{Y_1}\}$. \nl
[Why?  As there is a permutation $\pi$ of $X$ such that $\pi(a_{Y_1})
= a_{Y_2},\pi(b_{Y_1}) = Y_1,\pi(Y_3) = Y_1,\pi(Y_4) = Y_2$.  By
\sciteu{6.4} we get a contradiction to the choice of $f^*$.]
\ermn
By the choice of $Y^* \in {\Cal
P}_\pm$, we have (a), now (b) follows from (c) by \scite{15.12a} by
induction on $|Y \backslash Y^*|$,  so it is enough to prove (c).  Assume
$(\alpha), (\gamma)$ fail and we shall prove $(\beta)$.  So there is
$Z_1 \in {\Cal P}_\pm$ such that $Z_1 \notin \{Y^*,(X \backslash Y^*)
\cup \{a_{Y^*},b_{Y^*}\}\}$.  We can find $c_1,c_2 \in X \backslash
\{a_{Y^*},b_{Y^*}\}$ such that $c_1 \in Y^* \leftrightarrow c_2 \in
Y^*$ and $c_1 \in Z_1 \leftrightarrow c_2 \notin Z_1$. \nl
[Why?  if $Y^* \cup Z_1 \ne X$ any $c_1 \in X \backslash Y^*
\backslash Z_1,c_2 \in Z_1 \backslash Y^*$ will do; so assume $Y^*
\cup Z_1 = X$ so as $k+2 < |X|$ clearly
by \scite{15.12}, $|Z_1 \cap Y^*| \le 2$.  As not
case $(\gamma)$ of (c), that is by the choice of $Z_1$, 
necessarily $\{a_{Y^*},b_{Y^*}\} \nsubseteq Y^*
\cap Z_1$ and using $\pi \in \text{ Per}(X),\pi \restriction Z_1 =
\text{ id},
\pi(Y^*) = Y^*$ by \scite{15.11} we contradict \scite{15.2} + \scite{15.3}.]

Let $Z_2 = Z_1 \cup \{c_1,c_2\} \backslash (Z_1 \cap \{c_1,c_2\})$ so
$Z_1,Z_2 \in \binom Xk$ and 
$|Z_2 \cap Y^*| = |Z_1 \cap Y^*|$ and $Z_1 \cap \{a_{Y^*},b_{Y^*}\} =
Z_2 \cap \{a_{Y^*},b_{Y^*}\}$ hence by \scite{15.10}(4) we have $Z_2
\in {\Cal P}_\pm$ and clearly  $|Z_1 \backslash Z_2| = 1$.

By \scite{15.12a} we have $\{a_{Z_1},b_{Z_1}\} = \{a_{Z_2},b_{Z_2}\}$.
Similarly by $(*)$ we can prove by induction on 
$m = |Z \backslash Z_1|$ that $\{a_{Z_1},b_{Z_1}\}
\subseteq Z \in \binom Xk \Rightarrow Z \in {\Cal P}_\pm \and
\{a_{Z_1},b_2\} = \{a_{Z_1},b_{Z_1}\}$.  If
$(\beta)$ of (c) fails, then for some $Z_3 \in {\Cal
P}_\pm,\{a_{Z_1},b_{Z_1}\} \nsubseteq Z$.
Easily $\{a_{Z_3},b_{Z_3}\} \subseteq Z \in \binom Xk \Rightarrow Z
\in {\Cal P}_\pm \and \{a_2,b_2\} = \{a_{Z_3},b_{Z_3}\}$.  We 
can find $Y \in \binom Xk$, 
such that $\{a_{Z_1},b_{Z_1},a_{Z_3},b_{Z_3}\} \subseteq Y$;
contradiction.
${{}}$  \hfill$\square_{\scite{15.18}}$\margincite{15.18}
\enddemo
\bigskip

\proclaim{\stag{15.18} Claim}  1) The $(*)$ of \scite{14.4} holds. \nl
2) In \scite{15.16} clause (c), clause $(\alpha)$ holds.
\endproclaim
\bigskip

\demo{Proof}  1) Obvious by part (2) from $(\alpha)$. \nl
2) First assume $(\beta)$, so by \sciteu{6.17}, clause (b), 
\scite{15.2} + \scite{15.3} we have without loss of generality 
either $\{a,b\} = \{a_{Y^*},b_{Y^*}\} \subseteq Y \in \binom Xk
\Rightarrow f_Y(a,b) = b$ or $\{a_{Y^*},b_{Y^*}\}
\subseteq Y \in \binom XK \Rightarrow f_Y(a_{Y^*},b_{Y^*}) = b_{Y^*} =
f(b_{Y^*},a_{Y^*})$.  In both cases $f^*$ is simple and not a monarchy
contradiction to \scite{14.0a}. \nl
Second, assume $(\gamma)$.  Let $\langle \pi_i:i < i^* \rangle$ be a
list of the permutations $\pi$ of $X$ such that $\pi\{a_{Y^*}\} =
\{a_{Y^*},b_{Y^*}\}$.

Let $f^*_i$ be $f^*$ conjugated by $\pi_i$.  Now define $g^i$ for $i
\le i^*$ by induction on $i:g^0_Y(x_1,x_2) = x_1,
g^{i+1}_Y(x_1,x_2) = f^*_i(g^i_Y(x_1,x_2),x_2)$.  So $g^{i^*} \in
{\Cal F}_{[2]}$ and dom$_2(g^{i^*}) = \dbca_{i < i^*} \text{ dom}_2(f^*_i)$
where dom$_2(g) = \{(Z,a,b):a,b \in Z \in \binom Xk$ and $g_Z(a,b) = b
\ne a\}$, so dom$_1(g^{i^*}) = \dbcu_{i < i^*} \text{ dom}_1(f^*_i)$ hence
\mr
\item "{$(*)_1$}"  $g^{i^*}_Y(a_1,a_2) = a_2$ if $\{a_1,a_2\} =
\{a_{Y^*},b_{Y^*}\}$
\sn
\item "{$(*)_2$}"  $g^{i^*}_Y(a_1,a_2) = a_1$ if $\{a_1,a_2\} \ne
\{a_{Y^*},b_{Y^*}\}$.
\ermn
Now $g$ is simple but not a monarchial contradiction to \scite{14.0a}.
\hfill$\square_{\scite{15.18}}$\margincite{15.18}
\enddemo
\newpage

\head {\S7 The case $r \ge 4$} \endhead  \resetall \sectno=7
\bigskip

\demo{\stag{16.1} Hypothesis}  $r = r({\Cal F}) \ge 4$.
\enddemo
\bigskip

\proclaim{\stag{16.2} Claim}  1) For every $f \in {\Cal F}_r$ there is
$\ell(f) \in \{1,\dotsc,r\}$ such that
\mr
\item "{$\circledast$}"  if $Y \in \binom Xk,\bar a \in {}^r
Y,|\text{Rang}(\bar a)| < r$ (i.e. $\bar a$ is not one-to-one) \ub{then}
$f_Y(\bar a) = a_{\ell(*)}$.
\ermn
2) $r \le k$.
\endproclaim
\bigskip

\demo{Proof}  1) Clearly there is a two-place function $h$ from
$\{1,\dotsc,r\}$ to $\{1,\dotsc,r\}$ such that: if $y_1,\dotsc,y_r \in
Y \in \binom Xk,y_\ell = y_k \wedge \ell \ne k \Rightarrow f_Y(y_1,\dotsc,y_r) =
y_{h(\ell,k)}$; we have some freedom so let \wilog 
\mr
\item "{$\boxtimes$}"  $\ell \ne k \Rightarrow h(\ell,k) \ne k$.
\endroster
\enddemo
\bn
\ub{Case 1}:  For some $\bar x \in{}^r Y,Y \in \binom Xk$ and $\ell_1
\ne k_1 \in \{1,\dotsc,r\}$ we have 

$$
|\text{Rang}(\bar x)| = r-1
$$

$$
x_{\ell_1} = x_{k_1}
$$

$$
f_Y(\bar x) \ne x_{\ell_1},
$$
\mn
equivalently: $h(\ell_1,k_1) \notin \{\ell,k\}$. \nl
Without loss of generality $\ell_1 = r-1,k_1 = r,f_Y(\bar x) = x_1$
(as by a permutation $\sigma$ of $\{1,\dotsc,r\}$ we can replace $f$
by $f^\sigma:f^\sigma_Y(x_1,\dotsc,x_2) =
f_Y(x_{\sigma(1)},\dotsc,x_{\sigma(r)}))$.

We can choose $Y \in \binom Xk$ and $x \ne y$ in $Y$ so
$h(x,y,\dotsc,y) = x$ hence $\ell \ne k \in \{2,\dotsc,r\}
\Rightarrow h(\ell,k) = 1$. \nl
Now for $\ell \in \{2,\dotsc,r\}$ we have agreed $h(1,\ell) \ne \ell$,
(see $\boxtimes$), as $h \restriction \{(m,n):1 \le m < n \le r\}$ 
is not constantly 1 (as $f$ is not a monarchy)
for some such $\ell,h(1,\ell) \ne 1$ so \wilog \,  
$h(1,2) \ne 1,2$, so \wilog \, $h(1,2) = 3$ but as $r \ge 4$ we
have if $x \ne y \in Y \in \binom Xk$ 
then $f_Y(x,x,y,y,\dotsc,y)$ is $y$ as $h(1,2) = 3$ and is
$x$ as $h(3,4) = 1$, contradiction.  So
\mr
\item "{$\circledast$}"  $h \restriction \{(\ell,k):
1 \le \ell < k \le r\}$ is constantly 1.
\ermn

$$
\bar x \in {}^r X \text{ is with repetitions } \Rightarrow h(\bar
x) = x_1
$$
\mn
as required. 
\bn
\ub{Case 2}:  Not Case 1.
\nl
So $h(\ell,k) \in \{\ell,k\} \text{ for } \ell \ne k \in
\{1,\dotsc,r\}$.  Now let $Y \in \binom Xk,x \ne y \in Y$ and look at
$f_Y(x,x,y,y,\ldots)$ it is both $x$ as $h(1,2) \in \{1,2\}$ and $y$
as $h(3,4) \in \{3,4\}$, contradiction.
\nl
2) This follows as if $f \in {\Cal F}_{[r]}$ and $k < r({\Cal F})$
and $\ell(*)$ is as in part (1) then $f_Y(\bar x) = x_{\ell(*)}$
always, as $x_{\ell(*)}$ has repetitions by pigeon-hull.
\hfill$\square_{\scite{16.2}}$\margincite{16.2}
\bigskip

\definition{\stag{16.3} Definition}  $f = f_{r;\ell,k} = 
f^{r;\ell,k}$ is the $r$-place function 

$$
f_Y(\bar x) = \cases
x_\ell \qquad \bar x &\text{ is with repetitions} \\
  x_k \qquad &\text{ otherwise} \endcases
$$
\enddefinition
\bigskip

\proclaim{\stag{16.4} Claim}  1) If $f_{r;1,2} \in {\Cal F}$ \ub{then}
$f_{r;\ell,k} \in {\Cal F}$ for $\ell \ne k$. \nl
2) If $f_{r;1,2} \in {\Cal F},r \ge 3$ then $f_{r+1;1,2} \in {\Cal F}$.
\endproclaim
\bigskip

\demo{Proof}  1)  Trivial. \nl
2) Let $r \ge 5$ for $g(x_1,\dotsc,x_{r+1}) =
f_{r,1,2}(x_1,x_2,\tau_3(x_1,\ldots) \ldots \tau_r(x_1,\ldots))$ where
$\tau_m \equiv f_{r,1,m}(x_1,\ldots,x_m,x_{m+2},\dotsc,x_{r+1})$ that
is $x_{m+1}$ is omitted. \nl
Continue as in the proof of \scite{2.7}.  \hfill$\square_{\scite{16.4}}$\margincite{16.4}
\enddemo
\bigskip

\proclaim{\stag{16.5} Claim}  Assume $Y \in \binom Xk,\bar a \in {}^r Y$
is one-to-one.  There is $f = f_{Y,\bar a} \in {\Cal F}_r$ such that
$f_{Y,\bar a}(\bar a) = a_2$ and $f_{Z,\bar a}(\bar b) = b_1$ if 
$Z \in \binom Xk$ and $\bar b \in {}^r X$ 
is not one to one (so $(*)$ of \scite{14.4} holds).
\endproclaim
\bigskip

\demo{Proof}  Let $f \in {\Cal F}_r$ be not monarchical, and \wilog \,
$\ell(*)=1$ in \scite{16.2}.  By being not a monarchy, for some $Y,\bar a$
and some $k \in \{2,\dotsc,r\}$  we have $f(\bar a) = a_k \ne a_1$;
necessarily $\bar a$ is one-to-one.  Conjugating 
by $\pi \in \text{ Per}(X)$ we get any $Y,\bar a$ and permuting 
$[2,r]$, we get $f_{Y,\bar a}(\bar a) = a_2$.  \hfill$\square_{\scite{16.3}}$\margincite{16.3}
\enddemo
\newpage

     \shlhetal 

\newpage
    
REFERENCES.  
\bibliographystyle{lit-plain}
\bibliography{lista,listb,listx,listf,liste}

\enddocument

\enddocument